\documentclass[11pt]{article}
\usepackage{ae,lmodern} 
\usepackage[latin1]{inputenc}
\usepackage[T1]{fontenc}
\usepackage{amsmath}
\usepackage{amsthm}
 \usepackage{dsfont} 
 \usepackage{bm}
\usepackage{mathtools}
 \usepackage{amscd}
\usepackage{amssymb}
\usepackage{amsfonts}
\usepackage{hyperref}
\hypersetup{
	colorlinks   = true,
	urlcolor     = blue, 
	linkcolor    = blue, 
	citecolor   = red 
}
\usepackage{graphics}
\usepackage[dvipsnames]{xcolor}
\usepackage{mathtools}
 \usepackage{cancel}
\usepackage[noadjust]{cite}
\usepackage{enumitem}
\usepackage{soul} 
\usepackage[percent]{overpic}
\usepackage{grffile}

\setenumerate[2]{label=(\alph*)}
\setenumerate[1]{label=\textup{(\roman*)}}
\setenumerate[3]{label=(\arabic*)}

\newcommand{\normW}[1]{{rt\kern-0.25ex\vert\kern-0.25ex\vert #1 
		\vert\kern-0.25ex\vert\kern-0.25ex\vert}}
\newcommand{\normWbig}[1]{{\Big\vert\kern-0.25ex\Big\vert\kern-0.25ex\Big\vert #1 \Big\vert\kern-0.25ex\Big\vert\kern-0.25ex\Big\vert}}
\newcommand{\C}{\mathbb{C}}
\newcommand{\N}{\mathbb{N}}

\newcommand{\R}{\mathbb{R}}

\newcommand{\boA}{\mathcal{A}}
\newcommand{\boB}{\mathcal{B}}
\newcommand{\boC}{\mathcal{C}}
\newcommand{\boD}{\mathcal{D}}
\newcommand{\boE}{\mathcal{E}}
\newcommand{\boF}{\mathcal{F}}

\newcommand{\boI}{\mathcal{I}}

\newcommand{\boK}{\mathcal{K}}
\newcommand{\boL}{\mathcal{L}}
\newcommand{\boM}{\mathcal{M}}

\newcommand{\boS}{\mathcal{S}}
\newcommand{\boT}{\mathcal{T}}

\newcommand{\boV}{\mathcal{V}}
\newcommand{\boW}{\mathcal{W}}

\newcommand{\gc}{\mathfrak{c}}

\newcommand{\gs}{\mathfrak{s}}

\newcommand{\bu}{\boldsymbol{u}}
\newcommand{\boldeta}{\boldsymbol{\eta}}
\newcommand{\boldtheta}{\boldsymbol{\theta}}
\newcommand{\NV}{\mathcal{NE}(\mathbb{R})}
\DeclareFontFamily{U}{mathx}{}
\DeclareFontShape{U}{mathx}{m}{n}{<-> mathx10}{}
\DeclareSymbolFont{mathx}{U}{mathx}{m}{n}
\DeclareMathAccent{\widehat}{0}{mathx}{"70}
\DeclareMathAccent{\widecheck}{0}{mathx}{"71}
\newcommand{\wh}{\widehat }
\newcommand{\W}{\mathcal{W}}
\newcommand{\ve}{\varepsilon}
\newcommand{\vp}{\varphi}

\newcommand{\wc}{\widecheck }
\newcommand{\ptl}{{\partial}}

\newcommand{\Eodd}{\mathcal{E}_{\mbox{\tiny odd}}(\mathbb{R})}
\newcommand{\Heven}{H^2_{\mbox{\tiny e}}(\mathbb{R})}
\newcommand{\Leven}{L^2_{\mbox{\tiny e}}(\mathbb{R})}
 \newcommand{\TWclambda}{(S(\W_\lambda,c))}
\newcommand{\TWzlambda}{(S(\W_\lambda,0))}
 
\newcommand{\TWzdelta}{(S(\delta_0,0))}
\newcommand{\lambdazero}{\lambda_*}
\DeclareMathOperator{\supp}{{\rm supp}}

\newcommand{\loc}{\textup{\tiny{loc}}}
\providecommand{\abs}[1]{|#1 |}
\providecommand{\norm}[1]{\lVert#1 \rVert}

\renewcommand{\Re}{\operatorname{Re}}
\renewcommand{\Im}{\operatorname{Im}}
\newcommand{\sech}{\operatorname{sech}}

\newcommand{\wto}{\rightharpoonup}
\newcommand{\wstar}{\stackrel{\ast}{\rightharpoonup}}
\newcommand{\bq}{\begin{equation}}
	\newcommand{\eq}{\end{equation}}
\theoremstyle{plain}
\newtheorem{theorem}{Theorem}[section]
\newtheorem{proposition}[theorem]{Proposition}
\newtheorem{lemma}[theorem]{Lemma}
\newtheorem{corollary}[theorem]{Corollary}

\theoremstyle{definition}
\newtheorem{definition}[theorem]{Definition}

\newtheorem{remark}[theorem]{Remark}

\theoremstyle{remark}

\numberwithin{equation}{section}

\begin{document}
	\title{Gray and black solitons of nonlocal Gross--Pitaevskii equations: existence, monotonicity and nonlocal-to-local limit }                                 
	\author{
    Andr\'e de Laire\footnote{
    Univ.\ Lille, CNRS, Inria, UMR 8524 - Laboratoire Paul Painlev\'e, F-59000 Lille, France.\\	E-mail: {\tt andre.de-laire@univ-lille.fr}} 
    \and Salvador  L\'opez-Mart\'inez\footnote{
	Departamento de Matem\'aticas, Universidad Aut\'onoma de Madrid, Ciudad Universitaria de Cantoblanco, 28049, Madrid, Spain.
		E-mail: {\tt salvador.lopez@uam.es}}
        }
	\date{}
	\maketitle
\begin{abstract}
This article investigates the qualitative aspects of dark solitons of one-dimensional Gross--Pitaevskii equations with general nonlocal interactions, 
which correspond to traveling waves with subsonic speeds. Under general conditions on the potential interaction term,
we provide uniform bounds, demonstrate the existence of symmetric solitons, and identify conditions under which monotonicity is lost.
Additionally, we present new properties of black solitons. Moreover, we establish the nonlocal-to-local convergence, i.e.\ 
the convergence of the soliton of the nonlocal model toward the explicit dark solitons of the local Gross--Pitaevskii equation.
\end{abstract}   
	
	\maketitle
	
	\medskip
	\noindent{{\em Keywords:}
		  Gross--Pitaevskii equation, defocusing nonlocal Schr\"odinger equation, traveling wa\-ves, dark solitons, nonlocal  optical solitons,	nonzero boundary conditions at infinity
		
		\section{Introduction}\label{intro}
		\subsection{The problem}
In numerous physical models,  nonlinear effects manifest with spatially nonlocal interactions.  For instance, this is the case in nonlinear optics, where the correction to the refractive index at a specific spatial point depends on the light intensity around that point \cite{guo2015nonlocal,krolikowski2000}. 
Other examples include weakly interacting  Bose-Einstein condensates, 
as explained by Gross \cite{gross} and Pitaevskii \cite{pitaevskii}, 
and nematic liquid crystals \cite{PECCIANTI2012147}.  
Thus, we consider the following nonlocal Schr\"odinger equation, usually called the nonlocal Gross--Pitaevskii equation, that describes the dynamics in these physical settings \cite{koutsokostas2021universal}, and takes the (dimensionless) form
\begin{equation}
\label{GP-full}
i\partial_t\Psi+\partial_{xx}\Psi+g \Psi (\W* \abs{\Psi}^2)=0, 
\quad 	\text{in } \R\times \R.
\end{equation}
Here, $\Psi: \R\times \R\to \C$,   $\W$ is a real-valued distribution that represents the nonlocal response of the medium, $*$ is the convolution in space, and $g\in \{\pm 1\}$.
For instance,  $\W$ models the interaction between bosons in Bose-Einstein condensates.
We consider here the case $g=-1$, which corresponds to a defocusing nonlinearity, 
since we are interested in studying localized structures with a nontrivial background. 
Therefore, we suppose that  $\Psi$ satisfies the (normalized)  nonzero condition at infinity:
		\begin{equation}
			\label{nonzero}
			\lim_{\abs{x}\to \infty}\abs{ \Psi(x,\cdot)}=1,
		\end{equation}
representing the fact that the density is constant far away.

		To provide a clear mathematical context to the problem, it is useful to perform the change of variables  $\Psi\to e^{-it}\overline{\Psi} $, which leads to the equation 
		\begin{equation}
			\label{NGP}
			i\partial_{t}\Psi={\partial_{xx}\Psi}+\Psi(\W*(1-|\Psi|^{2})),\quad \text{ in }~\mathbb{R}\times\mathbb{R},
		\end{equation}
		where we assumed that $\W*1=1$.  We will also assume from now on that  $\W$ is an even tempered distribution. In this manner, \eqref{NGP} 
		is Hamiltonian and its energy
		\begin{equation*}
			E(\Psi(t))=\frac12 \int_{\R}\abs{\partial_x\Psi(t)}^2\,dx +\frac 14 \int_{\R}(\W*(1-\abs{\Psi(t)}^2))(1-\abs{\Psi(t)}^2)\,dx,
		\end{equation*}
		is formally conserved. The (renormalized) momentum and mass 
		\begin{equation*}
			p(\Psi(t))=\int_{\mathbb{R}}\langle i \partial_x \Psi(t),\Psi(t) \rangle\left(1-\frac{1}{|\Psi(t)|^2}\right)dx,\qquad 
   N(\Psi(t))=\int_\R (1-\abs{\Psi(t)}^2)dx, 
		\end{equation*}
		are formally conserved too whenever $\inf_{x\in\R}\abs{\Psi(x,t)}>0$, where 
  we used the notation $\langle z_1,z_2\rangle =\Re(z_1 \bar z_2)$, for $z_1,z_2\in \C$ (see \cite{de2010global}).

We will be interested in special solutions to \eqref{NGP} with boundary condition \eqref{nonzero}, and belonging to the energy space
and nonvanishing energy space:
  $$
  \boE(\R)=\{v \in H^{1}_{\loc}(\mathbb{R}) : 1-|v|^{2}\in L^2(\mathbb{R}), \ v' \in L^{2}(\mathbb{R})\}\quad \text{and}\quad 
     \NV=
		\{v \in \boE(\mathbb{R})  :  \inf_{\R}\abs{v}>0\},
  $$
  i.e.\  the so-called dark solitons. Roughly speaking, these are localized density notches that propagate without spreading \cite{kartashov07}. They have been observed for example in Bose--Einstein condensates \cite{denschlag2000,becker2008}. Dark solitons are typically classified as {\em black} or {\em gray}, depending on whether their moduli vanish at some point or remain away from zero. Mathematically, these two types of solitons require distinct approaches, as will be justified below.
  
  More precisely, dark solitons in our context will be nontrivial finite energy solutions to \eqref{NGP} of the form
		$$\Psi_c(x,t)=u(x-ct),$$
		which represents a traveling wave with profile $u:\R\to\C$ propagating at speed $c\in\R$. Hence, the soliton $u$ satisfies
		\begin{equation}
			\label{TWc}
			\tag{$S({\W,c})$}
			icu'+u''+u(\W* (1-|u|^2))=0,\quad\text{ in }\R.
		\end{equation}
		Notice that taking the complex conjugate of $u$ in equation \eqref{TWc}, we are reduced to the case $c\geq 0$. 
        By Lemma~\ref{lemma:regularity}, if $c\neq 0$, we can write $u=\sqrt{1-\eta}e^{i\theta}$, and we will use this notation throughout this paper.

In the case of a constant interaction given by a Dirac delta function $\W=\delta_0$, so that $\widehat \W=1$, the finite energy solutions to \eqref{TWc} are unique for $c\in[0,\sqrt2 )$,  up to translations and multiplications by complex constants of modulus one. These dark solitons are explicitly given by
	\begin{equation}
	\label{sol:1D}
	\bu(x;c)=\sqrt{\frac{2-c^2}{2}}\tanh\Bigg(\frac{\sqrt{2-c^2}}{2}x\Bigg)-i\frac{c}{\sqrt{2}}, \text{ so  }
 \boldeta(x;c)=\Big(1-\frac{c^2}{2}\Big)\sech^2\Bigg(\frac{\sqrt{2-c^2}}{2}x\Bigg),
\end{equation}
for $x\in \R$, where $\boldeta=1-\abs{\bu}^2.$ For brevity, we shall sometimes indicate the dependence on $c$ as $\bu_c$ and $\boldeta_c$ instead of $\bu(\cdot\,;c)$ and $\boldeta(\cdot\,; c)$, respectively. In view of \eqref{sol:1D}, the solitons for $c>0$ can be written in the polar form
\bq
\label{u:theta}
\bu(x;c)=\sqrt{1-{\boldeta(x;c)}}e^{i\boldtheta(x;c)}, \text{ with } \boldtheta(x;c)=\arctan
\Bigg(
\frac{\sqrt{2-c^2}}{c}
\tanh\Bigg( \frac{\sqrt{2-c^2}}{2} x
\Bigg)
\Bigg)-\frac{\pi}2.
\eq
 On the other hand, from \eqref{sol:1D}, we see, in particular, that $\boldeta$ is an even solution, strictly decreasing on $\R^+$. In addition, if $c\geq\sqrt{2}$, the only solutions are the constants of modulus one (see \cite{bethuel2008existence}).

The study of the existence of finite energy solutions for the nonlocal problem \eqref{TWc} was started in \cite{delaire-mennuni}, and improved  in \cite{dLMar2022} under the following general hypothesis, which we will assume to hold throughout this paper: 
\begin{enumerate}[label=(H\arabic*),ref=\textup{({H\arabic*})}]
	\setcounter{enumi}{-1}
	\item\label{H0} $\W\in\boM(\R)$ is  even, with $\wh \W(0)=1$,
 \end{enumerate}
 where $\boM(\R)$ denotes the space of real, signed and finite Borel measures (see Subsection~\ref{subsec:hyp} for more details). 
 The main existence result proved in \cite{dLMar2022}, and recast using the remarks in \cite{dLDuMar23}, is as follows.
\begin{theorem}[\cite{dLMar2022,dLDuMar23}]
\label{thm:dLMar}
  Assume that 
	\begin{equation}
		\label{def:sigma}
		\gs:=\smash{\inf_\R \big(     \widehat \W(\xi)+\xi^2/2\big)}\in (0,1].
	\end{equation}	
  Then, for almost every $c\in (0,\sqrt{2\gs})$,
  there exists a nontrivial solution $u\in \mathcal E(\R)$ to \eqref{TWc}.
\end{theorem}
From a physical point of view, Theorem~\ref{thm:dLMar} provides the existence of solutions for almost all $c$ below  the Landau speed  $c_L(\W)$, defined by
$$
c_L(\W){=}\inf_{\R}\frac{\omega (\xi)}{\abs{\xi}}, \quad \text{with }
\omega(\xi)=\Big(\xi^4+2\widehat \W(\xi) \xi^2\Big)^{1/2}.$$
Here $\omega$ is the dispersion relation associated with \eqref{NGP}. Thus, if $\gs$ given in \eqref{def:sigma} lies in $(0,1]$, 
 Theorem~\ref{thm:dLMar} provides the existence of solutions for almost all $c\in (0,c_L(\W))$. Notice that, since $\widehat \W(0)=1$, 
the Landau speed cannot exceed the sonic speed 
$$c_s(\W)=\smash{ \lim_{\xi\to 0^+} {\omega(\xi)}/{\xi} = \sqrt 2}.$$
To better understand the role of the nonlocal interactions in the behavior of solitons, 
we provided in \cite{dLDuMar23} numerical simulations for several kernels.
For instance,  consider the Gaussian function proposed in \cite{Lopez-Aguayo},
\bq
\label{pot:Gauss}
\W_\lambda(x)=\frac{1}{2\abs{\lambda}\sqrt{\pi}}e^{-\frac{x^2}{4{\lambda}^2}}, \ x\in \R, \quad \text{ so that }\quad  \wh\W_\lambda(\xi)=e^{-\lambda^2\xi^2},
 \ \xi\in \R,
 \eq
where $\lambda\in \R\setminus \{0\}$. Thus, for small $\lambda$, this potential represents a smooth approximation of the Dirac delta function.
By invoking Theorem~\ref{thm:dLMar}, for $\lambda \in[0,1/\sqrt{2})$,
we conclude the existence of solutions for almost every $c\in(0,\sqrt 2)$. 
Also, for $\lambda\geq 1/\sqrt 2$, we infer the existence of solutions for almost every 
$c\in (0,c_L(\lambda))$,  where the Landau speed is 
$c_L(\lambda)=\lambda^{-1}{\sqrt{1+\ln(2\lambda^2)}}.$
However, the numerical simulations in \cite{dLDuMar23} lead us to the conjecture of existence of solutions for speeds above $c_L(\lambda)$. On the other hand, the shapes of the solitons are very similar to \eqref{sol:1D}
if  $\lambda$ is small, but oscillations appear when $\lambda$ increases. To exemplify this behavior, we depict in Figure~\ref{fig:PotGaussienNew2-etaettheta} the profile $\eta$ and the phase $\theta$ of the numerical solutions to \eqref{TWc}
obtained in \cite{dLDuMar23} for the Gaussian potential \eqref{pot:Gauss}, for $\lambda=3$.
\begin{figure}[!ht]
	\centering
	\begin{tabular}{cc}
		\resizebox{0.5\textwidth}{!}{
			\includegraphics{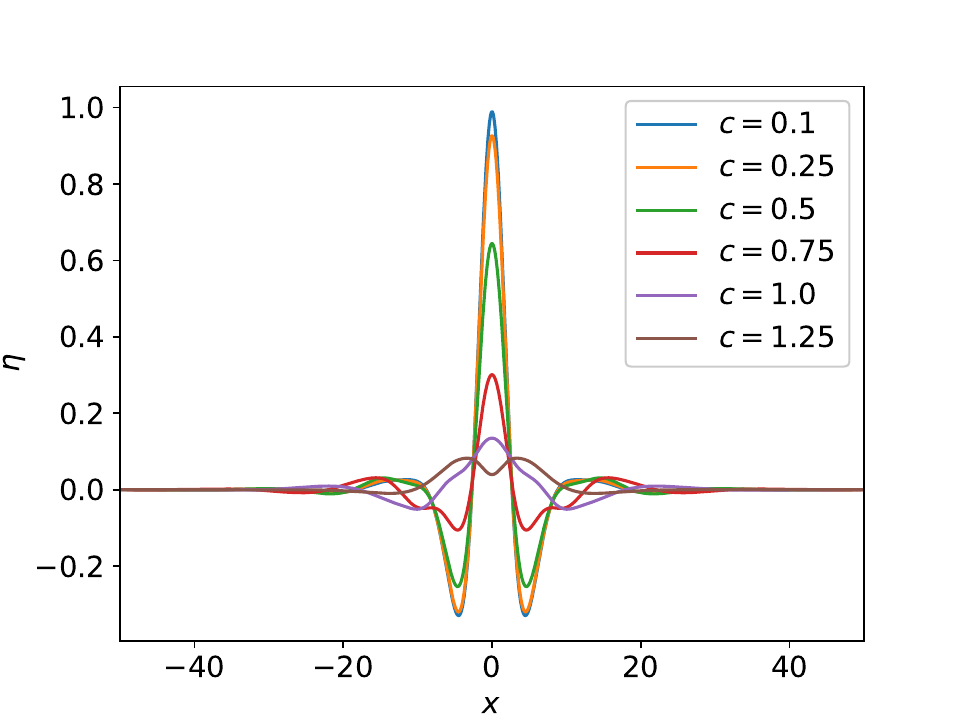}}   
			
	&                                                                            
	\resizebox{0.5\textwidth}{!}{\includegraphics{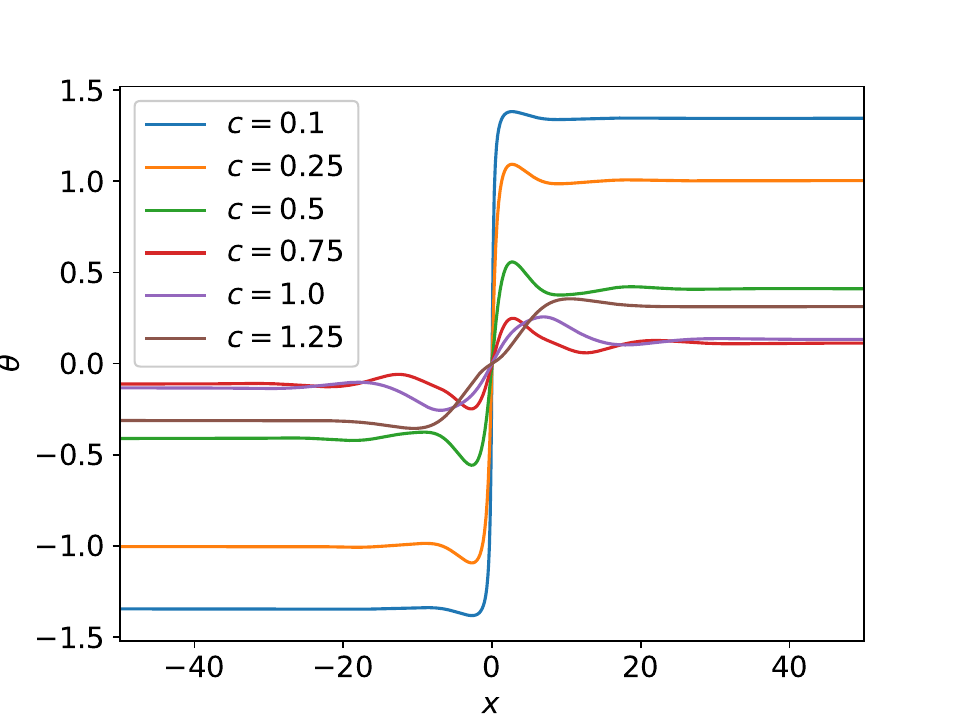}}
\end{tabular}
\caption{Numerically computed solitons for potential \eqref{pot:Gauss}
with $\lambda=3.0$, showing  $\eta$ (left panel) and ${\theta}_k$ (right panel) as function of $x$.
}
\label{fig:PotGaussienNew2-etaettheta}
\end{figure}

Another interesting potential  is
\begin{equation}
	 \label{pot:vander}
		{\mathcal W}_\lambda (x)
	=\frac{\beta}{\beta-2\lambda} \left(\delta_0 - \lambda {\rm e}^{-\beta |x|}\right), \ x\in\R,
	 \ \ \text {so that}\ \
	 	\widehat{\mathcal W}_\lambda(\xi)
	 = \frac{\beta}{\beta-2\lambda}\Big(1-\frac{2\lambda\beta}{\xi^2+\beta^2}\Big), \ \xi\in\R,
\end{equation}
for $\beta>0$  and $\lambda \in (-\infty,\beta/2)$,
for the study of dark solitons in a self-defocusing nonlocal Kerr-like medium.
If $\lambda>0$, the potential $\W_\lambda$ represents a strong repulsive interaction between particles that coincide in space, while the interaction becomes attractive otherwise, being this attraction more significant at short distances.
In contrast, for $\lambda<0$, the potential $\W_\lambda$ is purely repulsive.
As explained in \cite{dLDuMar23}, Theorem~\ref{thm:dLMar} provides the 
existence of solutions for a.e.\ $c\in (0,\sqrt{2})$, if
 $ \beta\geq\sqrt{2}$  and  $\lambda\in (-\infty,\beta/2)$,   or 
 $ \beta\in (0,\sqrt{2})$  and  $\lambda\in \Big(-\frac{\beta^3}{2(2-\beta^2)},\frac{\beta}{2}\Big)$.
We refer the reader to  \cite{dLDuMar23} for more detailed interpretations, where there are 
also many numerical simulations of solutions for several physically relevant potentials.
\begin{figure}[!ht]
	\centering
	\begin{tabular}[c]{cc}
		 \resizebox{0.5\textwidth}{!}{
   \includegraphics{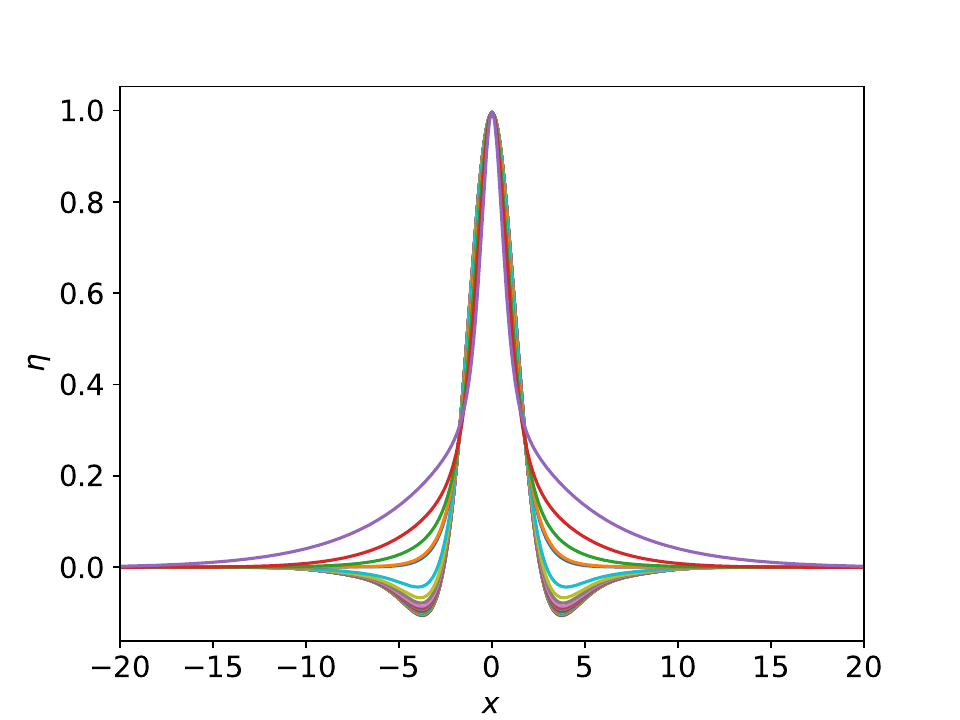}
 }
			&
    \resizebox{0.5\textwidth}{!}{
    \includegraphics{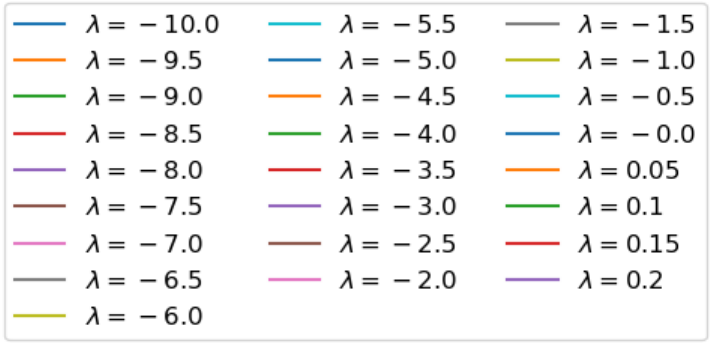}}
 \end{tabular}
\caption{Numerically computed solitons for potential \eqref{pot:vander}
  with $\beta=0.5$ and $c=0.1$ for several values of $\lambda$, plotted for $\eta=1-\abs{u}^2.$
  }
\label{fig:PotSalvadorNew-lambdabouge}
\end{figure}

The existence result provided by Theorem~\ref{thm:dLMar} does not handle the case of black solitons. Let us recall that black solitons are finite energy solutions vanishing at some point. In view of Lemma~\ref{lemma:regularity}, the only solutions having this property are the static ones, i.e.\ those with speed $c=0$, which can be taken as real-valued functions. In our previous works, we have avoided such solutions since they require an independent treatment.
For this reason, we introduce equation \eqref{TWc} for $c=0$ separately:
\begin{equation}\label{TW0}
	\tag{$S(\W,0)$} 
	u''+u(\W* (1-u^2))=0,\quad\text{ in }\R.
\end{equation}
Since we are assuming that  $\W*1=1$, this equation is equivalent to 
$-u''+u (\W*u^2)=u,$ which is sometimes called the (nonlocal) Ginzburg--Landau or Allen--Cahn equation. Recently, Lewin and  Nam \cite{lewin-thanh} have shown the existence of \textit{infinite ground state} solutions to this equation. However, these solutions have infinite energy.

Despite the significance of nonlocality in influencing the stability and dynamics of dark solitons in the physical literature, there are almost no mathematical results on the effect on the interaction potential, 
and only very particular cases of $\boW$ have been considered
\cite{nikolov2004,koutsokostas2021universal,guo2015nonlocal}.
 For instance, it is conjectured that the nonlocal effects can suppress the transverse instability of dark solitons \cite{Armaroli,dLGS-2,dLGS-3}.
A common approach, used for instance in \cite{koutsokostas2021universal,krolikowski2000}, is to study the so-called weakly nonlocal case, that is when we can replace $\boW$ by   $\boW_\ve(\cdot)=\boW(\cdot/\ve)/\ve$, for a small positive $\ve$. Then, performing a formal Taylor expansion, equation \eqref{NGP} reduces to the quasilinear equation 
		\begin{equation}
			\label{QGP}
   			i\ptl_t\Psi=\ptl_{xx}\Psi+\Psi(1-|\Psi|^2)+\lambda\Psi\ptl_{xx}(1-|\Psi|^2),\quad\text{ in }\R\times\R,
		\end{equation} 
  for some (small) real constant $\lambda$, depending on $\ve$.
This equation can be seen as a limit case of \eqref{NGP}, taking $\boV_{\lambda}=\delta_0+\lambda \partial_{xx}$, so that $\wh \boV_\lambda(\xi)=1-\lambda\xi^2$.
Of course, in this case, $\wh\boV_{\lambda}$ is not bounded, so does not satisfy \ref{H0}.

Only a few results deal with dark solitons for \eqref{QGP} when $\lambda\neq 0$. A branch of explicit dark solitons was found in \cite{krolikowski2000} for $\kappa\in [0,1/2)$. Also, in the setting of Korteweg models, Benzoni-Gavage, Danchin, Descombes and Jamin~\cite{benzoniStab} obtained existence of smooth dark solitons by using an ODE phase portrait analysis. Recently,  the study of localized solutions to \eqref{QGP} has been done in \cite{de2023exotic,lequiniou2024stability}, which leads to the existence of dark and antidark solitons, as well as more exotic localized solutions like dark cuspons, compactons, and composite waves, even for supersonic speeds.

This paper aims to provide rigorous results concerning the qualitative properties of dark solitons
under general conditions for $\boW$. In particular, to justify some behaviors detected numerically in \cite{dLDuMar23}.
More precisely:\\
(i) we obtain new a priori bounds for the $L^\infty$-norm of the solutions;\\
(ii) we show the {\em nonlocal-to-local} convergence, i.e.\ the convergence of a solution $u_{\lambda,c}$ to \eqref{TWc}
with a potential $\W_\lambda$, towards the explicit solution $\bu(\cdot\,; c)$ in \eqref{sol:1D}, as 
$\W_\lambda\to \delta_0$ in some sense;
\\
(iii) we prove the existence of symmetric solitons;\\
(iv) we provide conditions on $\W$ to characterize the loss of monotonicity of $\eta$ on  $\R^+$;\\
(v) we establish existence and new properties of \textit{black} solitons. 

As an application of (i), we can give explicit conditions on $\W$ to fill the speed gaps in Theorem~\ref{thm:dLMar}, so that we have existence of solutions \textit{for all} speed $c$ below the Landau speed, i.e.\ for 
all $c\in(0, \sqrt{2\gs})$. We also remark that the uniqueness, up to invariances, of nontrivial finite energy solutions to \eqref{TWc} for $c\in [0,\sqrt{2})$ is an open question. 


\subsection{Hypotheses and main results}
\label{subsec:hyp}
\subsubsection{Existence results}

Let us state the main hypotheses used in the following statements.

\begin{enumerate}[label=(H\arabic*),ref=\textup{({H\arabic*})}]
\setcounter{enumi}{0}
\item \label{W:derivative}
$\wh\W\in W^{1,\infty}(\R)$. In addition, either $\wh\W\in W^{2,\infty}(\R)$, or the map $\xi\mapsto \xi\big(\wh\W\big)'(\xi)$ is bounded and continuous a.e.\ in $\R$. 
\item\label{W:m} $\wh\W \in W_{\loc}^{1,\infty}(\R)$, $\widehat{ \mathcal  W}\geq 0$ in $\mathbb R$,    and there is $m\in [0,1)$ such that $\big(\wh\W\big)'(\xi)\geq-m\xi$ for a.e.\ $\xi>0$. 
\item\label{W:L1} $\W\in L^1(\R)$, $\W\geq 0$ in $\R$, and 
there is  $\tau>0$ such that $\sigma\coloneqq\inf_{x\in(-\tau,\tau)}\W(x)>0.$ 
\item\label{W:deltamu}  There is  $ \mu\in \boM(\R)$ such that $\W=A_\mu(\delta_0 + \mu),$
 with $\norm{\mu^-}_{\boM(\R)}<1$.
\end{enumerate}

These hypotheses cover a large class of physically relevant nonlocal potentials. In Subsection~\ref{subsec:examples} and 
in Section~\ref{sec:proofsexamples}, we will give several examples of potential satisfying the hypotheses introduced in this section.

Notice that for  $\W\in\boM(\R)$,  by the Jordan decomposition theorem, there are nonnegative measures $\W^+$ and $\W^-$ such that  $\W=\W^+-\W^-$, and $\W^+\perp\W^-$, i.e.\ they are mutually singular. Thus, the total variation is equal to
\begin{equation}
\label{def:M}
    \|\W\|_{\boM(\R)}=\int_\R d\W^+ + \int_\R d\W^-.
\end{equation}
 We will keep this notation for the rest of the paper. Moreover, the Fourier transform extends from 
$L^1(\R)$ to $\boM(\R)$ as 
\begin{equation}
\label{def:Fourier}
\wh \W(\xi)=\int_\R e^{-i x \xi}d\W(x).
\end{equation}
In particular, the Fourier transform is continuous everywhere, so the normalization condition $\wh  \W(0)=1$ in \ref{H0} is meaningful, 
and fixes the value of the constant $A_\mu$ in \ref{W:deltamu} as
\bq
\label{Au}
\frac{1}{A_\mu}=1+\wh{\mu^+}(0)-\wh{\mu^-}(0)=1+\norm{\mu^+}_{\boM(\R)}-\norm{\mu^-}_{\boM(\R)}.
\eq

In addition, \ref{H0} yields the following properties, which we use repeatedly, whose proofs can be found in Section 8.6 in \cite{folland}: 
 \begin{enumerate}
     \item $\wh\W\in L^\infty(\R)\cap\boC(\R)$ with 
     \bq
     \label{W:cota:medida}
     \|\wh\W\|_{L^\infty(\R)}\leq \|\W\|_{\boM(\R)}. 
     \eq
      \item For all $f\in L^p(\R)$, with $p\in [1,\infty]$, the convolution $\W*f$ is well-defined in 
     $L^p(\R)$ as 
     $(\W*f) (x)=\int_\R f(x-y) d\W(y) $, and the Young inequality holds:
          \begin{equation}
          \label{young}
              \| \W*f\|_{L^p(\R)}\leq\|\W\|_{\boM(\R)}\|f\|_{L^p(\R)}.
          \end{equation}
  \item By using that $\W$ is real-valued, even, and invoking Plancherel's identity, we have, for all $f,g\in L^2(\R)$,
 \begin{equation}
 \label{Plancherel}
 \int_\R f(\W*g)=\int_\R g(\W*f)=\frac{1}{2\pi}\int_\R \wh\W(\xi) \wh f(\xi)\overline{ \wh g}(\xi) d\xi. 
 \end{equation}
  \end{enumerate}

On account of Theorem~\ref{thm:dLMar}, the natural question is whether improving the \textit{almost all} in the speed interval is possible.
In \cite{dLDuMar23}, we have already established a general theorem that allows us to fill the gaps in the speeds, but it requires some a priori $L^\infty$-estimates on the solutions, which are not easy to obtain for nonlocal equations. For the sake of convenience, we recall the result now.

\begin{theorem}[Theorem 4.1 in \cite{dLMar2022}]
\label{thm:existence}
   Consider $\W$ satisfying \ref{W:derivative} and \ref{W:m}. 
  Assume that there exists a continuous function $ M: (0,\sqrt{2})\to (0,\infty)$ such that for any $u\in \NV$ solution to $(S(\boW,\tilde c))$, with  $\tilde c\in (0,\sqrt{2})$,  we have the estimate
$\norm{u}_{L^\infty(\R)}^2\leq  M(\tilde c)$.
If  $c\in (0,\sqrt{2})$ satisfies
\begin{equation}
\label{apriori}
m M(c)<1,
\end{equation}
where $m$ is given by \ref{W:m},
then there exists a nontrivial solution $u\in\NV$ to \eqref{TWc}.
\end{theorem}
This result is based on a refined study of the Palais--Smale sequences and the use of the profile decomposition theorem for bounded functions in $H^1(\R)$. Therefore, the problem of existence reduces to obtaining $L^\infty$-estimates for solutions to \eqref{TWc}, with good control of the constants, which is a difficult problem.
Also, we showed in \cite{dLDuMar23} that the necessary  $L^\infty$-estimate in 
Theorem~\ref{thm:existence} holds for a potential of the form given in \ref{W:deltamu} with 
\begin{equation}
    \label{dem:Mu}
M(c,\mu)=\Big(1+\frac{\|\mu^+\|_{\boM(\R)}}{1-\|\mu^-\|_{\boM(\R)}}\Big)\Big(1+\frac{c^2}{4}\Big).
\end{equation}
This yields the following result.
\begin{corollary}[\cite{dLMar2022}]
 \label{cor:existencegray}
 Let $c\in (0,\sqrt{2})$. Assume that $\W$ satisfy \ref{W:derivative}, \ref{W:m} and  \ref{W:deltamu}, with 
\begin{equation}\label{H:mmu}
    mM(c,\mu)<1,
\end{equation}
where $M(c,\mu)$ is given by \eqref{dem:Mu}. Then there is a nontrivial solution $u\in\NV$ to \eqref{TWc}.
\end{corollary}

Recall that we are assuming that \ref{H0} holds, so that the measure $\mu$ is even, with bounded Fourier transform. For instance, we can apply Corollary~\ref{cor:existencegray} to potential \eqref{pot:vander} for specific choices of $\beta$ and $\lambda$, as illustrated with Corollary~\ref{cor:example:vanderwaals} in subsection~\ref{subsec:examples}.

Our first result consists of establishing an $L^\infty$-bound for potentials satisfying \ref{W:L1}. In this manner, we can then apply Theorem~\ref{thm:existence} to obtain the existence of solutions to \eqref{TWc} in the whole subsonic interval, as follows.

\begin{theorem}\label{thm:existencegray}
 Let $c\in (0,\sqrt{2})$. Assume that $\W$ fulfills \ref{W:derivative}, \ref{W:m} and     \ref{W:L1},  such that $\tau$, $\sigma$ and $m$ satisfy
 \begin{equation}
 \label{H:msigmatau}
     \tau<\frac{\pi}{\sqrt{8+2c^2}},\quad m<\frac1{M(c,\tau,\sigma)},
      \text{ where }
 M(c,\tau,\sigma)=\frac{1}{\tau\sigma}\Big(1+\frac{c^2}{4}\Big)\sec\Bigg(\frac{\tau}{2}\sqrt{8+2c^2}\Bigg).
 \end{equation}
 Then there exists a nontrivial solution $u\in\NV$ to \eqref{TWc}.
\end{theorem}

We note that Theorem~\ref{thm:existencegray} implies that, for a fixed $\W$ satisfying \ref{W:derivative}, \ref{W:m} and \ref{W:L1} with $\tau<\pi/\sqrt{12}$ and $mM(\sqrt{2},\tau,\sigma)<1$, there exists a nontrivial solution $u\in\NV$ to \eqref{TWc} for every velocity $c\in (0,\sqrt{2})$. See  Subsection~\ref{subsec:examples} and 
in Section~\ref{sec:proofsexamples} for practical examples satisfying \eqref{H:msigmatau}.

So far we have excluded the existence of black solitons since the framework for their study is usually different. To be more specific, recall that, from Lemma~\ref{lemma:regularity}, gray solitons do not present vortices, so they can be lifted. This fact was the starting point of the work in \cite{dLMar2022}, since we could work in a simpler framework by considering the equation satisfied by the modulus of the solitons. In contrast, black solitons do present vortices (again by Lemma~\ref{lemma:regularity}), so the approach in \cite{dLMar2022} is no longer valid. Nevertheless, every black soliton is real-valued, modulo a change of phase, so we are allowed to work in the energy subspace of \emph{odd} real-valued functions, i.e.\ in the space
\[\Eodd=\{u\in \boE(\R): u\text{ is real-valued and odd}\}.\]
The minimization of the energy along $\Eodd$ gives the existence result. Adopting this approach, only the following condition on $\W$ is needed:
\begin{enumerate}[label=(H\arabic*),ref=\textup{({H\arabic*})}]
\setcounter{enumi}{4}
    \item\label{W:kappa} There exists $\kappa\geq 0$ such that $\widehat\W(\xi)\geq (1-\kappa\xi^2)^+$ for a.e. $\xi\in\R$.
\end{enumerate}
The precise statement reads as follows.

\begin{theorem}\label{thm:existenceblack}
	Let $\W$ satisfy \ref{W:kappa}. Then there exists a real-valued solution $u\in\Eodd$ to \eqref{TW0} which satisfies 
 \begin{equation}
 \label{oddmin}
  E(u)=\min\{E(v):\, v\in\Eodd\}.
  \end{equation}
\end{theorem}
In the sequel, we will call {\em odd minimizer} a function $u$ satisfying \eqref{oddmin} and the condition $u(\pm\infty)=\pm1$. In this manner, Theorem~\ref{thm:existenceblack} provides the existence of an odd minimizer black soliton solution.

\subsubsection{Nonlocal-to-local limit}

Once the existence of gray and black solitons is guaranteed, we will show that they satisfy the so-called \emph{nonlocal-to-local limit}.  The study of this kind of limit is a very natural question from a mathematical and physical point of view for several nonlocal problems \cite{ColomboEtAl,melchionna}. In our context, the idea is to prove that if the potential $\W$ is somehow close to the contact interaction $\delta_0$, then the solutions given by Corollary~\ref{cor:existencegray}, Theorem~\ref{thm:existencegray} and Theorem~\ref{thm:existenceblack} are also close in some sense to the (unique) solution to the local Gross--Pitaevskii equation. 

To make rigorous the idea of closeness to $\delta_0$, let us introduce the following hypotheses for a family of  potentials $\{\W_\lambda\}_{\lambda\in (0,\lambdazero)}\subset\boM(\R)$ satisfying \ref{H0}:

\begin{enumerate}[label=(H\arabic*),ref=\textup{({H\arabic*})}]
\setcounter{enumi}{5}
\item\label{Wlambda:aeconvergence} 
For some $\alpha\geq 0$, $\wh{\W_\lambda^+}(\xi)\to 1+\alpha,$ $\wh{\W_\lambda^-}(\xi)\to \alpha,$ as $\lambda\to 0$, for a.e. $\xi\in\R$.
\item\label{Wlambda:finitesup} 
The family is bounded in ${\boM(\R)}$, i.e.\ 
$\sup_{\lambda\in (0,\lambdazero)}\|\W_\lambda\|_{\boM(\R)}= M_{\lambdazero}$, for some $M_{\lambdazero}>0$. 
\end{enumerate}

Notice that assumption \ref{Wlambda:aeconvergence} implies that $\wh\W_\lambda(\xi)\to 1$ for a.e. $\xi\in\R$. This fact, combined with \ref{Wlambda:finitesup}, implies in turn that $\W_\lambda\to\delta_0$ vaguely
(see \eqref{conv:vague}).
We remark that we do not use the convergence of  $\W_\lambda$ to  $\delta_0$ in $\boM(\R)$, 
since it is, in general, too strong. Indeed, it does not hold for some classical potentials studied in Subsection~\ref{subsec:examples}, as showed in Lemma~\ref{lem:Wlam}.

The following result shows that the gray solitons obtained for the nonlocal equation \eqref{TWc} approximate the corresponding ones for the local equation. 

\begin{theorem}\label{thm:limittolocalgray}
Let $\lambdazero\in (0,\infty)$ and  $c\in(0,\sqrt{2})$.
Consider a family of even potentials $\{\W_\lambda\}_{\lambda\in (0,\lambdazero)}\subset\boM(\R)$ satisfying \ref{W:derivative}, \ref{W:m}, \ref{Wlambda:aeconvergence} and \ref{Wlambda:finitesup}, where the constant given by \ref{W:m} is denoted by $m_\lambda$ instead of $m$. Suppose also that one of the following two cases holds:
\begin{enumerate}
\item \label{alternativetau} For every $\lambda\in (0,\lambdazero)$, $\W_\lambda$
satisfies \ref{W:L1}, with constants $\tau=\tau_\lambda$ and $\sigma=\sigma_\lambda$, such that 
\begin{equation}
\inf_{\lambda\in (0,\lambdazero)}\tau_\lambda\sigma_\lambda>0  \quad  \text{ and }\quad 
\sup_{\lambda\in (0,\lambdazero)} m_\lambda M(c,\sigma_\lambda,\tau_\lambda)<1.
\end{equation}
\item \label{alternativemu} 
For every $\lambda\in (0,\lambdazero)$, $\W_\lambda$ is of the form $\W_\lambda=A_{\mu_\lambda}
(\delta_0+\mu_\lambda)$, satisfies \ref{W:deltamu} and 
\begin{equation}
\label{H:mmulambda}
 \sup_{\lambda\in (0,\lambdazero)} \frac{\|\mu_\lambda^+\|_{\boM(\R)}}{1-\|\mu_\lambda^-\|_{\boM(\R)}} <\infty
  \quad  \text{ and }\quad 
\sup_{\lambda\in(0,\lambdazero)} m_\lambda M(c,\mu_\lambda)<1.
\end{equation}
\end{enumerate}
Here, $M(c,\tau,\sigma)$ and $M(c,\mu)$ are defined in \eqref{H:msigmatau} and  in \eqref{dem:Mu}, respectively.
Let $u_\lambda \in \NV$ be the solution to \eqref{TWc} with $\W=\W_\lambda$, given by either Theorem~\ref{thm:existence} in case \ref{alternativetau}, or by Corollary~\ref{cor:existencegray} in case \ref{alternativemu}. Then, 
there exists a sequence $(x_\lambda,\theta_\lambda) \in\R^2$ such that 
\begin{equation}\label{eq:limittolocal}
u_\lambda(\cdot +x_\lambda)e^{i\theta_\lambda} \to \bu(\cdot \ ;c), \quad  \text{  as }\lambda\to 0, \text{ in } W^{k,\infty}_{\loc}(\R), 
\end{equation}
for every $k\geq 0$, where $\bu$ is the explicit soliton in \eqref{sol:1D}. Moreover, the convergence  in \eqref{eq:limittolocal} also holds weakly-$*$ in $L^\infty(\R)$, and weakly in
the homogeneous Sobolev space $\dot{H}^1(\R)$.
 \end{theorem}

\begin{remark} We think that the convergence stated in Theorem~\ref{thm:limittolocalgray} can be improved by refining the decay of solutions at infinity given in  \cite{dLMar2022}.
In a slightly different framework, a stronger convergence result will be given in Theorem~\ref{thm:even} below, for solutions constructed via the implicit function theorem.
\end{remark}

The proof of the previous result uses several arguments from \cite{dLMar2022}. 
In particular, the fact that the solutions are (limits of) mountain-pass critical points is essential to provide several key estimates.

We can deal with the nonlocal-to-local limit also for black solitons. In this case, we take advantage of the fact that the solutions given by Theorem~\ref{thm:existenceblack} are minimizers, as stated in the result.

\begin{theorem}\label{thm:limittolocalblack}
Let $\lambdazero\in (0,\infty)$.
Consider a family of even potentials $\{\W_\lambda\}_{\lambda\in (0,\lambdazero)}\subset\boM(\R)$ satisfying \ref{W:kappa}, \ref{Wlambda:aeconvergence} and \ref{Wlambda:finitesup}, where the constant $\kappa\geq 0$ is independent of $\lambda$. Let $u_\lambda\in \Eodd$ be the solution to \eqref{TW0}, with $\W=\W_\lambda$, given by Theorem~\ref{thm:existenceblack}. Then $u_\lambda\to \bu(\cdot\,;0)$ as $\lambda\to 0$ in $W^{k,\infty}_{\mbox{\tiny loc}}(\R)$ for every $k\geq 0$, where $\bu$ is the explicit soliton in \eqref{sol:1D}. Moreover, the convergence  in \eqref{eq:limittolocal} also holds weakly-$*$ in $L^\infty(\R)$, and weakly in the homogeneous Sobolev space $\dot{H}^1(\R)$.
\end{theorem}


\subsubsection{Existence of symmetric solitons for small \texorpdfstring{$\lambda$}{lambda}}

We have seen that the obtained black solitons are odd functions. By performing a perturbative argument, we can also establish the existence of gray solitons with \emph{even} modulus
 working in the Hilbert spaces of symmetric functions
\[\Heven=\{\eta\in H^2(\R): \eta \text{ is even}\} \quad \text{and } \quad\quad \Leven=\{\eta\in L^2(\R): \eta \text{ is even}\}.\]
This is certainly a natural property and was suggested by the numerical results in \cite{dLDuMar23}. Besides, a symmetry result for the solitons obtained by minimization at fixed momentum was proved in \cite{dLMar2022}, but this result applies only to some very particular potentials.

The proof will be based on a quantitative version of the implicit function theorem (see Theorem~\ref{thm:q:implict}). In this manner, we can construct solitons  bifurcating  from  $\bu(\cdot\,;c)$ in \eqref{sol:1D}. For applying this implicit function theorem, we only need the following general condition on a family $\{\W_\lambda\}_{\lambda\in (-\lambdazero,\lambdazero)}$ in $\boM(\R)$ for some $\lambdazero>0$:  

\begin{enumerate}[label=(H\arabic*),ref=\textup{({H\arabic*})}]
\setcounter{enumi}{7}
\item\label{Wlambda:C1} $\wh \W_0\equiv 1$ and, for a.e.\ $\xi\in\R$, the map $\lambda \in (-\lambdazero,\lambdazero)\mapsto\wh\W_\lambda(\xi)\in \R$, is of class $\boC^1$,
there exists $C>0$ such that,  for all $\lambda\in (-\lambdazero,\lambdazero)$, 
$\partial_\lambda\wh \W_\lambda(\cdot)$ is measurable in $\R$ and it satisfies
\begin{equation}
|\partial_\lambda\wh\W_\lambda(\xi)|\leq C(1+\xi^2),\quad\text{for a.e. }\xi\in\R.
\end{equation}
\end{enumerate}

Our result is as follows.
\begin{theorem}\label{thm:even}
Let us consider the open set of $H^2(\R)$
\[\Omega_e=\{\eta\in \Heven: \eta<1\,\,\text{ in }\R\}.\]
Assume that the family $\{\W_\lambda\}_{\lambda\in (-\lambdazero,\lambdazero)}$ in $\boM(\R)$ satisfies \ref{Wlambda:C1}, and fix $a\in (0,\sqrt 2)$. 
There exist $\lambda_a\in(0,\lambdazero)$ and an open set $\Omega_e^a\subset\Omega_e$ such that, for every $c\in\left[a,\sqrt{2-a^2}\right]$ and for every $\lambda\in (-\lambda_a,\lambda_a)$, there is $u_\lambda\in \NV$ such that $\eta_\lambda=1-|u_\lambda|^2\in\Omega_e^a$ and $u_\lambda$ is a solution  to \eqref{TWc} with $\W=\W_\lambda$. Moreover, the map
$\lambda \in  (-\lambda_a,\lambda_a) \mapsto \eta_\lambda \in \Omega_e^a$ is of class $\boC^1$, and for every $\lambda\in (-\lambda_a,\lambda_a)$, $u_\lambda$ is the unique solution to \eqref{TWc}, with $\W=\W_\lambda$, that belongs to $\Omega_e^a$. 
Furthermore, 
\bq
\label{TFI:conv}
u_\lambda\to \bu(\cdot\,; c)  \text{ in }\boC^{1}(\R)
\quad\text{ and }\quad \eta_\lambda \to \boldeta(\cdot\,; c) \text{ in } H^2(\R), \text{ as } \lambda\to 0.
\eq
Finally, if the family satisfies \ref{Wlambda:finitesup}, we have the stronger convergences
$u_\lambda\to \bu(\cdot\,; c)$ in $\boC^{k}(\R)$, and  $\eta_\lambda \to \boldeta(\cdot\,; c)$ in  $H^k(\R),$ 
for all $k\geq 0$,
as $\lambda\to 0$.
\end{theorem}

Let us remark that the classical implicit function theorem would require fixing first the speed $c$, to provide the existence of a number $\lambda_c>0$ and the existence of a family of solutions $\{u_\lambda\}_{\lambda\in (-\lambda_c,\lambda_c)}$. The difference in Theorem \ref{thm:even} is that given a (small) number $a$, we can obtain a family of solutions for all  $c\in \big[a,\sqrt{2-a^2}\big]$.

We stress that Theorem \ref{thm:even} establishes the existence of symmetric solutions to \eqref{TWc}, with $\W=\W_\lambda$,  satisfying the nonlocal-to-local limit with a strong convergence. However, notice that Theorem \ref{thm:even} does not provide solutions with speeds close to 0 or $\sqrt 2$, nor with large $\lambda$, contrary to Theorem~\ref{thm:existence}. Since the uniqueness is an open problem, it is not possible to guarantee whether the solutions in both results are the same or not.

Finally, let us recall that we are implicitly assuming that each $\W_\lambda$ satisfies \ref{H0}. Nevertheless, the implicit function argument  in Theorem~\ref{thm:even} can be done assuming only that $\wh\W_\lambda $ is even and that also satisfies $|\wh\W_\lambda(\xi)|\leq C(1+\xi^2)$, for $\lambda\in (-\lambdazero,\lambdazero)$. Therefore,  Theorem~\ref{thm:even}  could also be used to show the existence of symmetric dark soliton solutions to some quasilinear equations, like \eqref{QGP}, for small $\lambda$, constructed as a branch of solutions emanating from $\boldeta$.


\subsubsection{Monotonicity results}
Our last results concern the study of properties related to the monotonicity of solitons. As explained above, the numerical simulations in \cite{dLDuMar23} show that oscillations in the solitons appear when $\W$ ``is far''  from a Dirac delta function. Thus, we aim to show that, under some conditions on $\W$,  $\eta=1-|u|^2$ must have at least one oscillation on $\R^+.$
To motivate the conditions on the potential,  we perform a linearization of \eqref{eq:eta} around the trivial solution $\eta=0$, which yields: 

\begin{equation}\label{linearized}
	v''-2\W* v+c^2v=0,\quad\text{ in }\R.
\end{equation}
If $\W$ is given by an integrable function, we can seek solutions of the form $v(x)=e^{s x},$ with $s\in\R$, to 
get the characteristic equation
\begin{equation}
	\label{eq:car}
s^2-2\int_\R \W(y)e^{-s y}dy+c^2=0.
\end{equation}
Notice that, since $\W$ is even, we may take $s\geq 0$, without loss of generality. In \eqref{eq:car} we recognize  the  Laplace transform (also called bilateral or the two-sided Laplace transform), defined by 
\begin{equation}
	\widecheck\W(s)=\int_{\R}\W(y)e^{-s y}dy,
\end{equation}
whenever the integral is meaningful. Since this integral could diverge, we will consider it in a suitable generalized sense. To be precise, if $\W\in L^1(\R)$ and $\W\geq 0$ in $\R$, we define 
$\wc\W : [0,\infty) \to \R\cup\{\infty\}$ by
\begin{equation}
	\widecheck\W(s)=\lim_{R\to \infty}	\wc\W(s,R), \quad \text { where }\wc\W(s,R)= \int_{-R}^{\infty}\W(s)e^{-s y}ds.
\end{equation}
For instance, if $\W_\lambda(x)=\frac{1}{2\abs{\lambda}}e^{-{\abs{x}}/{\abs{\lambda}}}$, with $\lambda\neq 0$, its (generalized) Laplace transform is 
\begin{equation}\label{exponentialexample}
\wc\W_\lambda(s)=	\frac{1}{1-\lambda^2s^2},\quad\text{if } 
s\in [0,{1/\abs{\lambda} }),\quad \text{ and }\quad 
\wc\W_\lambda(s)=	\infty, \quad \text{if } s\in [{1}/\abs{\lambda},\infty).
	\end{equation}

Observe that evaluating the left-hand side of \eqref{eq:car} at $s=0$ leads to $-2+c^2$, which is negative in the subsonic regime $c\in [0,\sqrt{2})$. Thus, from \eqref{eq:car} we infer that the following condition 
\begin{equation}\label{H:conjecture}
s^2-2\wc \W(s)+c^2<0,\quad\text{for all } s\geq 0,
\end{equation}
 is sufficient to exclude the existence of exponential solutions to \eqref{linearized}. 
 If the function in \eqref{H:conjecture} is a polynomial (or a rational function), it 
  would be tempting to expect that it plays the role of the characteristic polynomial in  
linear second-order ODEs with constant coefficients, i.e.\  
    that the solutions to \eqref{linearized} present oscillations.
We will show that this can be proved rigorously for certain potentials such as the one in \eqref{exponentialexample}.
We refer to Theorem~\ref{thm:local:potentials} for more details. 
In the general case,  we will rigorously establish the presence of oscillations under a slightly more restrictive condition, as follows. 

\begin{theorem}\label{thm:nonmonotone} 
Let $c\in [0,\sqrt{2})$. Assume that  $\W$ satisfies \ref{W:L1} and  also the condition:
\begin{equation}\label{H:bumps}
s^2-2\wc \W(s)+c^2<-2cs, \quad \text{ for all }s\geq 0.
\end{equation}
Consider $u\in\boE(\R)$  a nontrivial solution to \eqref{TWc}. Then, for every $x_0\in\R$, neither of the following two conditions can hold: 
\begin{align}
\label{osc:contradiction}
 \eta \geq  0\  \text{ in } \R, \quad \text{ with } \quad \eta'\geq 0 \ \text{in }(-\infty,x_0),
 \\
\label{osc:contradiction2}
\eta \leq  0\  \text{ in } \R, \quad \text{ with } \quad \eta'\leq 0 \ \text{in }(-\infty,x_0),
\end{align}
 where $\eta=1-\abs{u}^2$. In particular, $\eta$ has an oscillation in the sense that there exist $x_1<x_2<x_3$ such that either 
$\eta'(x_1)>0$, $\eta'(x_2)<0$, $\eta'(x_3)>0$, or $\eta'(x_1)<0$, $\eta'(x_2)>0$, $\eta'(x_3)<0$. 
  \end{theorem}
Theorem~\ref{thm:nonmonotone}  establishes that $\eta$ has at least one oscillation. In addition, if, for instance, $\eta$ is positive, Theorem~\ref{thm:nonmonotone} implies that $\eta$ must have an infinite number of oscillations.

In the case $c=0$, the solutions are real-valued, by Lemma~\ref{lemma:regularity}. 
Thus, we can deduce that $u$ cannot be monotone.
\begin{corollary}
\label{cor:nonmonotone}
 Assume that  $\W$ satisfies \ref{W:L1} and also the condition in \eqref{H:conjecture} with $c=0$.
Consider $u\in\boE(\R)$  a nonconstant solution to \eqref{TW0} with $u(\pm \infty)=\pm 1$.
Then,  $u$ cannot be nondecreasing in $\R$.
\end{corollary}	

\begin{proof}
Assume by contradiction that  $u'\geq 0$ in $\R$, so that $ u(x)\in [-1,1]$, for all  $x\in \R$.  
By Lemma~\ref{eta:nonzero},  we deduce that $u(x)\in (-1,1)$ for all $x\in \R$, which gives  $\eta> 0$ in $\R$. Also, we deduce that there is $x_0<0$ such that $u(x)\in (-1,0)$, for $x<x_0$. Finally, since $\eta'=-2uu'$, it follows that  $\eta'\geq 0$ in 
$(-\infty,x_0)$.  This contradicts Theorem~\ref{thm:nonmonotone}.
\end{proof}

We conclude the analysis of the monotonicity by showing that, for a class of nontrivial potentials, increasing black solitons do exist. These potentials are of mixed attractive-repulsive nature, contrary to the purely repulsive ones satisfying \eqref{H:conjecture}. The result reads as follows.

\begin{theorem}\label{thm:monotonicity}
	Let $\W$ be of the form \ref{W:deltamu} with $\mu\in L^1(\R)$ even and nonpositive. Then every real-valued solution $u\in\boC^2(\R)$ to \eqref{TW0} such that $u(\pm\infty)=\pm 1$ satisfies
\begin{equation}\label{-1u1}
	-1<u<1,\quad\text{in }\R.
\end{equation}
Assuming in addition that $\mu$ is radially nondecreasing, there exists an odd minimizer solution to \eqref{TW0}, which is increasing on $\R$.
\end{theorem}

\subsection{Examples and applications}\label{subsec:examples}
We end this introduction by showing how the results stated above can be used to establish properties of dark solutions
in models considering long-range effects. We mainly consider physically relevant potentials studied in \cite{dLDuMar23}, so our results can be compared with their numerical simulations.

We start with the Gaussian potential introduced in \eqref{pot:Gauss}.
\begin{corollary}\label{cor:example:gaussian}
Consider $\W_\lambda$ the Gaussian potential  given in \eqref{pot:Gauss}, for some $\lambda>0$.  
\begin{enumerate}
 \item \label{item:gray} Let $c\in(0,\sqrt 2)$ and set 
\bq \label{gaussian:lambdac}
\lambda_c=\sup\Big\{\lambda\in\Big(0,\frac{\pi}{\sqrt{8+2c^2}}\Big) :\,\,
4e^{1/4}\sqrt{\pi}\Big(1+\frac{c^2}{4}\Big)\lambda^2\sec\big(\frac{\lambda}2\sqrt{8+2c^2}\big)<1\Big\}.
\eq
If $\lambda>\lambda_c$,
then there exists a nontrivial solution $u_\lambda\in\NV$ to $\TWclambda$.
Furthermore, there exist $\theta_0,x_0\in\R$ such that \eqref{eq:limittolocal} holds. 
\item \label{item:black} If $c=0$, there exists an odd minimizer solution $u_\lambda\in\Eodd$ to $\TWzlambda$ for every $\lambda\geq 0$. Furthermore, \eqref{eq:limittolocal} holds with $\theta_0=x_0=0$. 
\item \label{item:osc}
Let $c\in[0,\sqrt 2)$ and consider a nontrivial solution $u_\lambda\in\boE(\R)$ to $\TWclambda$, with $\lambda>\tilde \lambda_c$, 
where 
\bq
\label{gaussian:lambdactilde}
\tilde \lambda_c=\!\!\max_{s\geq \sqrt{2}-c} {s^{-1}}\big(\ln((s+c)^2/2)\big)^{1/2}.    
\eq
Then $\eta_\lambda=1-|u_\lambda|^2$ presents at least one oscillation (in the sense of Theorem~\ref{thm:monotonicity}). Moreover, if $c=0$ and $u_\lambda(\pm\infty)=\pm 1$, then $u_\lambda$ cannot be nondecreasing. 
\end{enumerate}
\end{corollary}

Equation \eqref{NGP} also appears in the study of the propagation of a polarized beam of light inside a planar cell filled with nematic liquid crystals, with the {\em nematic} potential  \cite{assanto2012nematicons}:
\bq
\label{pot:nematic}
\W_\lambda(x)=\frac{1}{2\abs{\lambda}}e^{-\frac{|x|}{\abs{\lambda} }},\ x\in \R, \quad \text{ so that }\quad \wh\W_\lambda(\xi)=\frac{1}{1+\lambda^2\xi^2}, \ \xi\in \R,
\eq
for $\lambda\neq 0$. In this context, dark solitons are called dark {\em nematicons}.

\begin{corollary}\label{cor:example:exponential}
Consider  $\W_\lambda$  given in \eqref{pot:nematic} for $\lambda>0$. For $c\in [0,\sqrt{2})$, we set 
\begin{align}
\label{exponential:lambdac}
\lambda_c&=\sup\Big\{\lambda\in\Big(0,\frac{\pi}{\sqrt{8+2c^2}}\Big):\,\, 4e\Big(1+\frac{c^2}{4}\Big)\lambda^2\sec\Big(\frac{\lambda}{2}\sqrt{8+2c^2}\Big)<1\Big\}, \\ 
\label{exponential:lambdactilde}
\tilde \lambda_c&=c^{-2}\, \sqrt{{4-c^2-\sqrt{8(2-c^2)}}}, \quad \text{ if } c\neq 0, \text{ and }\tilde\lambda_0=1/\sqrt{8}.   \end{align}
Then conclusions in \ref{item:gray}, \ref{item:black} and  \ref{item:osc}   of Corollary~\ref{cor:example:gaussian} hold.
\end{corollary}

\begin{remark}
To have an estimation of $\lambda_c$, it is sufficient to put $c=\sqrt 2$  in the inequality in its definition.
In the case of \eqref{gaussian:lambdac}, this leads to equation 
 $6\lambda^2 e^{1/4}\sec( \lambda \sqrt{3})=1$, whose first positive solution is  $\lambda_{\star}\approx 0.275$.
 Thus, for the Gaussian potential, we have the estimate $\lambda_c\geq \lambda_{\star},$ for all $c\in (0,\sqrt 2)$.
Similarly, for the nematic potential, we deduce that $\lambda_c$ in \eqref{exponential:lambdac} satisfies 
 $\lambda_c\geq 0.237,$ for all $c\in (0,\sqrt 2)$.
In addition, it is immediate to check that \eqref{gaussian:lambdactilde}, for $c=0$, yields $\tilde \lambda_0=1/\sqrt{2e}.$

Also, for the specific choice of $\W_\lambda$ in Corollary~\ref{cor:example:exponential}, it can be checked that condition $\lambda>\tilde \lambda_c$ is equivalent to \eqref{H:conjecture} (see Section \ref{sec:proofsexamples}). This shows that \eqref{H:bumps} is not sharp, in general.
\end{remark}

In the two examples above, the sign of $\lambda$ has no effect, so that we considered $\lambda>0$
for the sake of simplicity.
To also illustrate Theorem~\ref{thm:monotonicity}, we now consider the potential in \eqref{pot:vander}, where the sign of $\lambda$ plays an important role in the behavior of solitons.

\begin{corollary}\label{cor:example:vanderwaals}
Consider  $\W_\lambda$ given in \eqref{pot:vander}, with $\beta>0$
and $\lambda\in (-\infty,\beta/2)$. 
\begin{enumerate}
    \item Let $c\in (0,\sqrt 2)$ and set 
    \bq
    \label{def:lambdabeta}
    \lambda_{c,\beta}=\inf\{ \lambda<0 : \ \abs{\lambda}(2-\beta^2)<\beta^3/2, \ \text{ and } \ \abs{\lambda}<\beta^3/(4+c^2)\}.
    \eq
If $\lambda \in (\lambda_{c,\beta},\beta/2)$, then there exists a nontrivial solution $u_\lambda\in\NV$ to $\TWclambda$. Furthermore, there exist $\theta_0,x_0\in\R$ such that \eqref{eq:limittolocal} holds.
\item If $c=0$, there exists  an odd minimizer solution $u_\lambda\in\Eodd$ to $\TWzlambda$, for every $\lambda\in(-\infty,\beta/2)$. Also, \eqref{eq:limittolocal} holds with $\theta_\lambda=x_\lambda=0$. 
Furthermore, if $\lambda>0$, then $u_\lambda$ is increasing.

\item Let  $c\in [0,\sqrt 2)$ and $\lambda<0$. Assume that $u_\lambda\in\boE(\R)$  is a nontrivial solution to $\TWclambda$. If the number:
\begin{equation}\label{mon:waals}
\beta^2+2A_{\lambda}-c^2
-\big(
(\beta^2+2A_{\lambda}-c^2)^2
-4\beta^2(2A_{\lambda}-c^2)- 16\beta A_{\lambda} \lambda
\big)^{1/2}
\end{equation}
is either negative or has nonzero imaginary part,  
where $A_{\lambda}=\beta/(\beta-2\lambda)$,
then $\eta_\lambda=1-|u_\lambda|^2$ presents at least presents at least one oscillation. Moreover, if $c=0$ and $u_\lambda(\pm\infty)=\pm 1$, then $u_\lambda$ cannot be nondecreasing. 
\end{enumerate}
\end{corollary}
Notice that, for $c=0$, Corollary~\ref{cor:example:vanderwaals} rigorously captures the behavior obtained numerically in \cite{dLDuMar23} and shown in Figure~\ref{fig:PotSalvadorNew-lambdabouge}, i.e.\ that the black solitons $u_\lambda$ are monotone for $\lambda>0$. In addition, by using the number in \eqref{mon:waals}, with $\beta=1/2$ and $c=0$, it provides the existence of a strict negative minimum of $\eta_\lambda$, for all $\lambda<-1/4$.

Finally, for $\lambda\neq 0,$ let us consider  the rectangular potential, 
used to study   supersolids and nonlocal effects in optical media  \cite{aftalion-blanc,Kong-Wang-10}: 
\begin{equation}
	\label{pot:rect}
	\W_\lambda(x)=
		 \dfrac{1}{2\abs{\lambda}} \mathds{1}(x)_{[-\abs{\lambda},\abs{\lambda}]}
		\quad 	\text{ so that } \quad	\wh\W_\lambda(\xi)=\frac{\sin(\lambda\xi)}{\lambda\xi}, \ 
	\xi\in\R,
\end{equation}
  where $\mathds{1}$ is the indicator function. Let us also consider, for $\lambda\not=0$, the three Deltas potential, modeling a competition between repulsive and attractive interactions:
\begin{equation}
	\label{pot:deltas}
	\W_\lambda=2\delta_0-\frac{1}{2}(\delta_{-\lambda}+\delta_{\lambda}),
	\quad 	\text{ so that } \quad	\wh\W_\lambda(\xi)=2-\cos(\lambda \xi),\ \xi\in\R.
\end{equation}
   proposed in \cite{veskler2014} as a  model for interactions in a Bose--Einstein condensate. We showed in \cite{dLMar2022,dLDuMar23} that Theorem~\ref{thm:dLMar} applies to \eqref{pot:rect} and \eqref{pot:deltas}. However, these potentials do not satisfy the hypotheses in Theorem~\ref{thm:limittolocalgray} to complete the gaps in the speed interval. On the other hand, they fulfill \ref{Wlambda:C1}, so we can apply Theorem~\ref{thm:even} to deduce the existence of symmetric solutions as follows.
 
\begin{corollary}
\label{cor:ex:TFI}
Fix $a>0$ small. Let $\W_\lambda$ be given by either \eqref{pot:Gauss}, \eqref{pot:vander}, \eqref{pot:nematic}, 
 \eqref{pot:rect} or \eqref{pot:deltas}. There exists $\lambda_a>0$ such that for all $c\in [a,\sqrt{2-a^2}]$ and for all $\lambda\in (-\lambda_a,\lambda_a)$, there exists a nontrivial solution $u_\lambda\in\NV$ to $\TWclambda$.
 In addition, $\eta_\lambda=1-\abs{u_\lambda}^2$  is even and the convergences in \eqref{TFI:conv} hold,  as well as the stronger convergences for $\{u_\lambda\}$ in $\boC^{k}(\R)$,  and for $ \{ \eta_\lambda \} $ in $H^k(\R),$ for
 all $k\geq0$.
\end{corollary}

\subsubsection*{Organization of the paper}
In Section~\ref{sec:estimates} we prove the a priori bounds for the $L^\infty$-norm of the solutions to \eqref{TWc}. We establish 
the  nonlocal-to-local convergence in Section~\ref{sec:nonlocal}, and the existence of symmetric solitons by the implicit function argument in Section~\ref{sec:even}.
In Section~\ref{sec:monotonicity}, we study the loss of monotonicity of the solutions, while in Section~\ref{sec:blacksolitons} we analyze 
the existence and properties of black solitons. Finally, in Section~\ref{sec:proofsexamples}, we provide the proofs of the corollaries given in Subsection~\ref{subsec:examples}.

\section{A priori estimates}
\label{sec:estimates}

We begin by gathering properties for finite energy solutions to \eqref{TWc} in the following result, 
under our minimal assumption \ref{H0}.

\begin{lemma}\label{lemma:regularity}
Let $c\geq 0$ and  let $u\in \mathcal E(\R)$ be a solution to \eqref{TWc}. Then $u\in \boC^\infty(\R)\cap L^\infty(\R)$. Also, setting $\eta=1-\abs{u}^2$, we have $\eta \in H^{k}(\R)$ for all $k\in\N$, and 
\begin{equation}
	\label{limits-infty}
	\abs{u}(\pm \infty)=1,\  D^j u(\pm \infty)=D^j\eta(\pm\infty)=0, \ \ \text{for all } j\geq 1,
\end{equation}
where $D^j$ denotes the derivative of order $j$. In addition, $\eta$ satisfies the equation
\begin{equation}
\label{eq:etageneral}
 -\eta''+2{\mathcal W}* \eta - c^2\eta= 2|u'|^2 + 2({\mathcal W}* \eta)\eta,\quad\text{in }\R,
\end{equation}
which has the first integral
\begin{equation}\label{eq:firstintegral}
    c^2\eta^2+(\eta')^2=4|u'|^2(1-\eta),\quad\text{in }\R.
\end{equation}
In particular,
\begin{equation}
\label{eq:eta}
 -\eta''+2{\mathcal W}* \eta - c^2\eta= \frac{c^2\eta^2}{2(1-\eta)}
	+\frac{(\eta')^2}{2(1-\eta)} + 2({\mathcal W}* \eta)\eta,
\end{equation}
in the set $\{x\in\R : \eta(x)\neq 1\}$.  
Furthermore, depending on the value of $c$, we have:
\begin{enumerate}
    \item\label{item:realvalued} If $c=0$, then $u$ is real-valued up to a multiplication by a complex number of modulus one. In particular, $u$ presents a \emph{vortex}, i.e.\ there exists $x_0\in\R$ such that $u(x_0)=0$.
 Moreover, the unique continuation principle holds in the sense that if $\eta \equiv 0$ in some open interval, then $\eta\equiv 0$ in $\R$.
    \item If $c>0$, then $u$ belongs to $\NV$ and there exists a smooth lifting of $u$, i.e.\ there is a real-valued function $\theta\in\mathcal C^\infty(\R)$ such that $u=\abs{u}e^{i\theta}$ in $\R$. Moreover, $\theta'\in H^{k}(\R)$ for every $k\in\N$, and $D^j\theta(\pm \infty)=0$, for every $j\geq 1$.
\end{enumerate}
Finally, assuming in addition $\wh\W\geq 0$ on $\R$ and $\eta\in L^1(\R)$, it follows that either $\eta\equiv 0$ or
\begin{equation}\label{signofthemass}
    (2-c^2)\int_\R\eta(x)dx>0.
\end{equation}
\end{lemma}
\begin{proof}
The proof of Lemma~\ref{lemma:regularity} was performed in \cite{dLMar2022} except for item \ref{item:realvalued} and  \eqref{signofthemass}. 
Starting from this first point, let  $u\in\boE(\R)$ be a solution to \eqref{TW0}, and set $u=u_1 + iu_2$, for real-valued functions $u_1,u_2$.
Relation (2.6) in Proposition 2.2 of \cite{dLMar2022} establishes that $(\eta')^2=4|u'|^2|u|^2$, so that 
$|\langle u,u'\rangle|=|u'||u|$. Elementary algebra   yields 
$u'=\alpha(x) u$, whenever $|u|\not=0,$
where  $\alpha(x)={\langle u',u\rangle}/{|u|^2}.$ Hence, 
\bq
\label{dem:ODE:black}
u_1'=\alpha(x)u_1,\quad u_2'=\alpha(x)u_2,\quad\text{whenever }|u|\not=0.
\eq
Let $x_0\in\R$ be such that $u(x_0)\not=0$. Multiplying $u$ by a constant of modulus one, we may assume without loss of generality that $u_1(x_0)$ is real and that $u_2(x_0)=0$. Moreover, since $u$ is continuous, there exists $\delta>0$ such that $u(x)\not=0$ for every $x\in(x_0-\delta,x_0+\delta)$. Therefore, $\alpha$ is well-defined in this interval and by integrating  the second equation in \eqref{dem:ODE:black}, we obtain
\bq
\label{ODE:u2}
u_2(x)=u_2(x_0) e^{\int_{x_0}^x \alpha(y)dy}=0,\quad x\in (x_0-\delta,x_0+\delta).
\eq
In particular $u_2'(x_0)=0$.
On the other hand, $u_2$ satisfies the second order ODE:
\[-u''_2=q(x) u_2,\quad x\in\R,\]
where $q=\W*(1-|u|^2)\in L^\infty(\R)\cap \boC(\R)$. 
Since $u_2(0)=u_2'(0)=0$,  we can invoke the Cauchy--Lipshitz theorem to conclude that $u_2\equiv 0$ in $\R$, 
and therefore $u$ is a real-valued function.

It remains to establish \eqref{signofthemass}.  By integrating \eqref{eq:etageneral}, and using \eqref{limits-infty}, we get
		\[\int_\R (2\W*\eta-c^2\eta)=\int_\R (2|u'|^2+2\eta(\W*\eta)).\]
  Thus, if $\eta\not\equiv 0$, Plancherel's identity yields
  \[\int_\R (2\W*\eta-c^2\eta)=\int_\R (2|u'|^2+\frac{1}{\pi}\wh\W|\wh\eta|^2)>0.\]
  On the other hand, observe that
		\[\int_\R (\W*\eta)(x)dx=\int_\R d\W(x)\int_\R\eta(x)dx=\wh\W(0)\int_\R\eta(x)dx=\int_\R\eta(x)dx.\]
  Therefore, \eqref{signofthemass} holds.
  \end{proof}
We provide now properties for general smooth solutions to \eqref{TWc}.
\begin{lemma}\label{lemma:hamel}
 Let $c\geq 0$  and let $u\in \boC^2(\R)\cap L^\infty(\R)$ be a solution to \eqref{TWc}. Denote $M=\|u\|_{L^\infty(\R)}^2$. Then there exists $u_\infty \in \boC^2(\R)\cap L^\infty(\R)$ solution to \eqref{TWc} such that $v=|u_\infty|^2$ satisfies
\begin{align}
	-&v''\leq (2\W*(1-v)+{c^2}/{2})v,\quad\text{in }\R,\label{ineq:rhosquared}
	\\
	&M=\|v\|_{L^\infty(\R)}=v(0),\label{v0}
	\\
	&\big(\W*v\big)(0)\leq 1+{c^2}/{4}.\label{Wv0}
\end{align}
\end{lemma}

\begin{proof}
We start with a simple inequality satisfied by any bounded solution to \eqref{TWc}, in the spirit of \cite{bethuel}.
 Indeed, computing the derivatives of $|u|^2$ and applying Cauchy--Schwarz inequality, we obtain
	\begin{align*}
		(|u|^2)'' &= 2|u'|^2 + 2\langle u,u''\rangle =2|u'|^2 - 2\langle u, icu' + u\big(\W*(1-|u|^2)\big)\rangle
		\\
		&=2|u'|^2 - 2|u|^2\big(\W*(1-|u|^2)\big) - 2\langle cu,iu'\rangle
		\\
		&\geq -2|u|^2\big(\W*(1-|u|^2)\big) - {c^2}|u|^2/2.
	\end{align*}
This is precisely \eqref{ineq:rhosquared} for $v=|u|^2$.

Let $x_n\in\R$ be such that $|u(x_n)|^2\to M$ as $n\to\infty$. Let $u_n(x)=u(x+x_n)$. Clearly, $|u_n(0)|^2\to M$ as $n\to\infty$ and $|u_n(x)|^2\leq M$, for all $x\in\R$. In addition, since $u\in L^\infty(\R)$, arguing as in  Proposition 4.5 in \cite{dLMar2022}, it follows that $D^ku\in L^\infty(\R)$ for all $k\geq 1$. Hence, $\{u_n\}$ is bounded in $W^{k,\infty}(\R)$, for every $k\geq 0$. Thus, passing to a subsequence, there exists $u_\infty\in \boC^\infty(\R)$ such that $u_n\to u_\infty$ as $n\to\infty$ in $W^{k,\infty}_{\mbox{\tiny loc}}(\R)$, for all $k\geq 0$, which further satisfies $\|u_\infty\|_{L^\infty(\R)}^2=|u_\infty(0)|^2=M$. In other words, \eqref{v0} holds for $v=|u_\infty|^2$.  Moreover, since $u_n$ satisfies \eqref{TWc}, by using Lemma~4.19 in \cite{dLMar2022}, we can pass to the limit to deduce that $u_\infty$ also satisfies \eqref{TWc}. In particular,  \eqref{ineq:rhosquared} still holds for $v=|u_\infty|^2$.

Furthermore, since $M$ is a maximum, $v''(0)\leq 0$. Thus, evaluating \eqref{ineq:rhosquared} at zero yields
\[0\leq (2-2(\W*v)(0)+{c^2}/{2})M.\]
This implies directly \eqref{Wv0}.
\end{proof}
\begin{remark}
Using the notation in Lemma~\ref{lemma:hamel}, in the case of the explicit soliton $\bu(\cdot \, ;  c)$ in \eqref{sol:1D}, we see that $M=1$, and taking for instance $x_n=n$, we deduce that $u_\infty$   
is a constant function, so that  $v\equiv 1$. 
\end{remark}

From the previous lemma, we easily recover the estimates in Proposition~4.5 of \cite{dLMar2022}, but now being also valid for $c=0$ and for bounded solutions with infinite energy.

\begin{proposition}\label{prop:oldestimate}
Let $\W$ be of the form \ref{W:deltamu}, let $c\geq 0$ and let $u\in\boC^2(\R)\cap L^\infty(\R)$ be a solution to \eqref{TWc}. Then
\begin{equation}\label{est:deltamu}
	\|u\|^2_{L^\infty(\R)}\leq \left(1+\frac{\|\mu^+\|_{\boM(\R)}}{1-\|\mu^-\|_{\boM(\R)}}\right)\left(1+\frac{c^2}{4}\right).
\end{equation}
\end{proposition}

\begin{proof}
Thanks to Lemma~\ref{lemma:hamel}, we may assume without loss of generality that $v=|u|^2$ satisfies \eqref{v0} and \eqref{Wv0}. Thus, from \eqref{Wv0} we get
\[1+{c^2}/{4}\geq A_\mu(M+(\mu*v)(0))\geq A_\mu M(1-\|\mu^-\|_{\boM(\R)}).\]
Taking into account \eqref{Au}, this inequality leads directly to \eqref{est:deltamu}.
\end{proof}

The next result is new in the context of the nonlocal Gross--Pitaevskii equation. It provides uniform estimates for a class of potentials different from \ref{W:deltamu}, namely for those satisfying \ref{W:L1}. It is based on ideas in \cite{hamel-ryz}.

\begin{proposition}\label{prop:newestimate}
	Let $\W$ satisfy \ref{W:L1}, let $c\geq 0$ and let $u\in \boC^2(\R)\cap L^\infty(\R)$ be a solution to \eqref{TWc}. If $\tau<\pi/\sqrt{8+2c^2}$, then 
\begin{equation}\label{est:L1}
\|u\|^2_{L^\infty(\R)}\leq\frac{1}{\tau\sigma}\left(1+{c^2}/{4}\right)\sec\left(\tau\sqrt{2+{c^2}/{2}}\right),
\end{equation}
where $\tau,\sigma$ depend only on $\W$ and are given in \ref{W:L1}.
\end{proposition}

\begin{proof}
From \eqref{Wv0} we derive
\[1+\frac{c^2}{4}\geq \int_{-\tau}^0 \W(y)v(y)dy + \int_0^{\tau} \W(y)v(y)dy \geq \sigma\left(\int_{-\tau}^0 v(y)dy + \int_0^\tau v(y)dy\right).\]
The mean value theorem yields the existence of $y_-\in (-\tau,0)$ and $y_+\in(0,\tau)$ such that 
\[\frac{1}{\tau}\int_{-\tau}^0v(y)dy=v(y_-),\quad \frac{1}{\tau}\int_0^{\tau}v(y)dy=v(y_+).\]
Therefore,
\[v(y_{\pm})\leq \frac{1}{\tau\sigma}\Big(1+\frac{c^2}{4}\Big).\]
On the other hand, from \eqref{ineq:rhosquared} we deduce
	\[-v''-\left(2+{c^2}/{2}\right)v\leq 0,\quad\text{in }\R.\]
In sum, we have shown that $v$ is a subsolution to the  elliptic problem
\begin{equation}\label{linearprob}
	-\boL \phi = 0,\quad x\in (y_-,y_+),
	\quad 
	\displaystyle \phi(y_-)=\phi(y_+)=\frac{1}{\tau\sigma}\Big(1+\frac{c^2}{4}\Big),
\end{equation}
where
\[\boL \phi = \phi''+\big(2+{c^2}/{4}\big)\phi.\]

We will compare with a suitable supersolution later, but let us first check that $\boL$ satisfies the maximum principle in $(y_-,y_+)$, so that the comparison can be performed. Indeed, let us consider the eigenvalue problem
\begin{equation}\label{eigenvalueprob}
-\boL \phi = \lambda\phi,\quad x\in (-\tau,\tau),
	\quad 
	\displaystyle \phi(-\tau)=\phi(\tau)=0.
\end{equation}
It is immediate to verify that 
\[\lambda_1(\boL)=\frac{\pi^2}{4\tau^2}-2-\frac{c^2}{4},\quad \phi_1(\boL)=\cos\left(\frac{\pi}{2\tau}x\right),\]
where $\lambda_1(\boL)$ and $\phi_1(\boL)$ are the principal eigenvalue and eigenfunction, respectively, associated with \eqref{eigenvalueprob}. The choice $\tau<\pi/\sqrt{8+2c^2}$ assures that $\lambda_1(\boL)>0$, so that $\boL$ satisfies the maximum principle in $(-\tau,\tau)$. Actually, since $(y_-,y_+)\subset(-\tau,\tau)$, one has that the principal eigenvalue associated with the problem
\begin{equation*}
	-\boL \phi = \lambda\phi,\quad x\in (y_-,y_+),
\quad 
	\displaystyle \phi(y_-)=\phi(y_+)=0,
\end{equation*}
is also positive, that is, $\boL$ satisfies the maximum principle in $(y_-,y_+)$ too.

We are now in a position to compare $v$ with a supersolution to \eqref{linearprob}. To do so, let us consider the function
\[w(x)=\frac{1+{c^2}/{4}}{\tau\sigma\cos\left(\tau\sqrt{2+{c^2}/{2}}\right)}\cos\Big(x\sqrt{2+{c^2}/{2}}\Big).\]
It is straightforward to check that $w$ satisfies
\begin{equation*}
	-\boL w = 0,\quad x\in (-\tau,\tau),
	\quad
	\displaystyle w(-\tau)=w(\tau)=\frac{1}{\tau\sigma}\Big(1+\frac{c^2}{4}\Big).
\end{equation*}
Furthermore, since $w$ is even and radially decreasing for $x\in (-\tau,\tau)$, one has  
\[w(y_\pm)\geq\frac{1}{\tau\sigma}\Big(1+\frac{c^2}{4}\Big).\]
That is to say, $w$ is a supersolution to \eqref{linearprob}. Therefore, $v\leq w$ in $(y_-,y_+)$. Evaluating at zero leads exactly to \eqref{est:L1}.
\end{proof}

\begin{proof}[Proof of Theorem~\ref{thm:existencegray}]
It amounts to applying  Proposition~\ref{prop:newestimate} and Theorem~\ref{thm:existence}. 
\end{proof}	

\subsection{Decay and analyticity of finite energy solutions}
We briefly recall that the authors proved in \cite{dLMar2022} that several interesting properties can be established for finite energy solution to \eqref{TWc}, under general conditions for $\boW.$ In particular, these properties apply to the solution given by Theorem~\ref{thm:existencegray}.
The starting point is that equation \eqref{eq:etageneral} can be recast as:
\begin{equation}
	\label{eq:conv}
	\wh \eta(\xi)=L_c(\xi) \wh{F(\eta)}(\xi),\quad 
 \text{ so that }\eta=\boL_c *F(\eta), 
\end{equation}
$$
F(\eta)=\frac{c^2\eta^2}{2(1-\eta)}+\frac{(\eta')^2}{2(1-\eta)}+2\eta(\W\ast\eta),
\quad \text{ and }L_c(\xi)=\wh	\boL_c(\xi)=\frac{1}{\xi^2+ 2\wh\boW(\xi)-\xi^2}. 
$$
Indeed, as seen in the proof Theorem~2 in \cite{dLDuMar23},   under condition \eqref{def:sigma} with $\gs\in (0,1]$, for any $c\in (0,\sqrt{2\gs})$  there exists $\kappa\in (0,1/2)$ such that
$$
\wh \boW(\xi)\geq \tilde \gs-\kappa\xi^2, \quad \text{ for all }\xi\in\R, \text{ where }\tilde\gs= (c^2/2+\gs)/2.
$$
Therefore,
\begin{equation}
	\label{Lc:int}
\int_\R|L_c(\xi)|d\xi\leq\int_\R\frac{d\xi}{(1-2\kappa)\xi^2+2\tilde \gs-c^2}<\infty,
\end{equation}
and we conclude that $L_c\in L^1(\R)$, so that $\boL_c$ is a bounded continuous function on $\R$. 

Moreover, if the derivatives of $\boW$ are bounded, \eqref{eq:conv} allows deducing the algebraic decay and analyticity of the solutions, as follows.
\begin{lemma}[Corollary~5.9  and Theorem~5.12 in \cite{dLMar2022}]
\label{lem:decay}
Assume  that $\W$ satisfies \eqref{def:sigma} with $\gs\in (0,1]$ and that $\wh\W\in W^{k,\infty}(\R)$ for all $k\in \N$.
Let $c\in [0,\sqrt{2\gs})$. Then, $\boL_c\in L^\infty(\R)$ and $|\cdot|^\ell \boL_c\in L^2(\R)$ for all $\ell\in\N$. Moreover, for every $k, \ell \in\N$ and every $u \in \boE(\R)$ solution to \eqref{TWc}, it follows that 
	\[|\cdot|^\ell D^k\eta\in L^2(\R)\cap L^\infty(\mathbb{R}),\quad \lim_{x\to\pm\infty}|x|^\ell D^k \eta(x)=0, \]
 where $\eta=1-\abs{u}^2$.  In addition, if
 \begin{equation}
	\label{eq:analyticcond}
	\big|\big(\wh\W\big)'(\xi)\big|\leq C(|\xi|+1)\quad\text{a.e. }\xi\in\R,
\end{equation}
for some $C>0$,
 then $u$ (i.e.\  $\Re(u)$ and  $\Im(u)$) and $\eta$ are  real analytic on $\R$.	
\end{lemma}
As explained in \cite{dLMar2022}, \eqref{eq:conv} can also be used to establish the exponential decay of the solution.

\section{Nonlocal-to-local limit for gray solitons}
\label{sec:nonlocal}

In this section, we prove Theorem~\ref{thm:limittolocalgray} concerning solitons with a speed $c\in (0,\sqrt 2)$, so that they are gray soliton in the sense 
that their absolute value does not vanish. The proof will be a consequence of several lemmas. The first one provides sufficient conditions for a family of measures in $\boM(\R)$ to satisfy the property \eqref{Wlambda:tight}, which is reminiscent of the well-known concept of \emph{tight} family of measures.

\begin{lemma}\label{lemma:tight}
Let $\lambdazero\in(0,\infty)$. Consider a family of nonnegative measures $\{\boV_\lambda\}_{\lambda\in(0,\lambdazero)}\subset \boM(\R)$. Assume that there  is  $\alpha\geq 0$ such that 
\begin{align}
    &\wh\boV_\lambda(\xi)\to \alpha,\quad\text{for a.e. }\xi\in\R,\,\text{as }\lambda\to 0, \label{convFourier}
    \\
    &\|\boV_\lambda\|_{\boM(\R)}\to\alpha,\,\text{as }\lambda\to 0.  \label{convM}
\end{align}
For every $\varphi\in\boC^\infty_c(\R)$ with $\varphi\geq 0$ in $\R$, and for every compact set $K\subset\R$ such that $\supp(\varphi)\subset K$, one has
\begin{equation}\label{Wlambda:tight} \|\boV_\lambda*\varphi\|_{L^1(\R\setminus K)}\to 0, \text{ as }\lambda\to 0.
\end{equation}
In particular, if  $\boV\in L^1(\R)$ is such that  $\norm{\boV}_{L^1(\R)}=\alpha>0$ with  $\boV\geq 0$ in $\R$, then the family $\{\boV_\ve\}_{\ve>0}$ given by
\begin{equation}
    \label{mollifiers}
    \boV_\ve(x):={\ve^{-1}}\boV( x/ \ve  ), \quad \ve>0,
\end{equation}
 satisfies \eqref{Wlambda:tight}.
\end{lemma}

\begin{proof}
First, we use Proposition~{8.50} in \cite{folland} to deduce that \eqref{convFourier} implies that $\boV_\lambda$ converges vaguely to $\alpha\delta_0$, that is
\[\int_\R \psi d\boV_\lambda\to \alpha\psi(0),\quad\text{as }\lambda\to 0,\]
for every $\psi\in \boC_c^\infty(\R)$. In particular, fixing $x\in\R$ and taking $\psi=\varphi(x-\cdot)$ for some nonnegative $\varphi\in\boC_c^\infty(\R)$, we deduce
the pointwise convergence for the smooth function $\boV_\lambda*\varphi$, i.e.
\begin{equation}\label{pointwiseconv}
(\boV_\lambda*\varphi)(x)\to\alpha\varphi(x),\quad\text{as }\lambda\to 0.
\end{equation}

Let $K\subset\R$ be a compact set. Since $\boV_\lambda$ is nonnegative, we have
\bq
\label{decomposition}
\|{\boV_\lambda}*\varphi\|_{L^1(\R\setminus K)}=\|\boV_\lambda*\varphi\|_{L^1(\R)} - \|\boV_\lambda*\varphi\|_{L^1(K)}.
\eq
On the one hand, the properties of the convolution and \eqref{convM} yield
\[\|\boV_\lambda*\varphi\|_{L^1(\R)}=\|\boV_\lambda\|_{\boM(\R)}\|\varphi\|_{L^1(\R)}\to\alpha\|\varphi\|_{L^1(\R)},\quad\text{as }\lambda\to 0.\]
On the other hand, we aim to apply the dominated convergence theorem to
\begin{equation}
\label{L1K}
\|\boV_\lambda*\varphi\|_{L^1(K)}=\int_K (\boV_\lambda*\varphi)(x)dx.
\end{equation}
Since \eqref{pointwiseconv} holds, we only need to check the domination. Indeed, by Young's inequality \eqref{young},
\[(\boV_\lambda*\varphi)(x)\leq \|\varphi\|_{L^\infty(\R)}\sup_{\lambda\in(0,\lambdazero)}\|\boV_\lambda\|_{\boM(\R)}:=h,\quad\text{for all }x\in K,\]
where the constant  $h$ is trivially in $L^1(K)$. Thus, we conclude that $\boV_\lambda*\varphi$
converges to $\alpha \varphi$ in $L^1(K)$. Therefore, using \eqref{decomposition}, we infer 
\[\|\boV_\lambda*\varphi\|_{L^1(\R\setminus K)}\to \alpha\|\varphi\|_{L^1(\R)}-\alpha\|\varphi\|_{L^1(K)}=\alpha\|\varphi\|_{L^1(\R\setminus K)},\quad\text{as }\lambda\to 0.\]
If we choose $K$ such that $\supp(\varphi)\subset K$, then $\|\varphi\|_{L^1(\R\setminus K)}=0$, so
\[\|\boV_\lambda*\varphi\|_{L^1(\R\setminus K)}\to 0,\quad\text{as }\lambda\to 0,\]
which proves that \eqref{Wlambda:tight} holds.

Finally, considering  the family defined in \eqref{mollifiers}, 
we have $\norm{\boV_\ve}_{\boM(\R)}=\norm{\boV_\ve}_{L^1(\R)}=\alpha$, and $\wh\boV_\ve(\xi)=\wh\boV(\ve\xi)\to \norm{\boV}_{L^1(\R)}=\alpha$, as $\ve\to 0$, for all $\xi\in\R$. 
Therefore, \eqref{convFourier} and \eqref{convM} hold, and  $\boV_\ve$ satisfies the convergence in \eqref{Wlambda:tight}.
\end{proof}

The next result enables us to pass to the limit, as $\lambda\to 0$, in the nonlocal term of the weak formulation of \eqref{TWc}.

\begin{lemma}\label{lemma:limitnonlocalterm}
 Consider a family of  potentials $\{\W_\lambda\}_{\lambda\in (0,\lambdazero)}\subset\boM(\R)$ satisfying \ref{Wlambda:aeconvergence} and \ref{Wlambda:finitesup}, and a sequence $\{u_\lambda\}_{\lambda\in (0,\lambdazero)}$, with  $\sup_{\lambda\in(0,\lambdazero)}\|u_\lambda\|_{L^\infty(\R)}<\infty$, 
such that 
    \[u_\lambda\to u,\quad\text{as }\lambda\to 0,\,\,\text{ in }L^\infty_\loc(\R),\]
for some $u\in L^\infty(\R)$. Then, $u_\lambda (\W_\lambda*\eta_\lambda)$ converges to $u \eta$ in $\boD'(\R)$, where $\eta_\lambda=1-|u_\lambda|^2$ and $\eta=1-|u|^2$.    
\end{lemma}
    
\begin{proof}
   Let us fix $\varphi\in\boC^\infty_c(\R)$. We need to prove that 
   $$\int_\R\varphi u_\lambda\big(\W_\lambda*\eta_\lambda\big)\to\int_\R\varphi u \eta,\quad\text{as }\lambda\to 0.$$
   For this purpose, we write
    \begin{align*}
        \int_\R\varphi u_\lambda\big(\W_\lambda*\eta_\lambda\big)&=\int_\R\eta_\lambda\big(\W_\lambda*(\varphi u_\lambda)\big)=\int_K \eta_\lambda\big(\W_\lambda*(\varphi u_\lambda)\big) + \int_{\R\setminus K}\eta_\lambda\big(\W_\lambda*(\varphi u_\lambda)\big)
        \\  
        &=\int_K\big(\eta_\lambda-\eta\big)\big(\W_\lambda*(\varphi u_\lambda)\big) + \int_K\eta\big(\W_\lambda*(\varphi u_\lambda)\big) + \int_{\R\setminus K}\eta_\lambda\big(\W_\lambda*(\varphi u_\lambda)\big),
    \end{align*}
    where $K$ is a compact set such that $\supp(\varphi)\subset K$. We aim to analyze these last three integral terms. For the first one, we use the uniform convergence in compact sets of $\eta_\lambda$, the uniform bound on $u_\lambda$, and condition \ref{Wlambda:finitesup}, to derive
    \[\left|\int_K\big(\eta_\lambda-\eta\big)\big(\W_\lambda*(\varphi u_\lambda)\big)\right|\leq \|\eta_\lambda-\eta\|_{L^\infty(K)}\|\varphi\|_{L^\infty(\R)}\sup_{\lambda\in (0,\lambdazero)}\left(\|u_\lambda\|_{L^\infty(\R)}\|\W_\lambda\|_{\boM(\R)}\right)|K|\to 0,\]
as $\lambda\to 0$.  On the other hand, concerning the third term, we split it as
\[\int_{\R\setminus K}\eta_\lambda\big(\W_\lambda*(\varphi u_\lambda)\big)=\int_{\R\setminus K}\eta_\lambda\big(\W^+_\lambda*(\varphi u_\lambda)\big)-\int_{\R\setminus K}\eta_\lambda\big(\W^-_\lambda*(\varphi u_\lambda)\big).\]
Now we apply Lemma~\ref{lemma:tight} with $\boV_\lambda=\W_\lambda^+$, and again the uniform estimate on $u_\lambda$, to obtain
    \[\left|\int_{\R\setminus K}\eta_\lambda\big(\boV_\lambda*(\varphi u_\lambda)\big)\right|\leq \|\eta_\lambda\|_{L^\infty(\R)}\|u_\lambda\|_{L^\infty(\R)}\|\boV_\lambda*\abs{\varphi}\|_{L^1(\R\setminus K)}\to 0,\]
as $\lambda\to 0$. The same can be argued for $\boV_\lambda=\W_\lambda^-$. Therefore, the third term vanishes as $\lambda\to 0$.

It remains to consider the second term. Let us denote $f=\eta\chi_{K}$, $g_\lambda=\varphi u_\lambda$ and $g=\varphi u$. We claim that
\bq
\label{lim:claim}
\int_\R f\big(\W_\lambda*g_\lambda\big)\to\int_\R f g,\quad\text{as }\lambda\to 0,
\eq
    which will conclude the proof of the lemma. To prove the claim, notice first that $f\in L^2(\R)$ since it is a bounded function with compact support. Moreover, $g_\lambda\to g$ in $L^2(\R)$ since
    \[\|g_\lambda-g\|_{L^2(\R)}\leq \|u_\lambda-u\|_{L^\infty(\supp(\varphi))}\|\varphi\|_{L^2(\R)}.\]
    Now, since the Fourier transform is a continuous operator from $L^2(\R)$ to $L^2(\R)$, we also deduce that $\wh g_\lambda\to\wh g$ in $L^2(\R)$. Thus, by the partial converse of dominated convergence theorem, there exist a sequence $\{\lambda_n\}\subset (0,\lambdazero)$ converging to zero, and a function $h\in L^2(\R)$, such that 
    \[\wh g_{\lambda_n}(\xi)\to \wh g(\xi),\     \text{ and }\     \abs{\wh g_{\lambda_n}(\xi)}\leq h(\xi),\quad  \text{for a.e. }\xi\in\R.\]
    This implies, using \ref{Wlambda:aeconvergence}, that 
    \[\wh f(\xi)\wh\W_{\lambda_n}(\xi)\wh g_{\lambda_n}(\xi)\to \wh f(\xi)\wh{g}(\xi),\quad\text{for a.e. }\xi\in\R,\]
    and, furthermore, 
    \[|\wh f(\xi)\wh\W_{\lambda_n}(\xi)\wh g_{\lambda_n}(\xi)|\leq M|\wh f(\xi)||h(\xi)|,\quad\text{for a.e. }\xi\in\R,\,\,\text{for all }n,\]
    where $M=\sup_{\lambda\in(0,\lambdazero)}\| \wh\W_\lambda\|_{L^\infty(\R)}$ is finite by \ref{Wlambda:finitesup}. Since $\wh f h\in L^1(\R)$, the dominated convergence theorem  and \eqref{Plancherel} yield
     \[\int_\R f\big(\W_{\lambda_n}*g_{\lambda_n}\big)=\frac{1}{2\pi}\int_\R \wh f(\xi)\wh\W_{\lambda_n}(\xi)\overline{\wh g}_{\lambda_n}(\xi)d\xi\to \frac1{2\pi}\int_\R \wh f(\xi)\overline{\wh{g}}(\xi)d\xi= \int_\R f g.\]
     Finally, by the uniqueness of the limit, the last limit is valid for any sequence $\{\lambda_n\}\subset(0,\lambdazero)$ tending to zero, which proves the convergence in  \eqref{lim:claim}. 
\end{proof}
We can now establish the nonlocal-to-local limit, provided that the solution of the nonlocal equation satisfies some a priori estimates.

\begin{proposition}
\label{prop:nonlocal-to}
Let $c\in (0,\sqrt{2})$. Consider a family  $\{\W_\lambda\}_{\lambda\in (0,\lambdazero)}\subset\boM(\R)$ satisfying \ref{W:derivative}, \ref{W:m}, \ref{Wlambda:aeconvergence} and  \ref{Wlambda:finitesup}. For every $\lambda\in(0,\lambdazero)$, let $u_\lambda\in\NV$ be a nontrivial solution to \eqref{TWc} with $\W=\W_\lambda$. Assume that $\sup_{\lambda\in (0,\lambdazero)} \|u_\lambda\|_{L^\infty(\R)}<\infty$.
 Then there is a sequence $\{(\theta_\lambda,x_\lambda)\}\subset\R^2$,  such that, for every $k\geq 0$, up to a subsequence,
 \begin{equation}\label{eq:limittolocal:bis}
u_\lambda(\cdot +x_\lambda)e^{i\theta_\lambda} \to u_\infty, \quad  \text{  as }\lambda\to 0, \text{ in } W^{k,\infty}_{\loc}(\R),
\text{ and weakly-$*$ in }  L^\infty(\R), 
\end{equation}
 where $u_\infty \in L^\infty(\R)\cap \boC^\infty(\R)$, is a nonconstant solution to 
 \begin{equation}\label{eq:0}
	icu'+u''+u((1-|u|^2))=0,\quad \text{in }\R.
\end{equation}
Moreover,   $\inf_{\R}{\abs{u_\infty}}>0$ 
and the lifting $u_\infty=\rho_\infty e^{i\theta_\infty}$ satisfies $\rho'_\infty(0)=0$ and $\theta_\infty(0)=-\pi/2$.

Finally, if we additionally suppose that $u_\infty\in \boE(\R)$, then $u_\infty= \bu(\cdot \ ; c)$, 
where $\bu$ is the explicit soliton in \eqref{sol:1D}, and the convergence in \eqref{eq:limittolocal:bis} holds for the whole sequence.
 \end{proposition}

\begin{proof}
By Lemma~\ref{lemma:regularity}, we may consider the lifting $u_\lambda=\tilde\rho_\lambda e^{i\tilde\theta_\lambda}$. For every $\lambda\in(0,\lambdazero)$, since $\tilde \rho_\lambda\to 1$ as $\abs{x}\to \infty$, we can choose $(x_\lambda,\theta_\lambda)\in\R^2$ such that
$\abs{1-\tilde\rho_\lambda^2(x_\lambda)}=\max_{\R}\abs{1-\tilde\rho_\lambda^2}$ and   $\tilde\theta_\lambda(x_\lambda)+\theta_\lambda=-\pi/2$. Now we consider the translated and phase-shifted solution $v_\lambda=u_\lambda(\cdot+x_\lambda)e^{i\theta_\lambda}$. Clearly, denoting $\rho_\lambda=\tilde\rho_\lambda(\cdot+x_\lambda)$,  $\zeta_\lambda=\tilde\theta_\lambda(\cdot+x_\lambda)+\theta_\lambda$, 
and $\eta_\lambda=1-\rho_\lambda^2$, one has  
\begin{equation}\label{translatedsol}
v_\lambda=\rho_\lambda e^{i\zeta_\lambda},\quad\text{with }\rho_\lambda'(0)=0,\
\norm{\eta_\lambda}_{L^\infty(\R)}=\abs{\eta_\lambda(0)}, \ \zeta_\lambda(0)=-\pi/2,\,\text{ for all }\lambda\in(0,\lambdazero).
\end{equation}

 As explained in the proof of Proposition~4.5 in \cite{dLMar2022}, if $\|v_\lambda\|_{L^\infty(\R)}\leq C_0$, then $\| D^k v_\lambda \|_{L^\infty(\R)}\leq C_k$, for some constant $C_k$, for every $k\geq 1$. Indeed, by setting  
$w_\lambda(x)= v_\lambda(x) e^{\frac{ic}{2}x}$, for $x\in\R$, we have 
	\begin{equation}
		\label{eq:simply}
		-w_\lambda''=\Big(\frac{c^2}{4} +\W_\lambda\ast(1-\abs{w_\lambda}^2)\Big)w_\lambda\quad\text{ in }\R.
	\end{equation}
Then, by Young's inequality \eqref{young}, the fact that $\{\W_\lambda\}$ is bounded in $\boM(\R)$ implies that $\{w''_\lambda\}$
is bounded in $L^\infty(\R)$, and so is $\{ w'_\lambda \}$ by the Landau--Kolmogorov interpolation inequality. Hence, we can deduce that  
$\{ v_\lambda \}$ is bounded in $W^{2,\infty}(\R)$. The estimates for the higher order derivatives follow similarly, differentiating 
\eqref{eq:simply}. Therefore, for every $k\geq 0$, there exist a subsequence $\{v_{\lambda_n}\}$ and $u_\infty\in W^{k,\infty}(\R)$ such that 
for all $0\leq j\leq k,$
$
D^j v_{\lambda_n}  \wstar D^j  u_\infty$ in $L^\infty(\R)$,  as $n\to\infty$.
Moreover, by invoking the Ascoli--Arzela's theorem, the convergence also holds in $L^\infty_{\text{loc}}(\R)$.
Hence, by multiplying the equation 
\bq
 \label{eq:ulambda}
icv_\lambda'+v''_\lambda+v_\lambda\big(\W_\lambda*(1-|v_\lambda|^2)\big)=0,\quad \text{in }\R,
\eq
by a function in $\boC_c^\infty(\R)$, and using Lemma~\ref{lemma:limitnonlocalterm}, we can pass to the limit and deduce that $u_\infty$
is a solution to \eqref{TWc} with $\W=\delta_0$.

 In addition, hypotheses \ref{W:derivative} and \ref{W:m} allow us to apply Proposition~4.4 of \cite{dLMar2022}, so that,
 using also \eqref{young} and \ref{Wlambda:finitesup}, we get
$$(2-c^2)/4 \leq \|\W_\lambda*\eta_\lambda \|_{L^\infty(\R)}\leq M_{\lambdazero}\|\eta_\lambda \|_{L^\infty(\R)}.$$
 In view of  \eqref{translatedsol}, we conclude that 
 $$\abs{\eta_\lambda(0)}\geq  (2-c^2)/(4  M_{\lambdazero}).$$
 Thus, thanks to the convergence in $L^\infty_\loc$, we deduce that ${u_\infty(0)}\neq 1$
 so $u_\infty$ is a nontrivial solution.
In addition, since $\| D^k v_\lambda \|_{L^\infty(\R)}\leq C_1$, Proposition~4.7 of \cite{dLMar2022} implies that, 
for all $\lambda\in(0,\lambdazero)$, 
\bq
\label{cota:rho}
\rho_\lambda(x)= \abs{v_\lambda(x)}\geq 
\delta=\frac{ \sqrt{ 1+4c^2/C_1}-1}{\sqrt{ 1+4c^2/C_1}+1 }\in(0,1).
\eq
 Thus, by the convergence in $L^\infty_\loc$, we deduce that $\inf_\R\abs{u_\infty}\geq \delta,$ so that $u_\infty$ does not vanish, 
 ad we have the lifting $u_\infty=\rho_\infty e^{i\theta_\infty}.$ Furthermore, from \eqref{cota:rho}, we also deduce that $\{ \rho_\lambda \}$ is bounded in $W^{k,\infty}(\R)$, so  that $\{ \rho_\lambda \}$ must converge to $\rho_\infty$ in $W^{k,\infty}_\loc(\R)$,  
Also, by Proposition~{2.4} in \cite{dLMar2022}, we can write, in view of \eqref{expcontradiction},
\[\zeta_\lambda(x)=\frac{c}{2}\int_0^x\Big(\frac{1}{\rho_\lambda(y)^2}-1\Big)dy - \frac{\pi}{2}.\]
Thus, by the dominated convergence theorem, we deduce that 
$v_\lambda $ converges pointwise to $\rho e^{i\zeta}$, with 
$\zeta(x)=\frac{c}{2}\int_0^x(\rho_\infty(y)^{-2}-1)dy - \frac{\pi}{2}.$
Thus, we must have $\theta_\infty=\zeta+2\ell \pi$, for some $\ell \in  \N$, and we can take $\ell =0$.

In conclusion, using \eqref{translatedsol}, $u_\infty=\rho_\infty e^{i\theta_\infty}$ is a nontrivial solution  satisfying  $\rho'_\infty(0)=0$ and $\theta_\infty(0)=-\pi/2$. If we suppose additionally that $E(u_\infty)$ is finite, 
$u_\infty=\bu(\cdot\,;c)$, since it is the only nontrivial finite-energy solution to \eqref{TWc} (with $\W=\delta_0$), whose modulus has a critical point at zero and whose phase at zero is $-\pi/2$. In this case, the identification of the limit yields that the entire family $\{v_\lambda\}$ converges to $\bu(\cdot\,;c)$.
\end{proof}

The final ingredient to show  Theorem~\ref{thm:limittolocalgray} is to establish that the limit function $u_\infty$
has finite energy. This is the purpose of the following result.

\begin{lemma}\label{lemma:limitc}
Let $c\in (0,\sqrt{2})$. Consider a family of potentials $\{\W_\lambda\}_{\lambda\in (0,\lambdazero)}\subset\boM(\R)$ satisfying \ref{W:derivative}, \ref{W:m}, \ref{Wlambda:aeconvergence} and  \ref{Wlambda:finitesup}. For every $\lambda\in(0,\lambdazero)$, let $u_\lambda\in\NV$ be a nontrivial solution to \eqref{TWc} with $\W=\W_\lambda$.
Assume that 
$\sup_{\lambda\in (0,\lambdazero)}(\|u_\lambda\|_{L^\infty(\R)}+\|\rho'_\lambda\|_{L^2(\R)})<\infty$,  and that there exists $C>0$ such that
\begin{equation}\label{almostEpestimate}
		\int_\R \Big(1-\frac{c^2}{2\rho_\lambda^2}\Big)\eta_\lambda^2\leq C,\quad\text{for all }\lambda\in (0,\lambdazero),
\end{equation}
where $\rho_\lambda=|u_\lambda|$ and $\eta_\lambda=1-\rho_\lambda^2$. Then, the conclusion 
in Proposition~\ref{prop:nonlocal-to} holds, and the function $u_\infty$ belongs to $\boE(\R)$.
\end{lemma}

\begin{proof}
We start by performing the translation and phase-shift as in the proof of Proposition~\ref{prop:nonlocal-to} to 
define $v_\lambda$ as in \eqref{translatedsol}, and we keep the notation $\rho_\lambda= \abs{v_\lambda}$, for simplicity,
so that $v_\lambda$ converges, up to a subsequence, to $u_\infty$, for instance, in $W^{1,\infty}_\loc(\R)$.
Bearing in mind \eqref{cota:rho}, we have
\bq
\label{cota:rhoL}
\delta \leq \rho_\lambda(x) \leq C_0, \quad \text{for all }x\in \R,\, \lambda \in(0,\lambdazero),
\eq
with $\delta\in (0,\sqrt 2).$ 
Also, setting $\eta_\lambda=1-\abs{v_\lambda}^2$, and using \eqref{eq:firstintegral}, we have 
\bq
\label{cota-vL}
 \abs{v_\lambda'}^2=
{c^2 \eta_\lambda^2}/{(4 \rho_\lambda^2)}+(\rho_\lambda')^2.
\eq
Since we are supposing that $\norm{\rho'_\lambda}_{L^2(\R)}\leq C_1$, for some $C_1>0$,
it  is enough to establish that there is $C_2>0$ such that $\norm{\eta_\lambda}_{L^2(\R)}\leq C_2$
to conclude that  $u_\infty\in \boE(\R)$. Indeed, in this case,
by \eqref{cota:rhoL} and \eqref{cota-vL}, we deduce that $\{v_\lambda'\}$ is bounded in $L^2(\R)$, 
so that using also the convergence in  $W^{1,\infty}_\loc(\R)$, 
we infer, up to a subsequence, the weak convergences:
$$
  \eta_\lambda\wto \eta_\infty =1-\abs{u_\infty}^2  \ \text { and }  \ v'_\lambda\wto u'_\infty \quad \text{in }L^2(\R),
$$
and thus, 
$$
 \norm{\eta_\infty}_{L^2(\R)} \leq \liminf_{\lambda\to 0} \norm{\eta_\lambda}_{L^2(\R)}\leq C_2
\ \text { and } \ 
\norm{u'_\infty}^2_{L^2(\R)} \leq \liminf_{\lambda\to 0} \norm{v'_\lambda}^2_{L^2(\R)} \leq \frac{c^2 C_2^2}{4\delta^2}+C_1^2,
$$
i.e.\ $u_\infty$ has finite energy.

To show that $\{\eta_\lambda\}$ is bounded in $L^2(\R)$, 
we consider the set $$S_\lambda^r=\{x\in\R: \rho_\lambda(x)<r\}, \quad \text{ for }r>0.$$
Then, arguing as in the proof of Theorem~{4.1} of \cite{dLMar2022}, 
but invoking Lemma~\ref{lemma:limitnonlocalterm}  to pass to the limit in the nonlocal term of the equation,
we can show that there exists $\tilde r \in (c/\sqrt{2},1) $ such that the sequence $\{ \abs{S_\lambda^{\tilde r}} \}_{\lambda\in (0,\lambdazero)}$ is 
 bounded. Hence, from \eqref{almostEpestimate} and \eqref{cota:rhoL},  we obtain
	\begin{align*}
		C &\geq 
		\int_{\R\setminus S_\lambda^{\tilde r}}\left(1-\frac{c^2}{2\rho_\lambda^2}\right)\eta_\lambda^2+\int_{S_\lambda^{\tilde r}}\left(1-\frac{c^2}{2\rho_\lambda^2}\right)\eta_\lambda^2
\geq \left(1-\frac{c^2}{2{\tilde r}^2}\right) \int_{ \R \setminus S_\lambda^{\tilde r} }
\eta_\lambda^2-\frac{c^2(1-\delta^2)^2}{2\delta^2}|S^{\tilde r}_\lambda|.
	\end{align*}   
As a consequence, using also the upper bound in \eqref{cota:rhoL},
$$\norm{\eta_\lambda}^2_{L^2(\R)}\leq K_1(1+|S_\lambda^{\tilde r}|),\quad\text{ for all } \lambda\in(0,\lambdazero),$$
for some constant $K_1>0$ independent of $\lambda$, which completes the proof.
\end{proof}

\begin{proof}[Proof of Theorem~\ref{thm:limittolocalgray}]

First, let us denote $M_\lambda=M(c,\tau_\lambda,\sigma_\lambda)$ in case \ref{alternativetau} holds, while we denote $M_\lambda=M(c,\mu_\lambda)$ in case \ref{alternativemu} holds. In any of the cases  \ref{alternativetau} and \ref{alternativemu} we deduce, from \eqref{est:L1} and \eqref{est:deltamu} respectively, that $\{u_\lambda\}_{\lambda\in(0,\lambdazero)}$ is bounded in $L^\infty(\R)$.  Moreover, going through the proof of Corollary~{4.16} of \cite{dLMar2022}, one deduces (after taking limits $c_n\to c$, see the proof for more details) that
\begin{equation}\label{estimategamma}
\big(1-m_\lambda M_\lambda\big)\|\rho'_\lambda\|_{L^2(\R)}^2\leq \gamma_{\lambda,\gc}(b),
\end{equation}
for any $b\in (c,\sqrt{2})$ and any $\gc\in (0,c)$,  where $\gamma_{\lambda,\gc}(b)$ is the mountain pass level of the functional 
\[J_{\lambda,b}(1-\rho)=\frac12\int_\R(\rho')^2 + \frac14\int_\R \big(\W_\lambda*(1-\rho^2)\big)(1-\rho^2) - \frac{b^2}{8}\int_\R\frac{(1-\rho^2)^2}{\rho^2},\quad \rho\in 1+H^1(\R).\]
More precisely,
\[\gamma_{\lambda,\gc}(b)=\inf_{g\in\Gamma_\lambda(\gc)}\max_{t\in[0,1]} J_{\lambda,b}(g(t))>0,\]
where
\[\Gamma_\lambda(\gc)=\{g\in\boC([0,1],H^1(\R)):\,\, g<1\text{ in }[0,1],\,\, g(0)=0,\,\, g(1)=1-\phi_{\lambda,\gc}\},\]
and $\phi_{\lambda,\gc}\in 1+H^1(\R)$ is a fixed function such that $\phi_{\lambda,\gc}>0$ in $[0,1]$ and $J_{b,\lambda}(1-\phi_{\lambda,\gc})\in (-\infty,0).$

Recall that $\sup_{\lambda\in(0,\lambdazero)}\|\wh\W_\lambda\|_{L^\infty(\R)}\in (0,\infty)$. Hence, going through the proof of Lemma~{3.6} of \cite{dLMar2022} one can easily check that $\phi_{\lambda,\gc}=\phi_\gc$ can actually be chosen to be independent of $\lambda$. Thus, $\Gamma_\lambda(\gc)=\Gamma(\gc)$ is independent of $\lambda$ as well. Moreover, for every $\rho\in 1+H^1(\R)$ one has
\begin{align*}
J_{\lambda,b}(1-\rho) &=J_{0,b}(1-\rho)+\frac14\int_\R\big((\W_\lambda-\delta_0)*(1-\rho^2)\big)(1-\rho^2)
\\
&\leq J_{0,b}(1-\rho)+\|\wh\W_\lambda-1\|_{L^\infty(\R)}\|1-\rho^2\|^2_{L^2(\R)}.
\end{align*}
Taking maxima and infima in the previous inequality, we get
\[\gamma_{\lambda,\gc}(b)\leq \gamma_{0,\gc}(b) + \|\wh\W_\lambda-1\|_{L^\infty(\R)}\inf_{g\in\Gamma(\gc)}\max_{t\in[0,1]} \|g(t)(2-g(t))\|_{L^2(\R)}^2\leq C,\]
for some $C>0$ independent of $\lambda$.

In sum, from \eqref{estimategamma} we deduce that $\{\rho'_\lambda\}_{\lambda\in(0,\lambdazero)}$ is bounded in $L^2(\R)$. Finally, Lemma~{4.15} in \cite{dLMar2022} gives \eqref{almostEpestimate}. We can now apply Lemma~\ref{lemma:limitc} and conclude the proof.
\end{proof}

\section{Existence of symmetric  solitons for small \texorpdfstring{$\lambda$}{lambda} }\label{sec:even}

This section is devoted to the proof of Theorem~\ref{thm:even}, which establishes the existence of solutions to the nonlocal problem \eqref{TWc} bifurcating from the solutions to the local one. In order to have control over the bifurcation parameter independent of the velocity $c$, we will apply the following quantitative version of the implicit function theorem, as stated in \cite{bahsoun2024}.

\begin{theorem}[\cite{bahsoun2024}]
\label{thm:q:implict}
    Consider  $\boA,\boB$  Banach spaces,  $\Omega\subset\boB$  an open set, the interval $\Lambda=(-\lambdazero,\lambdazero)$ for some $\lambdazero>0$, and  $G:\Omega\times\Lambda\to \boA$ a $\boC^1$ function. Assume that there exists $x_0\in\Omega$ such that $G(x_0,0)=0$, $\partial_x G(x_0,0)$ is invertible and 
    \begin{align*}
        &\sup_{(x,\lambda)\in V_\delta}\|\boI-[\partial_x G(x_0,0)]^{-1}\partial_x G(x,\lambda)\|_{\boL(\boB)} \leq\frac{1}{2},
        \\
        &\sup_{(x,\lambda)\in V_\delta} \|[\partial_x G(x_0,0)]^{-1}\partial_\lambda G(x,\lambda)\|_{\boL(\R,\boB)}\leq M,
    \end{align*}
    for some $\delta,M>0$, where $\boI$ is the identity operator in $\boB$ and 
    \[V_\delta=\{(x,\lambda)\in\Omega\times\Lambda:\,\|x-x_0\|_{\boB}\leq \delta,\, |\lambda|\leq\delta\}.\]
    Then, denoting $\Lambda_{\delta,M}=(-\delta/2M,\delta/2M)$, there exists a $\boC^1$ function $g:\Lambda_{\delta,M}\to\Omega$ such  that all the solutions to the equation $G(x,\lambda)=0$ in the set 
    $\boK=\{(x,\lambda)\in \Omega\times \Lambda_{\delta,M} :\, \|x-x_0\|_{\boB}\leq\delta\}$ are given by $\{(\lambda,g(\lambda)):\,\lambda\in\Lambda_{\delta,M}\}$. 
\end{theorem}

In view of equation \eqref{eq:eta} in Lemma~\ref{lemma:regularity}, we introduce the open sets of $H^2(\R)$
$$
\Omega=\{\eta\in H^2(\R): \eta<1\,\,\text{ in }\R\}, \quad \Omega_e=\{\eta\in \Heven: \eta<1\,\,\text{ in }\R\},
$$
the interval  $\Lambda=(-\lambdazero,\lambdazero)$ for some $\lambdazero>0$,
and the functional 
$G_c:\Omega\times\Lambda\to L^2(\R),$ given by 
\begin{equation}\label{eq:operator}
G_c(\eta,\lambda)=\eta''-2\W_\lambda\ast\eta+c^2\eta+\frac{c^2\eta^2}{2(1-\eta)}+\frac{(\eta')^2}{2(1-\eta)}+2\eta(\W_\lambda\ast\eta),\quad\text{for all }(\eta,\lambda)\in \Omega\times\Lambda,
\end{equation}
for some $c\in (0,\sqrt 2)$ and $\{\boW_\lambda\}_{\lambda\in\Lambda}\subset\boM(\R)$. Notice that $G_c(\boldeta_c,0)=0$, where $\boldeta_c$ is given by \eqref{sol:1D}. Furthermore, by Proposition~2.4 in \cite{dLMar2022}, we can recover a solution  $u_\lambda$ to \eqref{TWc}, with $\W=\W_\lambda$, by determining a zero $\eta_\lambda$ of  $G_c(\cdot, \lambda)$. In such a case, $u_\lambda$ is given by
\begin{equation}
	\label{def:u}
	u_\lambda(x)=\sqrt{1-\eta_\lambda(x)}e^{i\theta_\lambda(x)}, \quad \textup{where}\quad  \theta_\lambda(x)=\frac{c}2\int_0^x \frac{\eta_\lambda(s)}{1-\eta_\lambda(s)}{\rm d}s-\frac \pi2.
\end{equation}

The next lemma establishes some properties to control $G_c$ quantitatively.

\begin{lemma}
\label{lem:C1}
 The functional  $G_c$ is $\boC^1(\Omega\times\Lambda, L^2_{}(\R))$ for every $c\in (0,\sqrt{2})$, with partial derivatives
 \begin{equation}\label{dlambdaJ}
\partial_\lambda G_c(\eta,\lambda)=2(\eta-1)\left((\partial_\lambda\W_\lambda)\ast\eta\right),
\end{equation}
and 
 \bq
\label{dem:dereta}
\begin{aligned}
\partial_\eta G_c(\eta,\lambda)(\sigma)&=\sigma''-2\W_\lambda\ast\sigma+c^2\sigma
\\
&+\frac{c^2\eta(2-\eta)+(\eta')^2}{2(1-\eta)^2}\sigma+\frac{\eta'}{1-\eta}\sigma'+2(\W_\lambda\ast\eta)\sigma+2\eta(\W_\lambda\ast\sigma), 
\quad \sigma\in H^2(\R).
\end{aligned}
\eq
In addition, there is $C>0$, independent of $c$ and $\lambda$, such that 
\bq
\label{est:derL}
\norm{\partial_\lambda G_c(\eta,\lambda)}_{\boL(\R,L^2(\R))}\leq C(1+\norm{\eta}_{H^2(\R)}^2).
\eq

  Moreover, for all $a \in (0,\sqrt 2)$, $\alpha >0$, $\beta \in (0,1)$ and $\ve >0$, there exists $\delta >0$ such that the following uniform continuity holds: For all $\eta,\bar \eta\in \Omega$ and $\lambda,\bar \lambda\in \Lambda$ satisfying
\begin{equation}\label{uniformcont}
     \|\eta-\bar{\eta}\|_{H^2(\R)}+|\lambda-\bar{\lambda}|<\delta,
     \quad 
    |\bar{\eta}'|\leq \alpha \,\text{ in }\R, \quad \text{ and } \quad  -\alpha\leq\bar{\eta}\leq 1-\beta\,\,\text{ in }\R,
\end{equation}
 we have, 
\bq
\label{dem:ineqeps}
\|\partial_\eta G_c(\eta,\lambda)-\partial_\eta G_c(\bar{\eta},\bar{\lambda})\|_{\boL(H^2(\R),L^2(\R))} <\varepsilon,
\eq
for all $c\in [a,\sqrt 2)$.
\end{lemma}

\begin{remark}
For our arguments, we do not need to verify the differentiability of the function $\lambda\in \Lambda\mapsto \W_\lambda\in \boM(\R)$. We only need to perform computations in the Fourier variable. 
For this reason, we only need to interpret the expression $\partial_\lambda \W_\lambda$ in \eqref{dlambdaJ} as a tempered distribution in the following sense.
By \ref{Wlambda:C1}, $\partial_\lambda\wh\W_\lambda$ is defined by
\bq
\label{def:derhatW}
\lim_{t\to 0}\frac{\wh\W_{\lambda+t}(\xi)-\wh\W_\lambda(\xi)}{t}=\partial_\lambda\wh\W_\lambda(\xi),\quad\text{for a.e. }\xi\in\R,
\eq
and the  growth condition implies that $\partial_\lambda \wh \W_\lambda$ is in $\boS'(\R)$, so that we can define $\partial_\lambda\W_\lambda$ in $\boS'(\R)$ as 
 $$
 \partial_\lambda\W_\lambda=\boF^{-1}(\partial_\lambda \wh \W_\lambda),
 $$
where $\boF$ denotes the Fourier transform in the space variable. 
 In this manner, 
\bq 
\label{dem:idF}
\boF\left((\partial_\lambda\W_\lambda)\ast\eta\right)=\partial_\lambda \wh \W_\lambda \cdot \wh  \eta,\quad\text{for all }\eta \in \boS(\R),
\eq
and this equality still holds for any $\eta \in H^2(\R)$, because of the growth hypothesis in \ref{Wlambda:C1}.

\end{remark}

\begin{proof}[Proof of Lemma~\ref{lem:C1}]
Let us check first that, for all  $(\eta,\lambda)\in\Omega\times\Lambda$, we have \eqref{dlambdaJ}.
 Indeed, by \eqref{dem:idF} and Plancherel's identity, we have the following estimate, for every $t\in(-\lambda-\lambdazero,-\lambda+\lambdazero)$,
\begin{align}
&\Big\|\frac{(\eta-1)((\W_{\lambda+t}  -\W_\lambda)\ast\eta)}{t} -(\eta-1)(\partial_\lambda\W_\lambda )\ast\eta\Big\|_{L^2(\R)}
 \leq  \label{dem:der1} \\
 &\|\eta-1\|_{L^\infty(\R)}\Big\|\Big(\frac{\W_{\lambda+t}-\W_\lambda}{t}-\partial_\lambda\W_\lambda\Big)\ast\eta\Big\|_{L^2(\R)}
=\frac{\|\eta-1\|_{L^\infty(\R)}}{2\pi}\Big\|\Big(\frac{\wh\W_{\lambda+t}-\wh\W_\lambda}{t}-\partial_\lambda\wh\W_\lambda\Big)\widehat{\eta}\Big\|_{L^2(\R)}.\nonumber
\end{align}

To apply the dominated convergence theorem in $L^2(\R)$, we see that, in view of hypothesis \ref{Wlambda:C1} and  the mean value theorem, there is $C>0$, independent of $\lambda$, such that
\[\left|\frac{\wh\W_{\lambda+t}(\xi)-\wh\W_\lambda(\xi)}{t}-\partial_\lambda\wh\W_\lambda(\xi)\right|\leq {C(1+\xi^2)},\quad\text{for a.e. }\xi\in\R,\]
then
\begin{equation*}
\Big| {\frac{\wh\W_{\lambda+t}-\wh\W_\lambda}{t}-\partial_\lambda\wh\W_\lambda} \Big|\cdot | \wh\eta|\leq C (|\wh\eta|+|\wh{\eta''}|)\in L^2(\R).
\end{equation*}
Bearing in mind \eqref{def:derhatW}, we may apply the dominated convergence theorem in \eqref{dem:der1}, so that \eqref{dlambdaJ} is proved. The continuity of $\partial_\lambda G_c:\Omega\times\Lambda\to L^2(\R)$ and the estimate in \eqref{est:derL} follow using the same arguments. 

On the other hand, standard computations show that $\partial_\eta G_c:\Omega\times\Lambda\to\boL(H^2(\R),L^2(\R))$  has the expression in \eqref{dem:dereta}.
Finally, it remains to show \eqref{dem:ineqeps}, which in particular gives the continuity of $\partial_\eta G_c$. Indeed, let 
 $\eta,\bar \eta\in \Omega$ and $\lambda,\bar \lambda\in \Lambda$ satisfy \eqref{uniformcont}  for some $\delta>0$ to be chosen small enough, and let $\sigma\in H^2(\R)$ with $\|\sigma\|_{H^2(\R)}=1$. We start analyzing the second term on the right-hand side of \eqref{dem:dereta}. From Plancherel's identity, the mean value theorem and \ref{Wlambda:C1}, it follows that for a.e. $\xi\in\R$ there is some $\tilde\lambda\in\Lambda$ such that 
\begin{align*}
    \|\W_\lambda\ast\sigma-\W_{\bar{\lambda}}&\ast\sigma\|^2_{L^2(\R)}=\frac{1}{2\pi}\int_\R\Big(\wh\W_\lambda(\xi)-\wh\W_{\bar{\lambda}}(\xi)\Big)^2|\wh\sigma(\xi)|^2d\xi
    \\
    &\leq \frac{|\lambda-\bar\lambda|^2}
    {2\pi}\int_\R |\partial_\lambda\wh\W_{\tilde\lambda}(\xi)|^2|\wh\sigma(\xi)|^2 d\xi 
    \leq  \frac{C|\lambda-\bar\lambda|^2}
    {2\pi}\int_\R (1+\xi^4)|\wh\sigma(\xi)|^2 d\xi 
    \\
    &= C|\lambda-\bar\lambda|^2\int_\R (|\sigma|^2 + |\sigma''|^2)\leq C|\lambda-\bar\lambda|^2<\ve,
\end{align*}
for $\delta^2<\ve/C$. We may also use the previous inequality to estimate the last term in \eqref{dem:dereta}, namely,
\begin{align*}
&\left\|\eta(\W_\lambda\ast\sigma)-\bar{\eta}(\W_{\bar{\lambda}}\ast\sigma)\right\|_{L^2(\R)} \leq \left\|\eta((\W_\lambda-\W_{\bar{\lambda}})\ast\sigma)\right\|_{L^2(\R)}+\left\|(\eta-\bar{\eta})(\W_{\bar{\lambda}}\ast\sigma)\right\|_{L^2(\R)}
\\
&\leq \|\eta\|_{L^\infty(\R)}\|(\W_\lambda-\W_{\bar{\lambda}})\ast\sigma\|_{L^2(\R)}+\|\eta-\bar{\eta}\|_{L^\infty(\R)}\|\W_{\bar{\lambda}}\ast\sigma\|_{L^2(\R)}<\ve,
\end{align*}
where we used that $\|\eta\|_{L^\infty(\R)}\leq\delta+\max\{\alpha,1-\beta\}$ and considered a smaller $\delta$. 
For the sixth term in \eqref{dem:dereta}, we  analogously get
\[\|(\W_{\lambda}\ast\eta)\sigma-(\W_{\bar{\lambda}}\ast\bar{\eta})\sigma\|_{L^2(\R)}<{\varepsilon}.\]

For the fourth term, we fix $\gamma\in(0,\beta)$ and consider the compact set 
\[\boK_\gamma=\{(x,y)\in\R^2: -\alpha-\gamma\leq x\leq 1-\beta+\gamma,\,\,|y|\leq\alpha+\gamma\}.\]
Observe that the functions $(x,y)\mapsto\frac{x(2-x)}{2(1-x)^2}$ and $(x,y)\mapsto\frac{y^2}{2(1-x)^2}$ are uniformly continuous in $\boK_\gamma$. Therefore, for any $S>0$, we may take $\delta>0$ even smaller (depending only on $a,\alpha,\beta,\gamma,S$) in such a way that, if $(x,y),(x_0,y_0)\in\boK_\gamma$ and $|x-x_0|+|y-y_0|<S\delta$, then
\[\left|\frac{c^2 x(2-x)+y^2}{2(1-x)^2}-\frac{c^2 x_0(2-x_0)+y_0^2}{2(1-x_0)^2}\right|\leq \ve.\]
In conclusion, bearing in mind \eqref{uniformcont} and the Sobolev embedding $H^2(\R)\hookrightarrow L^\infty(\R)$,  one may take $\delta>0$ small enough, depending only on $a,\alpha,\beta,\gamma$, such that 
\[\left\|\frac{c^2\eta(2-\eta)+(\eta')^2}{2(1-\eta)^2}\sigma-\frac{c^2\bar{\eta}(2-\bar{\eta})+(\bar{\eta}')^2}{2(1-\bar{\eta})^2}\sigma\right\|_{L^2(\R)}<{\varepsilon}.\]

The bounds for the fifth term on the right-hand side of \eqref{dem:dereta} follow similarly, which finishes the proof.
\end{proof}

With the purpose of applying Theorem~\ref{thm:q:implict}, we establish now that $\partial_\eta G_c(\boldeta_c,0)$ is invertible, when restricted to the even functions. For the sake of simplicity, let us denote
$\boT_c=\partial_\eta G_c(\boldeta(\cdot \, ;c),0) : H^2(\R)\to L^2(\R)$, whose expression is
\bq
\label{Tc}
\boT_c(\sigma)=\sigma''-(2-c^2)\sigma+\frac{c^2\boldeta_c(2-\boldeta_c)+(\boldeta_c')^2}{2(1-\boldeta_c)^2}\sigma+\frac{\boldeta_c'}{1-\boldeta_c}\sigma'+4\boldeta_c\sigma, \quad \sigma\in H^2(\R). 
\eq

\begin{lemma}
\label{lem:ker}
Let $c\in (0,\sqrt 2)$. We have  $\ker(\boT_c)=\{\boldeta_c'\}.$
\end{lemma}
\begin{proof}
It is immediate to check that  $\boT_c (\boldeta_c')=0$, so that $0$ is an eigenvalue. We use the following ODE argument to show that $0$ is simple in $H^2(\R)$.
Let $\sigma_1=\boldeta'_c$ and $\sigma_2 \in \ker(\boT_c),$  and set the Wronskian $W=\sigma_1\sigma_2'-\sigma'_1\sigma_2$. Since $\sigma_1,\sigma_2\in H^2(\R)$, we get 
$W(\infty)=0$. On the other hand, by \eqref{Tc} and the Abel identity,
$$W(x)=W(0)e^{-\int_0^x p(t)dt}, \quad \text{ with }p(t)={\boldeta_c'(t)}/{(1-\boldeta_c(t))}.
$$
Since $p\in L^1(\R)$ due to the exponential decay of the soliton, this implies that $W(0)=0$, and then $W\equiv 0$. Therefore, $\sigma_1,\sigma_2$ are linearly dependent and 
$\ker(\boT_c)=\{\boldeta_c'\}.$
\end{proof}

\begin{proposition}
\label{prop:invertible}
Let $c\in (0,\sqrt 2)$. The (restricted) linear operator   $\boT_{c,e}: \Heven\to \Leven$ is invertible.
Moreover, the map $c\in (0,\sqrt 2) \mapsto \boT^{-1}_{c,e}\in \boL(\Leven,\Heven)$ is continuous.
\end{proposition}
\begin{proof}
 Since  $\boldeta'_c$ is odd, it follows from Lemma~\ref{lem:ker} that  $\ker(\boT_{c,e})=\{0\}$. 
To show the bijectivity, we will use the Fredholm alternative in $\Heven$. For this purpose, we remark that the problem of finding $\sigma\in \Heven$ solution to 
$\boT_c(\sigma)=f$, for a given $f\in \Leven$, can be written  as
$$
\boL (\sigma)+\boB (\sigma)=f, \quad \text{with } \quad \boL(\sigma)= \sigma''-(2-c^2)\sigma \ \text{ and }\  \boB : \Heven\to \Leven. 
$$
Noticing that  $\boL : \Heven\to \Leven$ is a bijection, since $\boF(\boL(\sigma))(\xi)=-(\xi^2+(2-c^2))\wh\sigma(\xi),$ this can be recast as 
\bq
\label{dem:fred}
\sigma+\boK (\sigma)=\boL^{-1}(f), \quad \text{with } \quad \boK=\boL^{-1}\circ \boB : \Heven\to \Heven.
\eq
Because of the exponential decay of $\boldeta_c$, we infer that $\boB$ is a compact operator  (see e.g.\  Lemma~3.12 in \cite{stuart-elliptic}) and then so is  $\boK$.
By the Fredholm alternative in $\Heven$, \eqref{dem:fred} has a unique solution if and only if $\ker(\boI+\boK)=\{0\}$, which is equivalent to the fact that $\ker(\boT_{c,e})=\{0\}$. Therefore, there exists  $\boT_{c,e}^{-1}$, and belongs to $ \boL(\Leven,\Heven)$ by the Banach Isomorphism theorem.

To prove that $c\in (0,\sqrt 2) \mapsto \boT^{-1}_{c,e}\in \boL(\Leven,\Heven)$ is continuous, we will apply Theorem~1.16 in \cite{kato-linear}.
For notational simplicity, we will omit the subscripts in the operator norms in the sequel. Let us fix $c_0\in (0,\sqrt 2)$ and take $\ve>0$ arbitrary, but smaller than $1/\norm{\boT^{-1}_{c_0,e}}$. 
It is straightforward to check that $c\in (0,\sqrt 2) \mapsto \boT_{c,e}\in \boL(\Heven,\Leven)$ is continuous, so that there is $\delta>0$ such that, if $\abs{c-c_0}<\delta$, then 
\bq
\label{dem:kato}
\norm{\boT_{c,e}-\boT_{c_0,e}}< \ve.
\eq
Using the notation in Theorem~1.16 in \cite{kato-linear}, we take $T=\boT_{c_0,e}$, $S=\boT_{c,e}$, $A=S-T$. In view of \eqref{dem:kato}, $A$ is trivially $T$-bounded in the sense that $\norm{A\sigma}_{L^2(\R)}\leq \ve \norm{\sigma}_{H^2(\R)}$, for all $\sigma\in \Heven.$ Therefore, we deduce that 
$$\norm{S^{-1}-T^{-1}}\leq  \frac{\ve\norm{T^{-1}}^2}{1-\ve \norm{T^{-1}}},$$
which is the desired continuity of $\boT^{-1}_{c,e}$ at $c_0$.
\end{proof}

We are now in a position to complete the proof of the existence of solutions constructed via the quantitative implicit function theorem.

\begin{proof}[Proof of Theorem~\ref{thm:even}]
Let $a\in(0,\sqrt 2)$. For every $\delta>0$, let us define the set
\[V_{\delta,c}=\{(\eta,\lambda)\in \Omega_e \times\Lambda:\,\, \|\eta-\boldeta_c\|_{H^2(\R)}\leq\delta,\,\,|\lambda|\leq\delta\}.\]
One may choose $\delta>0$ independent of $c$ such that $V_{\delta,c}\subset \Omega\times\Lambda$. Indeed, let $(\eta,\lambda)\in V_{\delta,c}$. Then,
by the Sobolev embedding $H^2(\R)\hookrightarrow L^\infty(\R)$,  there is $S>0$ such that,
\[\eta\leq\|\eta-\boldeta_c\|_{L^\infty(\R)}+\|\boldeta_c\|_{L^\infty(\R)}\leq S\|\eta-\boldeta_c\|_{H^2(\R)}+\frac{2-c^2}{2}\leq S\delta+\frac{2-a^2}{2},\quad\text{in }\R.\]
Thus, for $\delta<{a^2}/(2S)$, we have  $\eta<1$, so it belongs to $\Omega$.

On the other hand, by Proposition~\ref{prop:invertible}, we deduce that $c\in [a,\sqrt{2-a^2}]\mapsto \norm{\boT_{c,e}^{-1}}_{\boL(\Leven,\Heven)} \in \R$ is uniformly continuous, so that there is 
a constant $\gamma_a>0$ such that, 
\bq\label{dem:gamma}
\|\boT_{c,e}^{-1}(f)\|_{H^2(\R)}\leq \gamma_a \|f\|_{L^2(\R)},\quad\text{for all }f\in \Leven,\quad \text{ for all }c\in [a,\sqrt{2-a^2}].
\eq

Let us fix now $(\eta,\lambda)\in V_{\delta,c}$. 
Taking 
$$\alpha=\max_{c\in [a,\sqrt{2-a^2}]}\norm{\boldeta'_c}_{L^\infty(\R)}\quad  \text{ and }\quad \beta=\min_{c\in [a,\sqrt{2-a^2}]}\norm{1-\boldeta_c}_{L^\infty(\R)}$$ in \eqref{uniformcont}, we conclude from the estimates in \eqref{dem:ineqeps} and \eqref{dem:gamma}  that
 \begin{align*}
 &\| \boI- {\boT_{c,e}^{-1}} \left(\partial_\eta G_c(\eta,\lambda)\right) \|_{\boL(\Heven,\Heven)} \\
  &\leq
  \norm{\boT_{c,e}^{-1}}_{\boL(\Leven,\Heven)} 
  \| \partial_\eta G_c(\boldeta,0)-\partial_\eta G_c(\eta,\lambda) \|_{\boL(\Heven,\Leven)} \leq \gamma_a \ve,    
 \end{align*}
 where $\boI$ is the identity in $\Heven$.
Finally, invoking \eqref{est:derL} in Lemma~\ref{lem:C1}, we deduce the existence of a constant $B_{\delta,a}$, such that, for all $(\eta,\lambda)\in V_{\delta,c}$, 
for all $c\in [a,\sqrt{2-a^2}]$,
$$
\| \partial_\lambda G_c(\eta,\lambda) \|_{\boL(\R,\Leven)}\leq C(1+(\delta+\norm{\boldeta_c}_{H^2(\R)})^2) \leq B_{\delta,a},
$$
so that,   for all $(\eta,\lambda)\in V_{\delta,c}$, 
 $$
\|\boT_{c,e}^{-1}\circ \partial_\lambda G(\eta,\lambda)\|_{\boL(\Heven)}\leq \gamma_a B_{\delta,a}. 
 $$
Invoking  Theorem~\ref{thm:q:implict} with $\boB=\Heven$, $\Omega=\Omega_e$, $\ve=1/(2\gamma_a)$ and $M=\gamma_a B_{\delta,a}$, 
and setting $\lambda_a=\delta/(2M)$ and $\Omega_e^a$ the intersection of $\Omega_e$ with the open ball (in $\Heven$) of center $\boldeta_c$ and radius $\delta$, 
we deduce the existence of the $\boC^1$-function $g : (-\lambda_a,\lambda_a) \to \Omega_e^a$, where $\lambda_a=\delta/(2M)$.
 Therefore, the conclusions in Theorem~\ref{thm:q:implict}  follow setting $\eta_\lambda=g(\lambda)$, and defining
 $u_\lambda=\sqrt{1-\eta_\lambda}\exp(i\theta_\lambda)$ according to \eqref{def:u}.
In addition, since $g$ is differentiable, 
$g(\lambda)=g(0)+o(1)$ in $H^2(\R)$. Therefore, since
 $g(0)=\boldeta_c$, we conclude that $\eta_\lambda\to\boldeta_c$ in $H^2(\R)$.
 Notice in particular that $\norm{\eta_\lambda}_{L^\infty(\R)}=\eta_\lambda(0)$, since $\eta_\lambda$ is even and  $\eta_\lambda\to\boldeta_c$ in $L^\infty(\R)$, by the Sobolev embedding. Moreover, 
 there is $\delta>0$ such that $1-\eta_\lambda(0)\geq \delta$, for all $\lambda$ small.

The convergence in $H^k(\R)$, for $k\geq 3$, follows by a bootstrap argument. Indeed, since $\eta_\lambda$ satisfies \eqref{eq:eta}, 
we can differentiate this equation to obtain an expression for $\eta_\lambda'''$ in terms of lower order derivatives.
Using \ref{Wlambda:finitesup} and arguing as in the proof of Lemma~\ref{lemma:limitnonlocalterm}, we deduce that for any sequence of functions $f_\lambda \to f\in L^2(\R)$, we have, for a subsequence $\{\lambda_n\}$, $\boW_{\lambda_n}*f_{\lambda_n} \to f\in L^2(\R),$
as $n\to\infty$. Thus, it is simple to verify that we can pass to the limit in all the lower order terms, to conclude that $\eta_{\lambda_n}'''$
converges to $\boldeta_c'''$ in $L^2(\R)$. The convergence of the whole sequence follows from the uniqueness of the limit. The same argument establishes the convergence in
$H^k(\R)$ for $k\geq 4$.

 Concerning $u_\lambda$, it is easy to see that $\theta_\lambda$ converges pointwise to $\boldtheta_c$ in \eqref{u:theta}.
 Let us show that this convergence is uniform in $L^\infty(\R)$. Since
\bq
\label{diff:eta}
\abs{\theta_\lambda(x)-\boldtheta_c(x)}\leq \frac1{(1-\eta_\lambda(0))(1-{\boldeta_c}(0))}\norm{\eta_\lambda-\boldeta_c}_{L^1(\R)},
\quad \text{ for all }x\in \R,
\eq
it is enough to show that $\eta_\lambda\to\boldeta_c$ in $L^1(\R)$. Indeed, 
by  the convolution equation in \eqref{eq:conv} and  the decay of $\boL_c$ in Lemma~\ref{lem:decay}, we have in particular that $\boL_c \in L^1(\R)$, and thus
\bq
\norm{\eta_\lambda-\boldeta_c}_{L^1(\R)}\leq \norm{\boL_c}_{L^1(\R)}\norm{F(\eta_\lambda)-F(\boldeta_c)}_{L^1(\R)}.
\eq
In addition, it is easy to check that the convergence of $\eta_\lambda$ in $H^2(\R)$ implies that $F(\eta_\lambda)$ converges to $F(\boldeta_c)$ in ${L^1(\R)}$.
Therefore, \eqref{diff:eta} gives  that $\theta_\lambda\to \boldtheta_c$ in $L^\infty(\R)$. Finally,  using that $\abs{1-e^{ix}}\leq \abs{x}$,
for all $x\in \R$, we conclude that
$$
\abs{e^{i\theta_\lambda(x)}-e^{i\boldtheta_c(x)}}=\abs{e^{-i\boldtheta_c(x)}}\abs{e^{i(\theta_\lambda(x)-\boldtheta_c(x))}-1}\leq \abs{\theta_\lambda(x)-\boldtheta_c(x))},
$$
for all $x\in \R$. Finally,  using \eqref{u:theta}, it is straightforward to conclude that $u_\lambda$ converges to $\bu_c$ in $\boC^{k}(\R)$, provided that $\eta_\lambda$ converges to $\boldeta_c$ in $H^{k+1}(\R)$.
\end{proof}

\section{Solitons with oscillations}\label{sec:monotonicity}
In this section will prove Theorem~\ref{thm:nonmonotone} concerning the non-monotonicity of the profile function $\eta=1-\abs{u}^2$.
Our arguments use some ideas introduced by Berestycki, Nadin,
 Perthame and  Ryzhik \cite{BerNadin} for the study of the nonlocal KPP-Fisher equation.
We start by establishing some a priori estimates on the linearized problem.
\begin{lemma}\label{lemma:estimateslinear}
 Assume that $\W$ satisfies \ref{W:L1}, for some constants $\sigma,\tau>0$.
Let $R>1$. Consider  $g\in L^\infty((-\infty,R))$ and $h\in W^{1,\infty}((-\infty,R))$ satisfying
$$\norm{g}_{L^{\infty}((-\infty,R))}+\norm{h}_{W^{1,\infty}((-\infty,R))}\leq L,$$ for some $L\geq 1$.
 Assume that   $w\in \boC^2(\R)$ is a solution to 
\begin{equation}\label{eq:w}
	-w''+2\W*w=g w + h w',\quad \text{ in } (-\infty,R),
\end{equation}
 with $w\geq 0$ in $\R$,  $w\leq 1$ in  $(-\infty,0)$, $w(0)=1$, 
and $w'\geq 0$ in $(-\infty,R)$.
Then, 
there are constants $K,N$, depending only on $\sigma$ and $\tau$, such that, 
\begin{equation}
\label{w:harnack}
w(t+1)\leq \ (L K)^N w(t),\quad\text{for all } t \in (-\infty, R-2).
\end{equation}
Moreover, for any $a,b,\tilde b\in\R$, with $a<b<\tilde b<R$, there is a constant $C$, depending only on $a,b,\tilde b$, such that 
\begin{equation}
\label{est:der:w}
	\norm{w'}_{L^2((a,b))}+ \norm{\W*w}_{L^1((a,b))}\leq C L { w(\tilde b)}.
 \end{equation}
 In particular, there is a constant $\tilde C$, depending only on $a,b,\tilde b$, such that
\begin{equation}\label{est:w12loc}
 \norm{w}_{W^{2,1}((a, b))}\leq \tilde C L{ w(\tilde b)}.
\end{equation}
 \end{lemma}

\begin{proof}
To show \eqref{w:harnack}, we consider $a,b,\ve>0$ such that 
\begin{equation}
\label{a:b:tau}
    0<a<b<b+\ve <\tau, \ \text{ and }\ -\infty< t < R-\ve. 
\end{equation}
Let $\varphi\in \boC^2_c(\R)$ be a cut-off function such that $\varphi>0$ in $(-b,-a)$, and $\varphi=0$ in $\R\setminus(-b,-a),$
so that $\vp_t(x)=\vp(x-t)$ is localized in $(t-b,t-a)$.
Multiplying \eqref{eq:w} by $\varphi_t$ and integrating by parts, we get
\begin{equation*}
	\int_{t-b}^{t-a} (\W*w)\varphi_t=\frac12\int_{t-b}^{t-a}(g\varphi_t-(h\varphi_t)'+\varphi_t'')w\leq 
	\frac{3L+1}{2}\norm{\varphi}_{W^{2,\infty}(\R)}\int_{t-b}^{t-a} w.
	\end{equation*}
Since $w$ is nondecreasing in $(-\infty,R)$, and since $L\geq 1$, we conclude that 
\begin{equation}
\label{proof:harnack}
\int_{t-b}^{t-a} (\W*w)\varphi_t\leq 2L\norm{\varphi}_{W^{2,\infty}(\R)}(b-a)w(t-a).
\end{equation}
To estimate the left-hand side of \eqref{proof:harnack},
 using again that  $w$ is nonnegative and  nondecreasing in $(-\infty,R)$, we see that
\begin{align*}
	\int_{t-b}^{t-a}&\varphi_t(x)\int_{-\infty}^\infty\W(x-y)w(y)dydx \geq 
	\int_{t-b}^{t-a}\varphi_t(x)\int_t^{t+\ve}\W(x-y)w(y)dydx
	\\
	&\geq w(t)\int_{t-b}^{t-a}\varphi_t(x)\int_{t}^{t+\ve}\W(x-y)dydx=w 
	(t)\int_{t-b}^{t-a}\varphi_t(x)\int_{x-t-\ve}^{x-t}\W(y)dydx.
\end{align*}
Since, by \eqref{a:b:tau}, for any $x\in[t-b,t-a]$, we have $x-t-\ve> -\tau$ and $x-t\leq -a$, 
we deduce that we can bound $\W(y)$ from below by $\sigma$ in the integral, so that
\begin{equation*}
\int_{t-b}^{t-a} (\W*w) \vp_t\geq \ve  \sigma \int_{t-b}^{t-a}\varphi_t(x)dx=\ve \sigma\norm{\vp}_{L^1(\R)}.
\end{equation*}	
Therefore, combining with \eqref{proof:harnack}, we conclude that 
\begin{equation}
\label{dem:har1}
w(t)\leq KL w(t-a), \quad \text{ for all }t\in(-\infty,R-\ve),
\end{equation}
with $K=2\norm{\varphi}_{W^{2,\infty}(\R)}(b-a)/(\ve \sigma \norm{\vp}_{L^1(\R)})$.
For instance, we can take  $\ve=b=\min\{\tau/4,1/4\}$ and $a=b/2$. From \eqref{dem:har1} we immediately deduce 
\begin{equation}
\label{dem:har}
w(t+1)\leq KL w(t+1-a), \quad \text{ for all }t\in(-\infty,R-2).
\end{equation}
Then, for $t<R-2$, writing $1=Na+\delta$ for some
$N\in\N $ and $\delta \in(0,a)$,
we conclude by iterating \eqref{dem:har} and by the monotonicity of $w$, that 
$$w(t+1)\leq (KL)^N w(t +1 - Na)=(KL)^N w(t-\delta)\leq  (KL)^N w(t).$$
This finishes the proof of \eqref{w:harnack}.

Let us now establish the local estimates in \eqref{est:der:w}.
Take  $\varphi\in \boC^1(\R)$ a cut-off function  such that
$\varphi\geq 0$ in $\R$, $\varphi=1$ in  $(a,b)$, and
$\varphi= 0$ in  $\R\setminus(\tilde a,\tilde b)$, with $\tilde a=a-1$.
By multiplying \eqref{eq:w} by $w\varphi^2$, dropping the (positive) nonlocal term and integrating by parts, we deduce 
\begin{equation}\label{eq:H1est}
	\int_\R (w')^2\varphi^2+2\int_\R w'w\varphi'\varphi\leq\int_\R gw^2\varphi^2+\int_\R 
	hw'w\varphi^2.
\end{equation}
On the one hand, Cauchy's inequality implies
\[\int_\R w' w\varphi(h\varphi-2\varphi')\leq \frac12\int_\R 
(w')^2\varphi^2+\frac12\int_\R w^2(|h|\varphi+2|\varphi'|)^2.\]
On the other hand, we see that 
\[\int_\R g w^2\varphi^2 +  \frac12\int_\R w^2(|h|\varphi+2|\varphi'|)^2\leq 
C_\vp L^2 \norm{w}_{L^\infty(\tilde a,\tilde b)}^2,
\]
where $C_\vp$ is a constant depending only on the norms of $\vp$ in $(\tilde a,\tilde b)$.
Plugging the last two inequalities in \eqref{eq:H1est}, we infer 
\[
\frac 12\int_{a}^b (w')^2\leq\frac12\int_\R(w')^2\varphi^2\leq C_\vp L^2 \norm{w}_{L^\infty(\tilde a,\tilde b)}^2,\]
 which shows the bound for the first term in \eqref{est:der:w}, since $\norm{w}_{ L^\infty(\tilde a,\tilde b)}=w(\tilde b)$.
To treat the nonlocal term, we multiply \eqref{eq:w} by 
the cut-off function $\varphi$ and  integrate by parts, so that
\begin{equation}\label{eq:West}
2\int_a^b \W*w\leq 2\int_\R(\W*w)\varphi=\int_\R 
w(\varphi''+g\varphi-(h\varphi)')\leq 
\tilde C_\vp L { w(\tilde b)},
\end{equation}
where the constant $\tilde C_\vp$ depends only on the norms of $\vp$.
This completes the proof of \eqref{est:der:w}.
\end{proof}

The following lemma assures that, if a function satisfies a sort of Harnack inequality, then it can be {\em uniformly} bounded from below by an exponential function.

\begin{lemma}\label{lemma:harnack}
Assume that $w\in\boC^1((-\infty,0])$ satisfies
$w(0)=1$,  and $w'\geq 0$ in $(-\infty,0]$, and 
\begin{equation}
\label{ineq:beta}
w(x)\leq \beta w(x-1),\quad\text{for all }x\leq 0,
\end{equation}
for some $\beta >1$. 
Then there exist $A,\gamma>0$, depending on $\beta$, but not on $w$, such that 
$w(x)\geq Ae^{\gamma x}$, for all $x\leq 0$.
\end{lemma}

\begin{proof}
Take $\gamma=\log(\beta)$ and $A<1/\beta$, and let $w\in\boC^1((-\infty,0])$. 
If $x\in [-1,0]$, we deduce from  \eqref{ineq:beta} that
\[w(x)\geq w(-1)\geq w(0)/\beta \geq A\geq Ae^{\gamma x}.\]
If $x\in [-2,-1]$, we can apply this equality with $x+1\in[-1,0]$ to obtain
\[w(x)\geq \frac{1}{\beta}w(x+1)\geq \frac{A}{\beta}e^{\gamma(x+1)}=Ae^{\gamma x}.\]
The conclusion follows by an iterative argument.
\end{proof}

We recall that the analyticity of solutions is established in Lemma~\ref{lem:decay}, but only for solitons with speed less than the Landau speed.
To avoid this technical problem,  we now prove a unique continuation result adapted to our setting.
\begin{lemma}
\label{eta:nonzero}
    Let $c\in [0,\sqrt{2})$. Assume that  $\W$ satisfies \ref{W:L1}. Consider $u\in\boE(\R)$  a  solution to \eqref{TWc}. Assume that $\eta\geq 0$ in $\R$,
    or that $\eta\leq 0$ in $\R$, where $\eta=1-|u|^2$. If there is an open nonempty interval $I\subset\R$ such that $\eta=0$ in $I$, then 
    $\eta=0$ in $\R$.
\end{lemma}

\begin{proof}
We give the proof in the case $\eta\geq 0$ in $\R$. Arguing by contradiction, assume $\eta$ is nontrivial, so that either $\sup I$ or $\inf I$ is finite. Let us assume that $b:=\sup I$ is finite, with the other case being analogous. Moreover, without loss of generality, we may consider that $b$ is in the closure of $I$, with $\eta(b)=0$
and $\eta>0$ in $(b,b+\varepsilon)$ for some $\varepsilon\in (0,\tau)$.

On the one hand, using that $\eta$ satisfies \eqref{eq:eta} in $I$, one deduces that $\W*\eta=0$ in $I$, so that $(\W*\eta)(b)=0$. On the other hand, using \ref{W:L1}, we derive
    \begin{align*}
        (\W*\eta)(b)&=\int_\R \W(y)\eta(b-y)dy \geq \int_{-\tau}^0 \W(y)\eta(b-y)dy
        \\
        & \geq \sigma\int_{-\tau}^0 \eta(b-y)dy = \sigma \int_{b}^{b+\tau}\eta(y)dy\geq \sigma \int_b^{b+\varepsilon}\eta(y)dy>0.
    \end{align*}
    This is a contradiction.
\end{proof}

To clarify the proof of Theorem~\ref{thm:nonmonotone}, we first prepare the setting that allows us to obtain rigorously the limit equation. For this purpose, let us consider $u$ a  solution to \eqref{TWc},  with $\eta=1-|u|^2$ satisfying  $\eta(-\infty)=0$ and $\eta'(-\infty)=0$. To get a contradiction, let us also assume that there exists $x_0\in\R$ such that \eqref{osc:contradiction} holds. By translation invariance, we may assume without loss of generality that $x_0=0$.
Therefore,  to obtain a contradiction, we can invoke  Lemma~\ref{eta:nonzero} to define a sequence $\{x_n\}$ of positive real numbers such that $x_n\to\infty$,
with $\eta(x_n)>0$, and define the functions
\bq
\label{def:wn}
\eta_n(x)=\eta(x-x_n),\quad w_n(x)=\frac{\eta_n(x)}{\eta_n(0)},\quad g_n=c^2+\W*\eta_n+\frac{c^2 \eta_n}{2(1-\eta_n)},\quad 
	h_n=\frac{\eta_n'}{2(1-\eta_n)}.
 \eq
	One easily checks that $w_n$ satisfies
	\begin{equation}\label{eq:wn}
		-w_n''+2\W*w_n=g_n w_n + h_n w_n',\quad x\in (-\infty,x_n).
	\end{equation}
 In this manner, $\eta_n$ and $\eta_n'$ converge pointwise to $0$ as $n\to\infty$, and so does $\W*\eta_n$, 
 by the dominated convergence theorem.
Furthermore, the following properties are easily  verified:
\begin{enumerate}
\item \label{gnhn} 
 $\{ g_n \}$ and $\{ h_n \}$ are bounded in $W^{k,\infty}((-\infty,R))$ for every $k\in\N$ and every $R>0$, 
  with $g_n\to c^2$ and $h_n\to 0$ pointwise as $n\to\infty$, 
\item \label{wnpos} $w_n\geq 0$ in $\R$, $w_n\leq 1$ in $(-\infty,0)$, $w_n(0)=1$ and  $w_n'\geq 0$ in $(-\infty,x_n)$. 
	\end{enumerate}

If $\eta$ satisfies \eqref{osc:contradiction2}, we proceed in the same way. Noticing that in this case $\eta_n$
is negative, we have that $w_n$ is positive in $\R$ , with $w_n'\geq 0$ in $(-\infty,x_n)$. Moreover, the properties in (i) and (ii) stated above are still satisfied, so it remains to study \eqref{eq:wn}.

To simplify the analysis of \eqref{eq:wn}, we absorb the linear term $c^2w_n$ by defining the function $z_n(x)=w_n(x) e^{cx}$, so that
	\begin{equation}\label{eq:zn}
		z_n''=2(\W*w_n)e^{cx }+2cw_n'e^{cx }+f_n, \ \text{ in } (-\infty,x_n), \quad \text{ where } f_n= -((g_n-c^2) w_n + h_n w_n')e^{cx}.
	\end{equation}
By using the estimates in Lemma~\ref{lemma:estimateslinear}, we can pass to limit in \eqref{eq:wn} in a weak sense, to obtain the following result.

\begin{lemma}
\label{lem:osc:bis}
 Let $w_n$ be as defined above, satisfying \eqref{eq:wn}. Then there is $w\in W^{1,q}_\loc(\R)$, for all $q\in[1,\infty)$
 such that, for a subsequence, 
  $w_n\to w$ in $L_\loc^\infty(\R)$, and $w_n'\to w'$ in $L^q_\loc(\R)$ and 
  $w_n'\to w'$ pointwise.
 In addition,  $w$ is nonnegative and increasing in $\R$ ,  $w(0)=1$ and $w\leq 1$ in $(-\infty,0]$, and there are $A_0,\gamma_0>0$, such that 
\begin{equation}
\label{ineq:exponential}
	w(x)\geq A_0e^{\gamma_0 x},\quad\text{for all }x\leq 0.
\end{equation}
Moreover,   
  $w$ satisfies the differential inequality 
   \begin{equation}\label{ineqw} 
	(w(x)e^{cx })'(x)\geq 2\int_{-\infty}^x (\W*w(t))e^{ct
}dt,\quad\text{for all }x\in\R.
\end{equation}
\end{lemma}

\begin{proof}
 By using \ref{gnhn} and \ref{wnpos}, we can iterate  inequality \eqref{w:harnack} in Lemma~\ref{lemma:estimateslinear}, to deduce that  $\{w_n\}$ is bounded in $L^{\infty}_{\loc}(\R)$, in $H^1_\loc(\R)$ and also in $W^{2,1}_{\loc}(\R)$, by invoking \eqref{est:der:w}. For any bounded interval $I$, we can thus invoke the compact Sobolev embeddings,
  $W^{ 2,1 }(I) \hookrightarrow W^{1,q}(I)$,  for all $q\in[1,\infty)$, and 
  $H^{2}(I)\hookrightarrow \boC^0(I)$, and the Helly selection theorem to  conclude the existence of 
  $w\in W^{1,q}_\loc(\R)$ and the stated convergences.
Thanks to properties in \ref{gnhn} and \ref{wnpos}, it is immediate that $w_n$ satisfies \eqref{w:harnack},
so \eqref{ineq:beta} holds for $w_n$, for some constant $\beta>1$ independent of $n$. Thus, Lemma~5.2 provides the existence of $A$ and $\gamma$ such that $w_n(x)\geq A_0e^{\gamma_0 x}$, for all $x\leq 0$, and \eqref{ineq:exponential} follows  invoking the pointwise  convergence.
 
 To prove \eqref{ineqw}, we set $z(t)=w(t)e^{ct}$, for $t\in \R$,  fix $x\in\R$ and consider  a sequence $\{\varphi_k\}$ of nondecreasing $\boC^2$-functions,
 such that
$$\varphi_k:(-\infty,x]\to [0,1], \quad 
\varphi_k(t)=1\,\,\text{for all }t\in[x-k,x],\quad \varphi_k(t)=0\,\,\text{for all }t\leq x-2k.$$
Taking into account that $\varphi_k'(x)=0$, so that $z_n'(x)\varphi_k(x)=(z_n\varphi_k)'(x)$, and
using \eqref{eq:zn}, we compute
\begin{align}
z_n'(x)\varphi_k(x)=\int_{-\infty}^x (z_n\varphi_k)''(t)dt &=\int_{-\infty}^x( 2(\W*w_n)\varphi_k e^{ct}  + 2cw_n'\varphi_ke^{ct } +f_n\varphi_k +2z_n'\varphi_k '+z_n\varphi_k '' )dt
 \nonumber
 \\
 \label{dem:osc1}
 	&=\int_{-\infty}^x( 2(\W*w_n)\varphi_k e^{ct} + 2cw_n'\varphi_ke^{ct } +f_n\varphi_k +z_n'\varphi_k ' )dt,
  \end{align}
where the last equality follows by integration by parts.

Now we take inferior limits (with respect to $n$) in this equality, and we bound the last term from below. In the first place, Fatou's lemma (applied twice) implies
\bq\label{dem:osc2}
\begin{aligned}
	\liminf_{n\to\infty} 
	\int_{-\infty}^x(\W*w_n)\varphi_k &e^{ct}dt \geq\int_{-\infty}^x\liminf_{n\to\infty}(\W*w_n)\varphi_k e^{ct}dt 
 \\
	&= \int_{-\infty}^x\liminf_{n\to\infty}\left(\int_\R\W(t-y)w_n(y)dy\right)\varphi_k(t) e^{ct}dt
	\\
 	&\geq \int_{-\infty}^x\left(\int_\R\W(t-y)w(y)dy\right)\varphi_k(t) e^{ct}dt
	=\int_{-\infty}^x(\W*w)\varphi_k e^{ct}dt.
\end{aligned}
\eq
In the second place, setting $I=(x-2k,x)$, the  convergence of $\{w_n'\}$ 
to $w'$ in $L^2(I)$ implies that, for a subsequence, $\{w_n'\}$ is dominated by some  function in $L^2(I)$, while 
 $w_n\to w$  in $L^\infty(I)$.
Using also that, $f_n\to 0$  pointwise as $n\to\infty$ and the estimates in \ref{gnhn},
the dominated convergence theorem implies that 
$$\int_{-\infty}^x(2cw_n'\varphi_k e^{ct }+f_n\varphi_k +z_n'\varphi_k ' )dt \to \int_{-\infty}^x (2cw'\varphi_k e^{ct }+ z'\varphi_k' )dt ,
\quad\text{as }n\to\infty.$$
 Combining with \eqref{dem:osc1} and \eqref{dem:osc2}, we end up with
\[z'(x)\varphi_k(x)\geq \int_{-\infty}^x(2 (\W*w)\varphi_k+w'\varphi_k' + 2cw'\varphi_k + cw\varphi_k')e^{ct}dt.\]
Since $w'\geq 0$ and $\varphi'_k\geq 0$,  we can apply Fatou lemma, this time in $k$, to conclude that 
\begin{equation*}
	z'(x)\geq \int_{-\infty}^x (2(\W*w)(t) + 2cw'(t))e^{ct}dt.
\end{equation*}
Dropping the last (positive) term on the right-hand side, we finish the proof  \eqref{ineqw}.
\end{proof}

\begin{remark}
In the case of the soliton $\boldeta(\cdot \, ; c)$   in \eqref{sol:1D}, the function $w_n$ in \eqref{def:wn}
can be computed explicitly, and we can check that $w_n$ converges pointwise to $ w(x)=\exp(\sqrt{2-c^2}x)$. 
Therefore, we expect, in general,  the function $w$ in Lemma~\ref{lem:osc:bis} to be unbounded.
\end{remark}

We can now complete the proof of our main theorem concerning the oscillating behavior of solutions to \eqref{TWc}.
\begin{proof}[Proof of Theorem \ref{thm:nonmonotone}]
Assume by contradiction that there exists $x_0\in\R$ such that \eqref{osc:contradiction} or \eqref{osc:contradiction2} holds. Inequality \eqref{ineq:exponential} in  Lemma~\ref{lem:osc:bis} enables us to define
\[\bar\gamma=\inf\{\gamma>0:\,\,\text{there exists }A>0\text{ such  that }w(x)\geq Ae^{\gamma 
x}\text{ for all }x\leq 0\}.\]
Indeed, by \eqref{ineq:exponential}, the set where the infimum is taken is nonempty. Let $\gamma\in(\bar\gamma,\bar\gamma+\delta/2)$, where $\delta>0$ will be chosen later. Let $A>0$ 
be such that \eqref{ineq:exponential} holds for $\gamma_0=\gamma$ and $A_0=A$. Hence, for any $R>0$ and any $x\leq -R$, it follows 
\[(\W\ast w)(x)\geq \int_{-R}^\infty\W(y)w(x-y)dy\geq Ae^{\gamma x}\int_{-R}^\infty \W(y)e^{-\gamma 
y}dy= Ae^{\gamma x}\wc\W (\gamma,R).\]
Plugging this into \eqref{ineqw} and integrating the resulting inequality, we derive
\bq
\label{dem:them:osc}
w(x)e^{cx}\geq \frac{2A}{(\gamma+c)^2}e^{(\gamma+c) x}\wc\W(\gamma,R),\quad\text{for all }x\leq -R.
\eq
By \eqref{H:bumps}, we have $$\lim_{R\to\infty}\frac{2\wc\W(\gamma,R)}{(\gamma+c)^2}>1,$$ so that there is $R_\gamma>0$ such that 
\[\delta:=\frac{1}{R_\gamma}\log\left(\frac{2\wc\W(\gamma,R_\gamma)}{(\gamma+c)^2}\right)>0.\]
It thus follows from \eqref{dem:them:osc} that
\[w(x)\geq A e^{\delta R_\gamma+\gamma x},\quad \text{for all }x\leq -R_\gamma.\]
Iterating the previous argument, we deduce 
\begin{equation}\label{expcontradiction}
	w(x)\geq A e^{\delta NR_\gamma+\gamma x},\quad \text{for all }N\in\N\cup \{0\},\,\,x\leq -NR_\gamma.
\end{equation}

Now let $y\leq 0$ and let $N_y\in\N$ be such that
$-(N_y+1)R_\gamma<y\leq -N_yR_\gamma.$
Using \eqref{expcontradiction} with $N=N_y+1$, and the fact that $w$ is nondecreasing, we get
\[w(y)\geq w(-(N_y+1)R_\gamma)\geq A e^{-(\gamma-\delta)(N_y+1)R_\gamma}\geq 
Ae^{-R_\gamma(\gamma-\delta)}e^{(\gamma-\delta)y}.\] 
Observe that the constant $\tilde A=Ae^{-R_\gamma(\gamma-\delta)}$ does not depend on $y$, so the previous inequality is valid for every $y\leq 0$. Therefore, the fact that $\bar\gamma$ is an infimum implies that 
$\gamma-\delta\geq \bar\gamma.$
This is a contradiction with the choice of $\gamma$. We conclude that \eqref{osc:contradiction} cannot happen.

Finally, let us show that $\eta$ has at least one oscillation. First, assume that $\eta$ changes its sign, i.e.\
there are points $y_1,y_2\in \R$, with $y_1<y_2$, such that $\eta(y_1)\eta(y_2)<0$. In the case  that $\eta(y_1)>0$, the 
boundary condition $\eta(\pm\infty)=0$ implies that there are $x_1<y_1<x_2<y_2<x_3$ such that $\eta'(x_1)>0$, $\eta'(x_2)<0$ and $\eta'(x_3)>0$.
The case  $\eta(y_1)<0$ is analogous. 

Assume now that $\eta\geq 0$ in $\R$. Since $\eta$ is nontrivial and $\eta(\pm\infty)=0$, we deduce that  $\eta$ admits a positive global maximum, 
that we suppose to be attained at the origin, for simplicity. Thus, there is $x_3<0$ such that $\eta'(x_3)>0$. In view of  \eqref{osc:contradiction}, there must exist 
$x_2\in (-\infty,x_3)$ such that $\eta'(x_2)<0$. This also implies that $\eta(x_2)>0$. Hence, the condition  $\eta(-\infty)=0$ provides the existence of $x_3\in (-\infty,x_2)$ such that $\eta'(x_3)>0.$ The case $\eta\leq 0$ in $\R$ is analogous. 
 \end{proof}

It is clear that  Theorem~\ref{thm:nonmonotone} applies for instance to the potentials \eqref{pot:Gauss} and  \eqref{pot:rect}, 
since the Laplace transform are given by 
$\wc \boW(s)=\sinh(\abs{\lambda}s)/(\abs{\lambda}s) $ and $\wc \boW(s)=e^{\lambda^2s^2}$, $s\in \R$,
respectively. In addition, we remark that Theorem~\ref{thm:nonmonotone} can be generalized and improved by using the properties of $w_n$ in \eqref{def:wn} for some specific potentials. We end this section by detailing this point for the potentials in  \eqref{pot:nematic} and  \eqref{pot:vander}.

Let  $\W_\lambda$ be the nematic potential in \eqref{pot:nematic}, with $\lambda>0$. Using its  Laplace transform given in \eqref{exponentialexample}, condition  \eqref{H:conjecture} reads:
\bq
\label{laplace1}
\frac{P_{\lambda,c}(s)}{\lambda^2s^2-1}<0, \quad \text { for all } s\in (0,1/{\lambda}), \text{  where }\ 
P_{\lambda,c}(s)=\lambda^2s^4+(c^2\lambda^2-1)s^2+2-c^2.
\eq
Since $P_{\lambda,c}\geq 2$ in $[1/{\lambda},\infty)$, we conclude that the condition in \eqref{laplace1} is equivalent to saying 
that $P_{\lambda,c}$ is positive in $\R$. Thus,  we need  the discriminant to be  negative (see Remark~\ref{rem:polynomials}), i.e.
\bq
\label{laplace2}
(c^2\lambda^2-1)^2-4\lambda^2(2-c^2)<0.
\eq
Finding the roots of this polynomial, we obtain the value $\tilde \lambda_c$ in \eqref{exponential:lambdactilde}, and 
 \eqref{laplace2} holds exactly for all $\lambda>\tilde \lambda_c.$

In the case of the potential in \eqref{pot:vander} with $\lambda<0$, we need to use the definition of the Laplace transform of a  distribution to 
justify that  $\wc \delta_0=1$, and argue as before to conclude that  \eqref{H:conjecture} is equivalent to saying that the polynomial 
\bq
\label{pol:cha}
  P_{\lambda,\beta,c}(s)=s^4 - (\beta^2+2A_\lambda-c^2)s^2 + \beta^2(2A_\lambda-c^2)+4\beta A_\lambda\lambda, \
\text{ with } A_\lambda=\beta/(\beta-2\lambda), 
\eq
is positive in $\R$. This reduces again to the study of the discriminant, but the analysis in this case is more involved, as seen in the proof of the next result.

\begin{theorem}
\label{thm:local:potentials}
Let $c\in [0,\sqrt 2)$. Consider one of the following cases: 
\begin{enumerate}
    \item $\W_\lambda$ is  the potential in \eqref{pot:nematic} and   $\tilde \lambda$ is defined in \eqref{exponential:lambdactilde}, and  $\lambda>\tilde\lambda_c$.
    \item  $\W_\lambda$ is  the potential in  \eqref{pot:vander}, with $\beta>0$ and $\lambda<0$  such that the number in \eqref{mon:waals} is either negative or 
    has nonzero imaginary part.
\end{enumerate}
If  $u\in\boE(\R)$ is  a nontrivial solution to $\TWclambda$, then it cannot exist $x_0\in\R$ such that either \eqref{osc:contradiction}  or \eqref{osc:contradiction2} hold, where $\eta=1-\abs{u}^2$. 
In particular, $\eta$ presents at least one oscillation (in the sense of Theorem~\ref{thm:nonmonotone}).
In addition, if $c=0$ and  $u(\pm \infty)=\pm 1$, then 
  $u$ cannot be nondecreasing in $\R$.
\end{theorem}

  \begin{proof}
Let $\W_\lambda=e^{-\abs{\cdot }/\lambda}/(2\lambda)$. We start as in the proof of Theorem~\ref{thm:nonmonotone}, arguing by contradiction and defining   
$\eta_n$, $w_n$, $g_n$ and $f_n$ as in \eqref{def:wn}. In addition, we set $v_n=\W_\lambda*w_n.$
Invoking  Lemma~\ref{lemma:estimateslinear}, we deduce that for any $R>0$, there is  $C>0$ such that 
\begin{equation}\label{estimateswnexample}
  \|w_n\|_{L^\infty(I)} +\|w_n\|_{W^{2,1}(I)}+\|w_n'\|_{L^2(I)}+\|v_n\|_{L^1(I)}\leq C,
 \end{equation}
for all $n$ large enough, where $I$ is a compact interval of $(-\infty,R)$. In addition, the choice of the potential allows us to recast   \eqref{eq:wn} as the system
\begin{equation}\label{localsystem}
-w_n''+2v_n=g_nw_n+h_nw_n', \quad v_n -\lambda^2 v_n''=w_n,\quad  \text{ in } (-\infty,R).
\end{equation} 
 From the estimates in \eqref{estimateswnexample}, we infer from the second equation in \eqref{localsystem}  
 that $\{v_n''\}$ is also bounded in $L^1(I)$. Thus, by interpolating $\norm{v_n'}_{L^1(I)}$, we conclude that $v_n$ is bounded in $W^{2,1}(I)$. As in  Theorem~\ref{thm:nonmonotone}, this enables us to deduce the existence of $v,w\in W^{1,2}_\loc(\R)$ such that $v_n\to v$ and $w_n\to w$ in $L_\loc^\infty(\R)$, and $v'_n\to v'$ and  $w'_n\to w'$ in $L_\loc^2(\R)$. 
 Then we can pass to the limit  in the weak formulation of \eqref{localsystem}, which results in
\begin{equation}\label{limitsystem}
-w''+2v-c^2w=0, \quad v -\lambda^2 v''-w=0, \quad \text{ in }\mathcal{\boD'(\R)}.
\end{equation} 
By elliptic regularity, we infer that $w,v\in H^2_\loc(\R)$, which implies that $(w,v)$  is a classical solution to the linear system of ODEs in \eqref{limitsystem}. Moreover, since  $w_n$ is nondecreasing with $w_n(0)=1$, it follows that $w$ is a nonzero nondecreasing function.

On the other hand, we can check that the characteristic polynomial of the linear system in \eqref{limitsystem} is given by $P_{\lambda,c}$
in \eqref{laplace1}. Thus, $P_{\lambda,c}$  has non-real roots for $\lambda>\tilde \lambda_c$. However, this implies that $w$ cannot be monotone, 
which is a contradiction.

The proof for the potential $\W_\lambda$  in  \eqref{pot:vander} is very similar,
since Lemma~\ref{eta:nonzero} remains true in this case, due to the fact that $\lambda<0$.
Indeed, setting $A_\lambda=\beta/(\beta-2\lambda)$, $\mu_\lambda=-\lambda e^{-\beta|\cdot|}$, we can write \eqref{eq:wn}  as
    \begin{equation}
     \label{eq:vander:dem}
		-w_n''+2A_\lambda\mu_\lambda*w_n=\tilde g_n w_n + h_n w_n', \  \text{ with }
  \tilde g_n=-2A_\lambda+c^2+\W_\lambda*\eta_n+\frac{c^2 \eta_n}{2(1-\eta_n)}, 
\end{equation}
 in $(-\infty,R)$, so that  $\tilde g_n\to c^2-2A_\lambda$ pointwise. 
Since $\lambda<0$, $\mu_\lambda$ satisfies \ref{W:L1}. Thus, we can invoke  Lemma~\ref{lemma:estimateslinear} as before to obtain the estimates in \eqref{estimateswnexample}, 
with  $v_n=\mu_\lambda*w_n$, and also  recast \eqref{eq:vander:dem} as:
\begin{equation*}
-w_n''+2A_\lambda v_n=\tilde g_n w_n+h_nw_n', 
\quad 
\beta^2 v_n - v_n'' = -2\beta\lambda w_n, \quad \text{ in } (-\infty,R).
\end{equation*} 
Arguing as before, we conclude that there are  nonzero functions $v,w$, with $w$ nondecreasing,  solutions to the  linear system 
\begin{equation}\label{limitsystem2}
-w''+ 2A_\lambda v + (A_\lambda - c^2) w = 0, \quad \beta^2 v - v''- 2\beta\lambda w =0.
\end{equation} 
The characteristic polynomial of the system is given by $ P_{\lambda,\beta, c}$
in \eqref{pol:cha}, with discriminant 
$$\triangle=
(\beta^2+2A_{\lambda}-c^2)^2
-4\beta^2(2A_{\lambda}-c^2)- 16\beta A_{\lambda} \lambda.
$$
Therefore,  the number in  \eqref{mon:waals} corresponds to the square of a root of 
$ P_{\lambda,\beta, c}$. Thus, if this number is negative or has a nonzero imaginary part, then $ P_{\lambda,\beta, c}$ has roots with nonzero imaginary parts, which contradicts the monotonicity of $w$.
   \end{proof}

\begin{remark}\label{rem:polynomials}
 The condition $\lambda>\tilde \lambda_c$ for the nematic potential cannot be improved as done for the potential in \eqref{pot:vander}. Indeed,  some algebraic computations show that if the discriminant of $P_{\lambda,c}$ is positive, then 
 $P_{\lambda,c}$ cannot have complex solutions with nonzero imaginary part, for $c\in[0,\sqrt 2).$  For this reason, we could expect the conditions in Theorem~\ref{thm:local:potentials} to be sharp to detect oscillations of dark solitons.
   \end{remark}


\section{Analysis of black solitons}\label{sec:blacksolitons}

Recall that by Lemma~\ref{lemma:regularity} black solitons must be real-valued.
We are now ready to prove the existence result for black solitons.

\begin{proof}[Proof of Theorem~\ref{thm:existenceblack}]
	Since $\wh\W\geq 0$, it follows from Plancherel's identity in \eqref{Plancherel} that $E(u)\geq 0$ for every $u\in\Eodd$. Then, 
 there exists a sequence $\{u_n\}\subset\Eodd$ such that $E(u_n)\to\inf\{E(v):\,\,v\in\Eodd\}$, as $n\to\infty$. In particular, $\{E(u_n)\}$ is bounded. As a consequence, $\{u_n'\}$ is bounded in $L^2(\R)$ and, by Lemma~{2.1} of \cite{delaire-mennuni}, $\{\eta_n\}$ is bounded in $L^2(\R)$ and in $L^\infty(\R)$. Therefore, passing to a subsequence, there exists $u\in\Eodd$ such that, as $n\to\infty$,
\begin{align*}
	u_n\to u,\ \text{in }L^\infty_{\mbox{\tiny loc}}(\R),
	\quad 
	u_n'\wto u',\ \text{in }L^2(\R),
	\quad 
	\eta_n\wto \eta,\ \text{in }L^2(\R),
\end{align*}
by the Rellich-Kondrakov theorem.  Furthermore,
\bq
\label{dem:c0:1}
\|u'\|_{L^2(\R)}\leq \liminf_{n\to\infty}\|u_n'\|_{L^2(\R)}.
\eq
In addition, since $\wh\W\geq 0$ in $\R$, the functional $f\mapsto\int_\R(\W*f)f$ is convex in $L^2(\R)$ and,  in turn,
\bq
\label{dem:c0:2}
\int_\R\big(\W*\eta\big)\eta\leq \liminf_{n\to\infty}\int_\R\big(\W*\eta_n\big)\eta_n.
\eq
From \eqref{dem:c0:1} and \eqref{dem:c0:2}, it is now clear that
\[E(u)\leq\liminf_{n\to\infty} E(u_n)= \inf\{E(v):\,\, v\in\Eodd\}\leq E(u).\]
Therefore, $E(u)=\min\{E(v):\,\, v\in\Eodd\}$. 

It remains to show that $u$ is a solution to \eqref{TW0}. Using the fact that $u$ is a minimizer among odd functions, and bearing in mind the commutative property in \eqref{Plancherel}, it is standard to show that
	\[\int_\R u'\varphi' = \int_\R (\W*(1-u^2))u\varphi, \quad \text{ for every odd function }\varphi\in\boC^1_c(\R).
 \]
Furthermore, since $\W$, $1-u^2$ and $u'$ are even,  it follows that 
	\[\int_\R u'\varphi' = 0 = \int_\R (\W*(1-u^2))u\varphi, \quad \text{ for every even function }\varphi\in\boC^1_c(\R).\]
Hence, we conclude that $u$ is a solution to \eqref{TW0} by recalling that every function $\varphi\in\boC^1_c(\R)$ can be decomposed as
	\[\varphi(x)=\frac{\varphi(x)+\varphi(-x)}{2}+\frac{\varphi(x)-\varphi(-x)}{2}\coloneqq\varphi_{\mbox{\tiny even}}(x)+\varphi_{\mbox{\tiny odd}}(x),\]
	where both $\varphi_{\mbox{\tiny even}}$ and $\varphi_{\mbox{\tiny odd}}$ obviously belong to $\boC^1_c(\R)$.	
\end{proof}

We finally include the proof of the nonlocal-to-local limit theorem in the case $c=0$.

\begin{proof}[Proof of Theorem~\ref{thm:limittolocalblack}]
Since $u_\lambda$ is a minimizer, applying Plancherel's identity, we derive
\[E_\lambda(u_\lambda)\leq E_\lambda(v)\leq (M+1)E_0(v),
\quad \text {where } M=\sup\{ \|\wh\W_\lambda\|_{L^\infty(\R)} : \lambda\in (0,\lambdazero)\},
\]
for any $v\in\Eodd$.  As a consequence, $\{u_\lambda'\}_{\lambda\in (0,\lambdazero)}$ is bounded in $L^2(\R)$ and, by Lemma~{2.1} of \cite{delaire-mennuni}, $\{\eta_\lambda\}_{\lambda\in(0,\lambdazero)}$ is bounded in $L^2(\R)$ and in $L^\infty(\R)$. Passing to a subsequence, there exists $u_0\in\Eodd$ such that, as $\lambda\to 0$,
\begin{align*}
	u_\lambda\to u_0,\text{ in }L^\infty_{\mbox{\tiny loc}}(\R),
	\quad 
	u_\lambda'\wto u_0',\ \text{in }L^2(\R),
	\quad 
	\eta_\lambda\wto \eta_0,\ \text{in }L^2(\R).
\end{align*}
Finally, we apply Lemma~\ref{lemma:limitnonlocalterm} to pass to the limit in the weak formulation of the equation 
\[-u_\lambda''=u_\lambda\big(\W_\lambda*(1-u_\lambda^2)\big),\quad\text{in }\R,\]
and conclude that $u_0$ is a solution to $\TWzdelta$. By uniqueness of nontrivial finite-energy solution for equation  $\TWzdelta$, necessarily there exist $t\in\R$ and $\theta\in\R$ such that \eqref{eq:limittolocal} holds with $c=0$. Moreover, since $u_0$ is real-valued, $\theta=0$, and since it is odd, $x_0=0$, so $u_0=\bu(\cdot\,;0)$.
\end{proof}


\subsection{Monotone black solitons}

This subsection is devoted to the proof of Theorem~\ref{thm:monotonicity}. It is based on elementary theory on \emph{symmetric decreasing rearrangements.} In order to be self-contained, let us define them and state some basic properties that we will use. More details can be found, for instance, in the classical book \cite{LiebLoss2001}  or in the notes \cite{Burchard2009}.

\begin{definition} 
	On the one hand, given a measurable function $v:\R\to\R$, we say that it \emph{vanishes at infinity} if $|\{|v|>t\}|<\infty$ for all $t>0$. On the other hand, if $v$ vanishes at infinity, its \emph{symmetric decreasing rearrangement} is defined as
\[v^\star(x)=\int_0^\infty\chi_{\{|v|>t\}^\star}(x)dt,\]
where $\{|v|>t\}^\star$ is the open interval centered at zero whose measure coincides with $|\{|v|>t\}|$.
\end{definition}

One deduces straightaway from the definition that $v^\star$ is a nonnegative, even and lower semicontinuous function which is also nonincreasing in $(0,+\infty)$. Furthermore, one of the fundamental properties of the symmetric decreasing rearrangements is that, for every $t>0$,  the following identity holds,
\begin{equation}\label{sdr:fundamentalprop}
	\{v^\star>t\}=\{|v|>t\}^\star,
\end{equation}
which in turn implies
\[|\{v^\star>t\}|=|\{|v|>t\}|.\]
As a consequence, one has the following relation.

\begin{proposition}\label{sdr:prop1}
Let $G:[0,+\infty)\to\R$ be a function such that $G=G_1-G_2$, where $G_1$ and $G_2$ are both nondecreasing with $G_1(0)=G_2(0)=0$ and, for at least one $j\in\{1,2\}$, one has that $\int_\R G_j(|v(x)|)dx$ is finite. Then, for every measurable function $v:\R\to\R$ that vanishes at infinity, one has
\begin{equation}\label{sdr:integralidentity}
	\int_\R G(|v(x)|)dx=\int_\R G(v^\star(x))dx.
\end{equation}
\end{proposition}

Let us state and prove a further property that is still elementary but not usually listed in the classical references.

\begin{proposition}\label{sdr:prop2}
Let $v$ be an even function, nonincreasing in $(0,+\infty)$, vanishing at infinity. Then $v^\star\leq v$ in $\R$.
\end{proposition}

\begin{proof}
The fact that $v$ is even and nonincreasing in $(0,+\infty)$ implies that $\{v>t\}$ is an interval centered at zero. Therefore,  $\{v>t\}^\star\subset\{v>t\}$ and thus
\[v^\star(x)\leq \int_0^\infty \chi_{\{v>t\}}(x)dt=\int_0^{v(x)}dt=v(x).
\qedhere  
\]
\end{proof}

We will also need a generalization of the classical Riesz's rearrangement inequality. The statement and proof can be found in \cite[Theorem 2.2]{AlmgrenLieb1989}. We introduce here the statement adapted to our setting.

\begin{theorem}\label{thm:AlmgrenLieb}
	Let $G:[0,+\infty)\to[0,+\infty)$ be a function of class $\boC^1$ such that
	\[G(0)=0,\quad\text{and}\quad G'(s)\geq 0,\,\,\text{for all }s\geq 0.\]
	Let $\sigma\in L^1(\R)$ with $\sigma\geq 0$ in $\R$. Then, for every measurable function $v:\R\to\R$ that vanishes at infinity, one has
	\[\int_\R\int_\R G(|v(x)|)G(|v(y)|)\sigma(x-y)dydx\leq \int_\R\int_\R G(v^\star(x)) G(v^\star(y))\sigma^\star(x-y)dydx.\]
\end{theorem}

We finally recall the well-known P\'olya-Szeg\H{o} inequality.

\begin{theorem}\label{thm:PolyaSzego}
If $v\in W^{1,p}(\R)$ for some $p\in [1,\infty)$, then $v^\star\in W^{1,p}(\R)$ and
\[\int_\R|(v^\star)'(x)|^p dx\leq \int_\R|v'(x)|^p dx.\]
\end{theorem}

\begin{proof}[Proof of Theorem~\ref{thm:monotonicity}]
	Let $u\in\boC^2(\R)$ be a solution to \eqref{TW0} such that $u(\pm\infty)=\pm 1$. Since $\mu^+=0$, it follows from \eqref{est:deltamu} that $\|u\|_{L^\infty(\R)}\leq 1$, i.e. $-1\leq u\leq 1$ in $\R$.  It remains to show the strict inequality \eqref{-1u1}. In order to do so, assume by contradiction that there exists $x_0\in\R$ such that $u(x_0)=-1$ (the case $u(x_0)=1$ is analogous). Clearly, $u(x_0)=\min_{x\in\R}u(x)$. Thus,
	\begin{align*}
		0\geq - u''(x_0)=u(x_0)\big(\W*(1-u^2)\big)(x_0)=-A_\mu \big(\mu*(1-u^2)\big)(x_0)\geq 0.
	\end{align*}
	As a result, 
	\[0=\big(\mu*(1-u^2)\big)(x_0)=\int_\R \mu(x-x_0)(1-u(x)^2)dx.\]
	Finally, since both $\mu$ and $1-u^2$ have a sign, and $\mu$ is nontrivial, we conclude that $u^2\equiv 1$, i.e. $u\equiv -1$, a contradiction with $\lim_{x\to\infty}u(x)=1$.

We assume now that $\mu$ is nondecreasing in $(0,+\infty)$. We already know by Theorem~\ref{thm:existenceblack} that there exists $\tilde{u}\in\Eodd$ such that $E(\tilde{u})=\min\{E(u):\, u\in\Eodd\}$. Let us show that there exists a possibly different minimizer which is increasing.

At this point, we introduce the auxiliary functional 
\[J(v)=\frac{1}{2}\int_\R (v')^2 + \frac14\int_\R (\W*G(|v|))G(|v|),\quad v\in H^1(\R),\]
where $G:[0,+\infty)\to\R$ is of the form $G=G_1-G_2$, where
\[G_1(s)=2s,\quad G_2(x)=\min\{s^2,1\},\quad s\in\R.\]
It is easy to see that the function $v=1-|\tilde u|$ satisfies $v\in H^1(\R)$ and $J(v)=E(\tilde u)$. Indeed, on the one hand, we have
	\[||\tilde u|'|=|\tilde u'|\text{ if }\tilde u\not=0,\quad |\tilde u|'=0\text{ otherwise.}\]
Therefore, $v'\in L^2(\R)$ and $\int_\R (\tilde u')^2=\int_\R(v')^2$. On the other hand, the identity $(1-|\tilde u|)^2=(1-\tilde u^2)^2(1+|\tilde u|)^{-2}$ implies that $v\in L^2(\R)$. Hence we deduce that $v\in H^1(\R)$. Finally, since $|\tilde u|<1$ in $\R$, it follows that $0<v<1$ in $\R$, so $G(|v|)=G(v)=v(2-v)=1-\tilde u^2$. Thus the identity  $E(\tilde u)=J(v)$ trivially holds.
	
	Now, by virtue of Proposition~\ref{sdr:prop2}, $(-\mu)^\star\leq - \mu$ in $\R$. Therefore, Theorem~\ref{thm:AlmgrenLieb} can be applied with $\sigma=-\mu\geq 0$ to deduce
	\begin{equation}\label{eq:Riesz}
		\int_\R(\mu*G(v))G(v)\geq \int_\R(\mu*G(v^\star))G(v^\star).
	\end{equation}
	Next, Theorem~\ref{thm:PolyaSzego} and Theorem~\ref{sdr:prop1}, combined with \eqref{eq:Riesz}, lead to the inequality $J(v)\geq J(v^\star)$.
	
	Let us now define
	\[u_0(x)=\left\{\begin{array}{ll} 1-v^\star(x) &\text{ if } x\geq 0,\\ -1+v^\star(x) &\text{ if }x<0.\end{array}\right.\]
	Let us check that $u_0$ is continuous at zero. Indeed, on the one hand, we recall that $\tilde u(0)=0$, so that $\max_{x\in\R} v(x)=v(0)=1$. On the other hand, it is clear that $v^\star$ achieves its maximum at zero and, from \eqref{sdr:fundamentalprop}, we derive  $\max_{x\in\R}v^\star(x)=\max_{x\in\R}v(x)$. In sum, $v^\star(0)=1$, so $u_0$ is continuous at zero.
	
It is clear now that $u_0$ is a nondecreasing odd function that satisfies $1-u_0^2=G(v^\star)$. Moreover, taking into account the continuity of $u_0$ at zero and the fact that $v^\star\in H^1(\R)$, it is a simple exercise to prove that $u_0\in H^1_{\mbox{\tiny loc}}(\R)$ with  $u_0'=|(v^\star)'|$. Therefore, $u_0\in\Eodd$ and
	\[E(u_0)=J(v^\star)\leq J(v)=E(\tilde{u}).\]
	In conclusion, $E(u_0)=E(\tilde{u})=\min\{E(u):\, u\in\Eodd\}$.
\end{proof}

\section{Proofs of the examples}\label{sec:proofsexamples}
Let us start by studying the convergences of the potentials in the examples given in Subsection~\ref{subsec:examples},
  in the strong  and   vague (or  weak-$*$) topology of $\boM(\R)$. 
The analysis is contained 
in the following result.   
\begin{lemma}
\label{lem:Wlam}
 Let $\W_\lambda$ be given by either \eqref{pot:Gauss}, \eqref{pot:nematic}, \eqref{pot:vander}, 
 \eqref{pot:rect} or \eqref{pot:deltas}.  In all the cases, we set $\lambdazero>0$ arbitrary, except for $\eqref{pot:vander}$ where we fix $\lambdazero$ in the interval $(0,\beta/2)$. Then, the family $\{\boW_\lambda \}_{\lambda\in (-\lambdazero,\lambdazero)}$ satisfies   \ref{H0}, 
 \ref{W:derivative}, \ref{Wlambda:aeconvergence} and \ref{Wlambda:finitesup}. In particular, $\boW_\lambda$ converges vaguely to $\delta_0$, i.e.
 in the  weak-$*$ sense: 
 \bq
 \label{conv:vague}
 \int_\R \vp d\W_\lambda \to \vp(0), \text{ for all  }\vp\in \boC_0(\R), \text{ as }\lambda\to 0.
\eq
Finally, in the case \eqref{pot:vander}, $\boW_\lambda$ converges to $\delta_0$ in $\boM(\R)$, as $\lambda\to 0$. In all the other cases, 
$\boW_\lambda$ {\em does not} converge  in $\boM(\R)$.
\end{lemma}

\begin{proof}
It is trivial that the potentials satisfy \ref{H0} and  \ref{W:derivative}. 
We consider now the cases \eqref{pot:Gauss}, \eqref{pot:nematic} and \eqref{pot:rect}, where 
 the potentials are given by nonnegative functions. Hence, the family $\{\boW_\lambda \}_{\lambda_\in (-\lambdazero,\lambdazero)}$ is bounded, since
$$\norm{\boW_\lambda}_{\boM(\R)}=\norm{\boW_\lambda}_{L^1(\R)}=1,$$
for all $\lambda$.
In view of the explicit Fourier transform of these potentials, it is clear that  \ref{Wlambda:aeconvergence} is satisfied with $\alpha=0$.
Thus, by invoking  Proposition~{8.50} in \cite{folland}, we infer that $\boW_\lambda$ converges vaguely to $\delta_0$. 
Arguing by contradiction, let us assume that $\boW_\lambda$ converges strongly in $\boM(\R)$, which implies that it converges weakly-$*$. Thus, we deduce that $\boW_\lambda$ must converge to $\delta_0$ in $\boM(\R)$, so that,
by \eqref{W:cota:medida}, $\wh \boW_\lambda$ must converge to 1 in $L^\infty(\R)$. This is a contradiction, since the decay at infinity of the potential implies that $\norm{\wh \boW_\lambda-1}_{L^\infty(\R)}=1.$

In the case \eqref{pot:vander}, we use the triangle inequality to get 
$$\norm{\boW_\lambda-\delta_0}_{\boM(\R)} \leq \frac{\abs{\lambda}}{\beta-2\lambda}(2+ \norm{e^{-\beta\abs{\cdot }}}_{L^1(\R)})\leq  \frac{\abs{\lambda}}{\beta-2\lambdazero}(2+2/\beta).$$
Thus, $\boW_\lambda$ converges to $\delta_0$ strongly, and thus vaguely, as $\lambda\to 0$.

Finally, consider $\boW_\lambda$ given by the potential in \eqref{pot:deltas}. Its Jordan decomposition is $\boW_\lambda^+=2\delta_0$ and $\boW_\lambda^-=(\delta_{\lambda}+\delta_{-\lambda})/2.$ 
Indeed, taking $B=\R\setminus \{0\}$, we have  $\boW_\lambda^-(B)=\boW_\lambda^+(B^c)=0$, so that $\boW_\lambda^+$ and $\boW_\lambda^-$ are 
 are mutually singular. Therefore,  \ref{Wlambda:aeconvergence}  holds with $\alpha=1$ and 
the vague convergence is deduced as before. To see that there is no strong convergence,  we compute 
\begin{equation*}
\norm{\boW_\lambda-\delta_0}_{\boM(\R)}=\norm{\delta_0}_{\boM(\R)}+\frac12 \norm{\delta_{\lambda}+\delta_{-\lambda}}_{\boM(\R)}=2.
\qedhere
\end{equation*}
\end{proof}

We can now complete the proofs of the applications of our results for the potentials discussed in the introduction.

\begin{proof}[Proof of Corollary~\ref{cor:example:gaussian}]
By Lemma~\ref{lem:Wlam}, \ref{H0},  \ref{W:derivative}, \ref{Wlambda:aeconvergence} and \ref{Wlambda:finitesup} are satisfied.
It is immediate to check that $\big(\wh\W_\lambda\big)'(\xi)\geq-m_\lambda\xi$, for all $\xi>0$, with  
$m_\lambda=2\lambda^2$.
In particular, integrating this inequality, one derives
\begin{equation}\label{proof:gaussian:kappa}  
\wh\W_\lambda(\xi)\geq 1-\lambda^2 \xi^2,\quad \text{ for all }\xi\in\R, \lambda >0,
\end{equation}
so \ref{W:kappa} is satisfied. In addition, if we consider the interval $\Lambda=(0,1/\sqrt 2)$, then we take $\kappa=1/2$ in \ref{W:kappa}, for all  $\lambda \in \Lambda$.
Thus, the hypotheses in   Theorems~\ref{thm:existenceblack}
and \ref{thm:limittolocalblack} are fulfilled, and the conclusion of \ref{item:black} follows.

To prove \ref{item:gray}, we will apply Theorem~\ref{thm:existencegray}.
Remark that the computation above shows that  \ref{W:m} is fulfilled  for all $\lambda \in (0,1/\sqrt{2})$, and 
that  \ref{W:L1} holds for any fixed $\tau>0$, with 
\begin{equation}\label{proof:gaussian:sigma}
\sigma=\sigma_\lambda=\W_\lambda(\tau)=\frac{1}{2\lambda\sqrt{\pi}}e^{-\frac{\tau^2}{4\lambda^2}}.
\end{equation}
Hence, to check the conditions in \eqref{H:msigmatau}, we choose $\tau=\tau_\lambda=\lambda$ in \eqref{proof:gaussian:sigma}, so that 
\[
\tau_\lambda\sigma_\lambda=\frac1{2e^{1/4}\sqrt{\pi} }\quad \text{ and }\quad \frac{m_\lambda}{\tau_\lambda\sigma_\lambda}=4e^{1/4}
{\sqrt \pi }\lambda^2.
\]
Thus \eqref{H:msigmatau} holds for every $\lambda\in (0,\lambda_c)$, with $\lambda_c$ defined in  \eqref{gaussian:lambdac}, and we deduce the existence of $u_\lambda$. In addition, since
$$ \lim_{\lambda\to 0} m_\lambda M(c,\tau_\lambda,\sigma_\lambda)=0,$$
and  $\pi/\sqrt{8+2c^2}<1/\sqrt{2}$,
we can invoke  Theorem~\ref{thm:limittolocalgray} to infer that the convergence in  \eqref{eq:limittolocal} holds.

To establish  \ref{item:osc}, bearing in mind that $\wc\W_\lambda(\sigma)=e^{\lambda^2\sigma^2}$, it is easy to check that the condition \eqref{H:bumps} is equivalent to $\lambda>\tilde \lambda_c$, with $\tilde \lambda_c$ defined in \eqref{gaussian:lambdactilde}. 
Hence, the conclusion follows by using  Theorem~\ref{thm:nonmonotone} and Corollary~\ref{cor:nonmonotone}.
\end{proof}

\begin{proof}[Proof of Corollary~\ref{cor:example:exponential}]
The proof follows the same lines as in the previous corollary, except for \ref{item:osc}, which was already proven in Theorem~\ref{thm:local:potentials}.
Indeed, by differentiating $\wh \W_\lambda$, we obtain as before $m_\lambda=2\lambda^2$, so that \eqref{proof:gaussian:kappa} still holds.
Since, for any fixed $\tau>0$, we have
$
\sigma_\lambda=\W_\lambda(\tau)=e^{-\frac{\tau}{\lambda}}/({2\lambda}),
$
we can choose  $\tau_\lambda=\lambda$ so that 
\[
\tau_\lambda\sigma_\lambda=\frac{1}{2e}, \quad \frac{m_\lambda}{\tau_\lambda\sigma_\lambda}=4e\lambda^2\quad 
\text{ and }
 \lim_{\lambda\to 0} m_\lambda M(c,\tau_\lambda,\sigma_\lambda)=0.
 \]
The conclusion follows as in  Corollary~\ref{cor:example:gaussian}.
\end{proof}

\begin{proof}[Proof of Corollary~\ref{cor:example:vanderwaals}]
Item \ref{item:osc} was established in Theorem~\ref{thm:local:potentials}. As before, by Lemma~\ref{lem:Wlam}, \ref{H0},  \ref{W:derivative}, \ref{Wlambda:aeconvergence} and \ref{Wlambda:finitesup} are satisfied, 
and by differentiating $\wh \W_\lambda$, we obtain 
\bq
\label{dem:cor:van}
m_\lambda=4(-\lambda)^+/(\beta^2(\beta-2\lambda)),
\ \text{ and }\ 
\wh\W_\lambda(\xi)\geq 1-\kappa_\lambda \xi^2,\quad \text{ for all }\xi\in\R,
\eq
where $\kappa_\lambda=m_\lambda/2$. 
     Therefore,  \ref{W:kappa} is satisfied and Theorem~\ref{thm:existenceblack} implies the existence of an odd black soliton $u_\lambda$.
     If  $\lambda\in (0,\beta/2)$, we invoke Theorem~\ref{thm:limittolocalblack} to get the convergence in \eqref{eq:limittolocal}.
     The limit as $\lambda\to 0^-$ is proved analogously by considering the family 
      $\{ \W_{ -\lambda } \}_{\lambda\in (0,\lambdazero)}$, for a fixed $\lambdazero>0$. To conclude the proof of \ref{item:black}, we invoke Theorem~\ref{thm:monotonicity} to get the monotonicity of the black soliton when $\lambda>0$.
     
We will apply  Corollary~\ref{cor:existencegray} to prove \ref{item:gray}. 
Notice that \ref{W:deltamu} holds with $\mu_\lambda= -\lambda e^{-\beta \abs{x}}$. Indeed,  $\mu_\lambda^-=0$ if $\lambda<0$, and 
$\norm{\mu_\lambda^-}_{\boM(\R)}=2\lambda/\beta<1 $, if $\lambda>0$. 
In addition, in view of \eqref{dem:cor:van}, we see that $m_\lambda=0$ if $\lambda>0$. Therefore, if $\lambda>0$ then \ref{W:m} and the condition in \eqref{H:mmu} trivially hold, and the conclusion follows.

For $\lambda<0$, using \eqref{dem:cor:van}, the condition $m_\lambda\in [0,1)$ is equivalent to  $\abs{\lambda}(2-\beta^2)<\beta^3/2$.
Hence, we assume in the sequel that $\lambda$ satisfies this inequality, so that \ref{W:m} is satisfied.
In conclusion, we have
\bq
\label{dem:cor:van2}
 \norm{\mu_\lambda^+}_{\boM(\R)}=\frac{2\abs{\lambda}}{\beta}, \quad \quad \norm{\mu_\lambda^-}_{\boM(\R)}=0, 
\text{ and }
M(c,\mu_\lambda)=\frac{(\beta+2\abs{\lambda})}{\beta}\Big(1+\frac{c^2}{4}\Big),
\eq
so that the condition $m_\lambda M(c,\mu_\lambda)<1$ is equivalent to $\abs{\lambda}<\beta^3/(4+c^2)$. Therefore, 
Corollary~\ref{cor:existencegray} also provides the  existence of $u_\lambda$ for $\lambda \in (\lambda_{c,\beta},0)$.
Finally, in view of \eqref{dem:cor:van} and  \eqref{dem:cor:van2}, we can invoke  
Theorem~\ref{thm:limittolocalgray}-(ii) to conclude the nonlocal-to-local limit as $\lambda\to 0$  
(arguing  for $\W_{-\lambda}$ if $\lambda$<0).
\end{proof}

\begin{proof}[Proof of Corollary \ref{cor:ex:TFI}]
By Lemma~\ref{lem:Wlam}, \ref{H0} and  \ref{Wlambda:finitesup} are satisfied.
 For the potential \eqref{pot:nematic}, we have,  for all $\abs{\lambda}\leq 1$,
$$\partial_\lambda \wh \boW_\lambda=-\frac{2\lambda\xi^2}{(1+\lambda^2 \xi^2)^2}, \quad \text {so that }\quad \abs{\partial_\lambda \wh \boW_\lambda}\leq 2\abs{\xi}^2.$$
Thus \ref{Wlambda:C1} holds. Notice that in this case $\norm{\partial_\lambda \wh \boW_\lambda}_{L^\infty(\R)}=1/(2\abs{\lambda})$, so we need to allow the quadratic growth in \ref{Wlambda:C1}  to have a uniform bound for $\lambda$ close to $0$.
It is straightforward to verify that \ref{Wlambda:C1} is also fulfilled for the other potentials. Therefore, the conclusion follows from Theorem~\ref{thm:even}.
\end{proof}

\subsubsection*{Declarations} 
{\bf Funding}. Partial financial support was received from the Labex CEMPI (ANR-11-LABX-0007-01).
S.~L\'opez-Mart\'inez was also supported by the Grant PID2021-122122NB-I00 funded by MCIN/AEI/ 10.13039/501100011033 and by ``ERDF A way of making Europe''.

\noindent{\bf Ethics approval and Conflict of interest.}  
The authors confirm that they have followed the ethical guidelines of the journal, 
that all the results they obtained are new, and that they have no competing interests to declare that are relevant to the content of this article.

\noindent{\bf Data availability} No data was used for the research described in the article.

\bibliographystyle{abbrv}

\end{document}